\RequirePackage{fix-cm}

\documentclass[preprint,12pt]{elsarticle}

\usepackage{epsfig}
\usepackage{fancybox}
\usepackage{multicol}
\usepackage{multirow}
\usepackage{color}
\usepackage{graphicx}
\usepackage{url}
\usepackage{floatrow}
\usepackage{subfig}

\usepackage{amsmath}
\usepackage{amssymb}

\usepackage{amsthm}
\usepackage{mathrsfs} 
\usepackage{empheq}

\newtheorem{hypo}{Hypothesis} 

\usepackage{lineno}

\makeatletter
\newcommand{\Rmnum}[1]{\expandafter\@slowromancap\romannumeral #1@}
\makeatother

\newtheorem{mydef}{Definition}

\usepackage{multirow}

\usepackage{comment}
\usepackage[]{algorithm2e}

\begin{document}

\begin{frontmatter}

\title{PDE-induced connection of moving frames for the Atlas of the cardiac electric propagation on 2D atrium}

\author{Sehun Chun}
\ead{sehun.chun@yonsei.ac.kr}
\address{Underwood International College, Yonsei University, South Korea}

\author{Chris Cantwell}
\ead{c.cantwell@imperial.ac.uk}
\address{Imperial College London, London, United Kingdom}

\begin{abstract}
As another critical implementation of moving frames for partial differential equations, this paper proposes a novel numerical scheme by aligning one of three orthogonal unit vectors at each grid point along the direction of a wave propagation to construct an organized set of frames, called a connection. This connection characterizes the geometry of wave propagation depending on (1) the initial point, (2) type of wave, and (3) shape of the domain with conduction properties. The constructed connection is differentiated again to derive the Riemann curvature tensor of orthonormal bases corresponding to important physical and biological meanings in wave propagation. As a practical application, the proposed scheme is applied to diffusion-reaction equations to obtain the Atlas, or a geometric map with connections, of an atrium with cardiac fibers, for the quantitative and qualitative analysis of cardiac action potential propagation, which could contribute to the clinical and surgical planning of atrial fibrillation.
\end{abstract}

\begin{keyword}
Diffusion-reaction equations \sep Moving frames \sep Connection form \sep Riemann curvature tensor \sep Cardiac electric propagation \sep Atrial fibrillation
\end{keyword}

\end{frontmatter}

\section{Introduction}

In the previous works on adapting moving frames to obrain the numerical solution of PDEs on curved surfaces such as conservational laws \cite{MMF1}, diffusion equation \cite{MMF2}, shallow water equations \cite{MMF3}, and Maxwell's equations \cite{MMF4}, the direction of the moving frames is given at random on the tangent plane in the absence of anisotropy. Because the first moving frames, denoted as $\mathbf{e}^1$, should be differentiable within each element, the direction of moving frames cannot be independent of the direction of neighboring moving frames to ensure differentiability. However, their average direction is still random. Let us call this type of distribution as \textit{moving frames of order zero} according to the Cartan's original definitions \cite{Cartan4}. In this case, moving frames are nothing more than an \textit{orthonormal local reference} reflecting the geometry of the tangent plane.

The \textit{moving frames of first order} are naturally obtained by aligning $\mathbf{e}^1$ in the direction of the tangent vector along the trajectory of a wave. In Maxwell's equations, for example, the first moving frame $\mathbf{e}^1$ is in the same direction as the electric or magnetic field, and in the diffusion-reaction equation, it is in the same direction as the gradient of the variable that represents an electric potential or material concentration. The \textit{moving frames of second order} and \textit{moving frames of third order} are obtained by additionally aligning $\mathbf{e}^2$ and $\mathbf{e}^3$ in the binormal and normal directions, respectively. In many cases, the first-order moving frames are naturally the second and third-order moving frames. Thus, we may refer to such moving frames of non-zero order collectively as the \textit{Darboux frame} \cite{HCartan}. 

\begin{figure}[ht]
\centering
\subfloat[Moving frames of order zero] {\label{MF0} \includegraphics[ width=5cm]{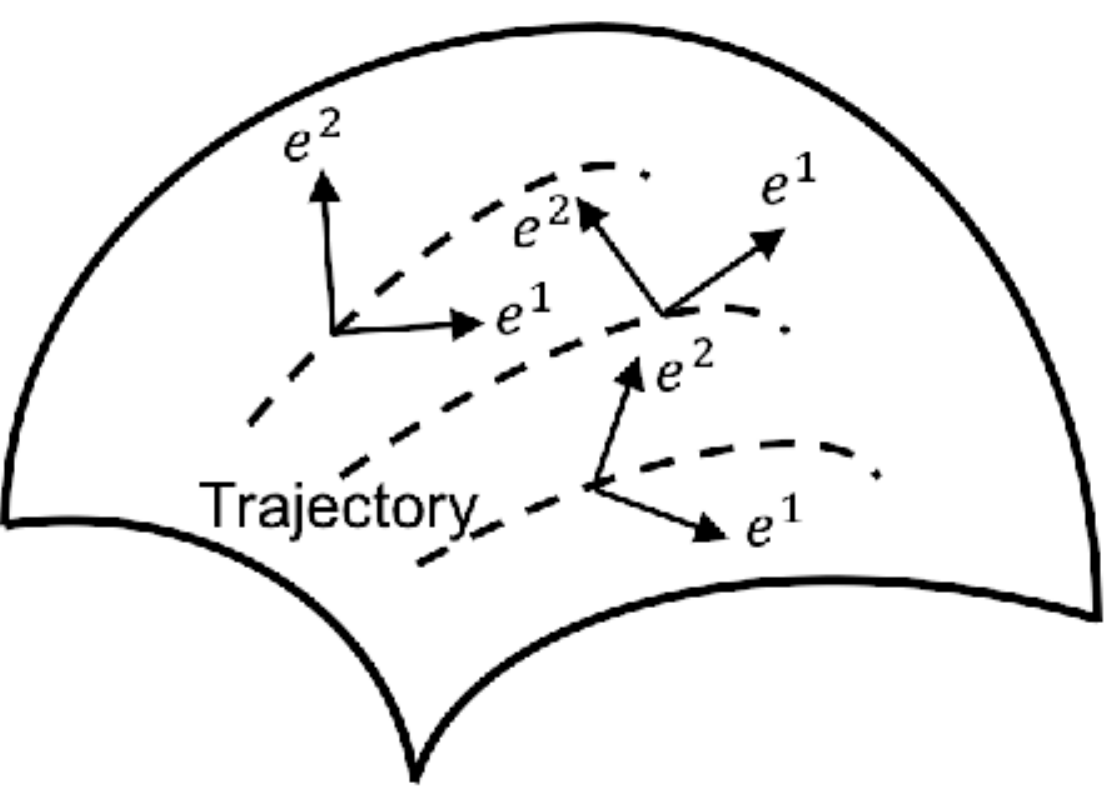} }  
\subfloat[Moving frames of first order] {\label{MF1} \includegraphics[ width=5cm]{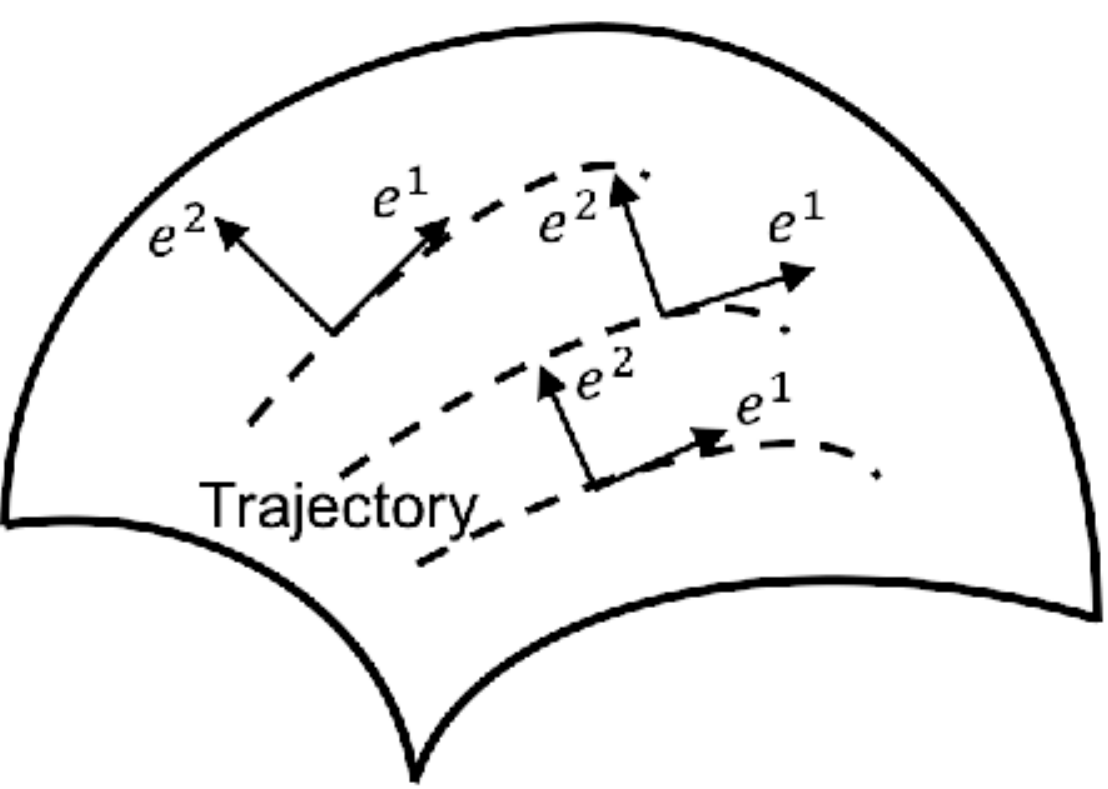} }  
\caption{Illustrations of moving frames of different kinds. The dashed line is the trajectory. $\mathbf{e}^1$ and $\mathbf{e}^2$ are the first and second moving frames, respectively.}
\label {ActZone}
\end{figure}

\subsection{Riemann curvature tensor of wave propagation}

The main objective of aligning the first moving frame in the direction of propagation is to analyze the trajectory of wave propagation qualitatively. This can be achieved by computing the Riemann curvature tensor of trajectory. This curvature is related to the \textit{shape} of the domain, but the shape is not the sole factor of the Riemann curvature. For example, the plane is \textit{flat} with zero curvature. However, the trajectory of wave propagation can have non-zero curvature; a planar propagation with flat wavefront has zero curvature because the trajectory is equally spaced and parallel. However, a point-initiated wave propagation has a trajectory that is aligned along the radial direction in the polar coordinate axis. Thus, the curvature of the trajectory is the same as the polar coordinate axis with non-zero curvature. The trajectory can have various curvature depending on many factors, particularly (1) initiation type and location, (2) type of wave, (3) shape of the domain, and (4) property of the medium.

\begin{figure}[ht]
\centering
\includegraphics[height = 4cm, width=8cm]{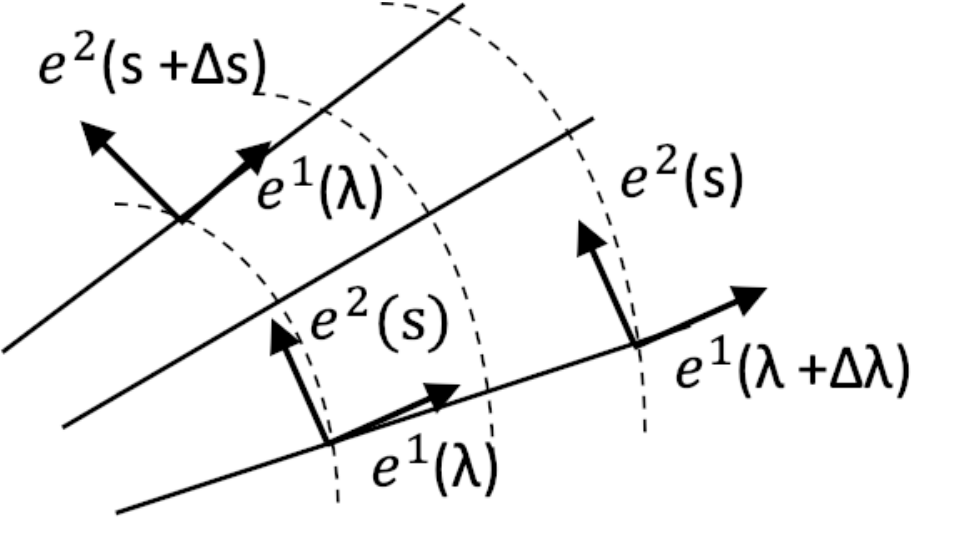} 
\caption{Illustrations of the separation vector of trajectory and its parameter.}
\label {RiemannCurvature}
\end{figure}

The Riemann curvature tensor has twenty-one components in three dimensions to indicate how moving frames change in each direction. One of the most critical components in this tensor is the component of $\mathbf{e}^2$ in the direction of the wavefront when $\mathbf{e}^1$ is aligned along the propagational direction, as shown in Fig. \ref{RiemannCurvature}. If $\mathbf{e}^2$ is parameterized with $s$, then $\partial \mathbf{e}^2 / \partial s^2$ indicates the acceleration of $\mathbf{e}^2$ along the wavefront. Note that $\mathbf{e}^2$ is the normalization of the \textit{separation vector} \cite{Misner}; thus, $\partial \mathbf{e}^2 / \partial s^2$ corresponds to the degree that the trajectory diverges or converges depending on its sign and magnitude.

A significant magnitude of the Riemann curvature component always has critical meaning in physics and biology. For example, in spacetime physics, the non-trivial magnitude of this curvature component for electromagnetic wave propagation can suggest the presence of \textit{mass} \cite{Misner, Dray}. In the heart, the non-trivial magnitude and the positive of this curvature component may indicate the \textit{stopping condition for the biological electric flow}, as suggested in ref. \cite{MyBIOP}. The objective of creating an \textit{Atlas} of the geometry map of the electric current flow on the atrium is to qualitatively and quantitatively analyze the stopping condition of the cardiac electric propagation on a patient-specific atrium structure to reveal a unidirectional block in order to cause atrial reentry and consequently fibrillation. Moreover, the Atlas may explain the role of the cardiac fiber in the heart, which has been unknown to the cardiology community and may promote more efficient surgical planning.

\subsection{Motivations for cardiac electric propagation analysis}
This paper describes the general numerical scheme to derive the Atlas of wave propagation with the subsequent connection form and Riemann curvature tensor of orthonormal bases while also providing an exemplary use for the analysis of the cardiac electric signal propagation. The motivation of this practical application will be described briefly using relevant literature in this subsection.

Cardiac electric signal propagation is often identified by a series of wavefronts \textit{in silico} or \textit{in vivo} observations. Then, the analysis of the stopping conditions or conditions for spiral wave or fibrillation is conducted by studying the shape of the wavefront, otherwise known as kinematics analysis \cite{Mikhailov, Zykov}; alternatively, spatiotemporal analysis by phase mapping after introduced in the 1980s by Winfree, can be used \cite{Winfree, Gray, Clayton}. However, both methods seem to fail to reflect the multidimensional anisotropic features of cardiac electric signal propagation in the atrium and ventricles completely. This is because the kinematic analysis cannot incorporate both the three-dimensional shape and anisotropy into the analysis \cite{MyBIOP}, whereas the phase mapping analysis has other limitations \cite{Umapathy}.

The proposed novel scheme uses the trajectory, instead of the wavefront, to analyze the spatial-temporal propagational pattern of cardiac electric signal propagation. One advantage of this approach is that the geometrical factors, such as the three-dimensional shape, direction, and strength of the cardiac fiber, can be conveniently incorporated into the Atlas in the propagation. Thus, the direct causal relation between the geometry and the propagational pattern can be more easily identified. Another advantage is that, similar to space-time physics, all different types of cardiac physiological phenomena that affect the electric signal propagation, such as cardiac restitution, cardiac memory, or the presence of infarction or fibrosis, can be incorporated into this geometric model to be interpreted considering geometric deformation.

\subsection{Goals}
This paper may provide a more significant result that can be immediately used for the cardiac electric signal propagation in the complex multidimensional and anisotropic heart. However, as the first step, this paper focuses on the description of the numerical scheme, validations, and exemplary use of this numerical scheme for cardiac electrophysiological analysis. The goals of this paper are as follows: (1) to develop an algorithm to efficiently align the first moving frames to the gradient of a variable in order to achieve the Darboux frames; (2) to compute the connection form and Riemann curvature tensor using Cartan's structure equation for diverse propagation on two-dimensional curved surfaces; (3) to apply the scheme to an actual atrium with cardiac fiber to understand the role of cardiac fiber in the propagation of the cardiac action potential.

The remainder of this paper is organized as follows. Chapter 2 describes the formulation to compute the connection form, and Chapter 3 explains the importance of a component of the connection form, $\omega_{212}$, with respect to the Riemann curvature tensor. Chapter 4 elaborates this importance in cardiac electric propagation. Chapter 5 describes the numerical scheme to derive the propagational direction from a diffusion-reaction equation. In Chapter 6, the numerical scheme to solve a diffusion-reaction equation is briefly described. Chapters 7 and 8 provide the exponential convergence of moving frames, connection form, and Riemann curvature for the test in the plane and on the sphere, respectively. Chapter 9 illustrates the meaning of anisotropy in trajectory tracking for cardiac electric propagation. Chapter 10 describes the pre-processing of raw fiber data. In Chapter 11, the two-dimensional cardiac fiber is analyzed to provide a geometric meaning to the propagation. In Chapter 12, the connection from the numerical simulation of the cardiac electric propagation is presented and explained for future clinical applications. The discussion is provided in Chapter 13.

\section{Connection Form}

Consider a set of orthonormal vectors, known as \textit{moving frames} $\mathbf{e}^i,~1 \le i \le 3$ at every grid point in a curved surface $\Omega$. Let two frames $\mathbf{e}^1$ and $\mathbf{e}^2$ lie on the tangent plane of the surface and $\mathbf{e}^3$ be aligned along the surface normal direction. Because moving frames are orthonormal, they satisfy the following conditions
\begin{equation*}
\mathbf{e}^i \cdot \mathbf{e}^j = \delta_{ij},~~~~ \|  \mathbf{e}^i \| = 1, ~~~~~ 1 \le i,j \le 3
\end{equation*}
where $\delta_{i j}$ is the Kronecker delta. The only freedom that the moving frames have is with respect to the direction of the first moving frame $\mathbf{e}^1$ on the tangent plane. The directions of other frames are subsequently determined because of the orthogonality of the frames. However, refs. \cite{MMF1} and \cite{MMF2} proved that the direction of $\mathbf{e}^1$ can affect the integration and differentiation error in differential operators such as divergence and Laplacian operators.

In this paper, $\mathbf{e}^1$ is aligned such that the connectivity between the moving frames is considered; this is termed as a \textit{specifically connected and ordered set of frames}. We may refer to this configuration of frames in the domain as a \textit{connection}. There are infinitely many ways to align moving frames. However, we choose $\mathbf{e}^1$ to be aligned along the direction of the wave propagation, which is the same as the direction of a vector or the gradient of a variable. Then, the pattern of the \textit{trajectory}, i.e., the convergence or divergence of the bundle of trajectories, can be expressed by the \textit{connection form}.

At every point in the domain, the set of moving frames is expressed in the following matrix form: for unit vectors $(\mathbf{x}, \mathbf{y}, \mathbf{z})$ along the Cartesian coordinate axes $(x,y,z)$, respectively, we obtain
\begin{align}
\left [
\begin{array}{c}
\mathbf{e}^1 \\
\mathbf{e}^2 \\
\mathbf{e}^3 
\end{array}
\right ] = \left [
\begin{array}{ccc}
{e}^1_x & {e}^1_y & {e}^1_z \\
{e}^2_x & {e}^2_y & {e}^2_z \\
{e}^3_x & {e}^3_y & {e}^3_z 
\end{array}
\right ] \left [
\begin{array}{c}
\mathbf{x} \\
\mathbf{y} \\
\mathbf{z} 
\end{array}
\right ] , 
\end{align}
Or,
\begin{align}
 \widehat{\mathbf{e}} =  \mathbf{A} \widehat{\mathbf{x}},   \label{mf1}
\end{align}
where we introduced a new tensor $\widehat{\mathbf{e}} = \left [ \mathbf{e}^1,~\mathbf{e}^2,~ \mathbf{e}^3  \right ]^T$, $\widehat{\mathbf{x}} = \left [\mathbf{x},~\mathbf{y}, ~\mathbf{z} \right ]^T $. The matrix $ \mathbf{A}$ represents the orientation of moving frames; thus, it is known as the \textit{orientation matrix} or \textit{attitude matrix} \cite{ONeil}. By using the fact that $\widehat{\mathbf{x}}$ is fixed everywhere and by computing the infinitesimal displacement of $\widehat{\mathbf{e}}$ on a surface, the 1-form $d \widehat{\mathbf{e} }$ is expressed as follows.
\begin{equation*}
d \widehat{\mathbf{e} } = \mathcal{W} \widehat{\mathbf{e}} ,~~~~~\mbox{where}~~~ \mathcal{W} \equiv  ( dA  ) A^T
\end{equation*}
where we used Eq. \eqref{mf1} and the fact that $A$ is an orthonormal matrix. Briefly speaking, the 1-form is the general type of derivative such that it has a scalar value when a vector is chosen. For example, let $v$ be a scalar; then, the exterior derivative of $v$, or $dv$, is the 1-form. Consequently, $dv \langle \mathbf{v} \rangle $ is a scalar value corresponding to the directional derivative of $v$ in the direction of $\mathbf{v}$, or $\nabla v \cdot \mathbf{v}$ \cite{Misner}.

Then, $\mathcal{W}$ is a 1-form with nine components $[w^i_j]~~1 \le i,j \le 3$ and is called the \textit{connection form} \cite{Cartan1, Cartan2}. In words, $w^i_j \langle \mathbf{v} \rangle$ represents the amount of rotation of $\mathbf{e}^i$ with respect to $\mathbf{e}^j$ when it moves in the direction of $\mathbf{v}$ \cite{Piuze2015}.  However, the orthonormality of moving frames yields the following equalities.
\begin{align}
&d \mathbf{e}^i \cdot \mathbf{e}^i = 0, \label{mf21}\\
&d \mathbf{e}^i \cdot \mathbf{e}^j + d \mathbf{e}^j \cdot \mathbf{e}^i = 0 .  \label{mf22}
\end{align}
Eq. \eqref{mf21} yields $\omega^i_{i}= 0$ for $1 \le i \le 3$ and Eq. \eqref{mf22} yields $\omega^i_j + \omega^j_i=0$. Thus, $\mathcal{W}$ only has three independent components, namely
\begin{equation*}
\mathcal{W}= \left [
\begin{array}{ccc}
0 & \omega^2_1 & \omega^3_1 \\
- \omega^2_1 & 0 &  \omega^3_2 \\
- \omega^3_1 & -\omega^3_2 & 0 
\end{array}
\right ] , 
\end{equation*}
Or, we have
\begin{align*}
d \mathbf{e}^1 &= ~~~~~~~~   \omega^2_1 \mathbf{e}^2 + \omega^3_1 \mathbf{e}^3, \\
d \mathbf{e}^2 &=  -\omega^2_1 \mathbf{e}^1 +~~~~~~   \omega^3_2 \mathbf{e}^3, \\
d \mathbf{e}^3 &=  -\omega^3_1 \mathbf{e}^1 - \omega^3_2 \mathbf{e}^2  . ~~~~~ ~~~
\end{align*}
Additionally, 1-form is linear such that $\omega^i_j \langle \mathbf{v} \rangle = v^1 \omega^i_j \langle \mathbf{e}^1 \rangle + v^2 \omega^i_j \langle \mathbf{e}^2 \rangle$, it would be sufficient to compute $\omega^i_j \langle \mathbf{e}^k \rangle$ for $1 \le k \le 2$. For simplicity, the following notation is used for connection form such that the magnitude of $\omega_{ij}$ is the same as that of $\omega^i_j$, but possibly with opposite signs, and $\omega_{ijk} \equiv \omega^i_j \langle \mathbf{e}^k \rangle$. Then, the connection form $\mathcal{W}$ can be obtained as
\begin{equation*}
\mathcal{W} \langle \mathbf{e}^k \rangle = d A \langle \mathbf{e}^k \rangle A^T,
\end{equation*}
Thus,
\begin{equation}
\omega_{ijk} = [e^i_x ~ e^i_y~ e^i_z ] \left [
\begin{array}{ccc}
\dfrac{d e^j_x}{dx} & \dfrac{d e^j_x }{d y}  & \dfrac{d e^j_x }{d z}  \\
\dfrac{d e^j_y}{dx} & \dfrac{d e^j_y}{d y}  & \dfrac{d e^j_y}{d z}  \\
\dfrac{d e^j_z}{dx} & \dfrac{d e^j_z }{d y}  & \dfrac{d e^j_z}{d z}  \\
\end{array}
\right ]  \left [
\begin{array}{ccc}
e^k_x \\
e^k_y \\
e^k_z
\end{array}
\right ]  ,  \label{omegaijk}
\end{equation}
where the $3 \times 3$ matrix is the Jacobian matrix for $\mathbf{e}^j$. There are nine components of the connection form, three for each of $\omega_{1}^2$, $\omega_{1}^3$, and $\omega_{2}^3$. Out of the three, $\omega_{1}^3$ and $\omega_{2}^3$ are related to the Gaussian curvature ($\mathcal{K}$) and mean curvature ($\mathcal{H}$) because they can be also expressed in terms of connections as follows \cite{Treibergs}
\begin{align*}
\mathcal{K} &= \omega_{311} \omega_{322} - \omega_{312} \omega_{321},   \\
\mathcal{H} &= \frac{1}{2} ( \omega_{311} + \omega_{322} )   .
\end{align*}
In a two-dimensional plane, $\omega^3_1$ and $\omega^3_2$ are zero along all the directions. However, various connections of moving frames can create a non-zero $\omega^2_1$ along a certain direction. In this paper, we are particularly interested in the 1-form $\omega^2_1$ along the direction of the second frame $\mathbf{e}^2$ because the corresponding quantity $ \omega_{212}$ is related to an important component of the Riemann curvature tensor as will be explained in the following section.

\section{Importance of $\omega_{212}$}
One of the benefits of deriving the connection form of moving frames is obtaining a Riemann curvature tensor representing important features of a physical or biological phenomena. However, computing the Riemann curvature tensor requires the metric tensor of the curved domain, or the length of the tangent vector of the curved axis of the surface. We refer to this tensor as the \textit{Riemann curvature tensor of coordinate basis} (the \textit{Riemann curvature tensor of holonomic basis}). However, computing this Riemann curvature tensor of coordinate basis is computationally challenging for the general curved surface, and most of all, is not always required.

Instead, the length of the tangent vector of the trajectory is normalized and the Riemann curvature tensor is computed for the normalized trajectory. This yield the \textit{Riemann curvature tensor of orthonormal basis}. Geometrically, this Riemann curvature describes the same qualitative behavior for the divergence and convergence of the trajectory in the orthogonal direction of propagation, or in the direction of the separation vector, but its magnitude may be changed proportionally. Note that the Riemann curvature tensor is exactly same as this Riemann curvature tensor of orthonormal basis if the metric tensor is constant.

To derive the Riemann curvature tensor of orthonormal basis, the second Cartan structure equation is recalled as follows \cite{Cartan1, Cartan2}.
\begin{equation*}
\Omega^j_i = d \omega^j_i  + \omega^j_k \wedge \omega^k_i  ,
\end{equation*}
where $d \omega^j_i$ is a 2-form, $\wedge$ is the wedge product, and $\Omega^j_i$ is known as the \textit{curvature 2-form} that is defined by $d^2 \mathbf{e}_i= \Omega^j_i \mathbf{e}_j$ \cite{Dray}. The \textit{Riemann curvature tensor} $\overline{\mathscr{R}}^j_{imn}$ is defined in relation with the curvature 2-form as follows.
\begin{equation*}
\Omega^j_i = \overline{\mathscr{R}}^j_{imn} \omega^m \wedge \omega^n,
\end{equation*}
where the 1-form $\omega^m$ representing an \textit{orthonormal dual frame} satisfies $\omega^i ( \mathbf{e}^j ) = \delta^i_j$ where $\delta^i_j$ is the Kronecker delta. The \textit{overline} notation is introduced to indicate the corresponding Riemann curvature computed with orthonormal bases. For the connection component $\omega^2_1$, $\overline{\mathscr{R}}^j_{imn}$ can be computed in the following manner.
\begin{equation}
\overline{\mathscr{R}}^2_{1mn}  \equiv \Omega^2_1 \langle  \mathbf{e}^m, \mathbf{e}^{n} \rangle = d \omega^2_1  \langle  \mathbf{e}^m, \mathbf{e}^{n} \rangle + \sum_{k=1}^2 \omega^2_k \wedge \omega^k_1 \langle  \mathbf{e}^m, \mathbf{e}^{n} \rangle,   \label{Riemannmain}
\end{equation}
where we used the property $\omega^j_i = - \omega^i_j$ and $\omega^i_i = 0$. The numerical evaluation of $d \omega^j_i  \langle  \mathbf{e}^m, \mathbf{e}^{n} \rangle$ can be obtained as follows. Consider $\omega^j_i$ expressed on a two-dimensional surface as follows.
\begin{equation*}
\omega^j_i = \omega^j_i \langle \mathbf{e}^1 \rangle \omega^1 + \omega^j_i \langle \mathbf{e}^2 \rangle \omega^2  .
\end{equation*}
Then, the derivative of this 1-form is obtained in the form of directional derivative by using the variable $\omega_{ijk}$ from Eq. \eqref{omegaijk} as follows.
\begin{equation}
d \omega^j_i \langle \mathbf{e}^1, \mathbf{e}^2 \rangle =   \nabla \omega_{ji1} \cdot  \mathbf{e}^2  -  \nabla \omega_{ji2} \cdot  \mathbf{e}^1      \label{Compute2form}
\end{equation}
The exact calculation of $d \omega^j_i$ should be obtained based on $d \xi_i$, the axis aligned along the same direction as $\mathbf{e}^i$ but with the non-unit magnitude of the tangent vector. Here, $d \omega^j_i$ is differentiated based on the unit axis rather than the underlying curved axis. This yields the Riemann curvature tensor of orthonormal bases, which differs from that of coordinate bases.

The second component in the right hand side of Eq. \eqref{Riemannmain} is known as the \textit{Gaussian torsion} \cite{Cartan1} and can be computed as follows.
\begin{equation}
\omega^j_k \wedge \omega^k_i \langle  \mathbf{e}^m, \mathbf{e}^{n} \rangle = \det
\left |
\begin{array}{cc}
\omega^j_k \langle  \mathbf{e}^m \rangle & \omega^j_k \langle  \mathbf{e}^n \rangle  \\
 \omega^k_i \langle  \mathbf{e}^m \rangle &    \omega^k_i \langle  \mathbf{e}^n \rangle   
\end{array}
\right |  \label{Computewedge}
\end{equation}
However, there exists a property that $\omega^i_i = 0$ reduces this Gaussian torsion to null. Substituting Eqs. \eqref{Compute2form} and \eqref{Computewedge} into Eq. \eqref{Riemannmain}, the following formula are obtained for the nontrivial components of the Riemann curvature tensor of orthonormal basis. 
\begin{equation}
\overline{\mathscr{R}}^2_{121}  =   \nabla \omega_{211} \cdot \mathbf{e}^2 -  \nabla \omega_{212} \cdot \mathbf{e}^1 .  \label{Riemanfinal21}
\end{equation}
Note that $\omega_{211}$ is negligible because the rotation of $\mathbf{e}^2$ with respect to $\mathbf{e}^1$ when moving toward $\mathbf{e}^1$ remains fixed orthogonally. Then, $\overline{\mathscr{R}}^2_{121}$ is approximately the same as the change of $\omega_{212}$ in the direction of $\mathbf{e}^1$ on a surface.

\begin{figure}[ht]
\centering
\subfloat[$D=0$, $\overline{\mathscr{R}}^2_{121} = 0$ ] {\label{R1} \includegraphics[height = 3cm,
 width=4cm]{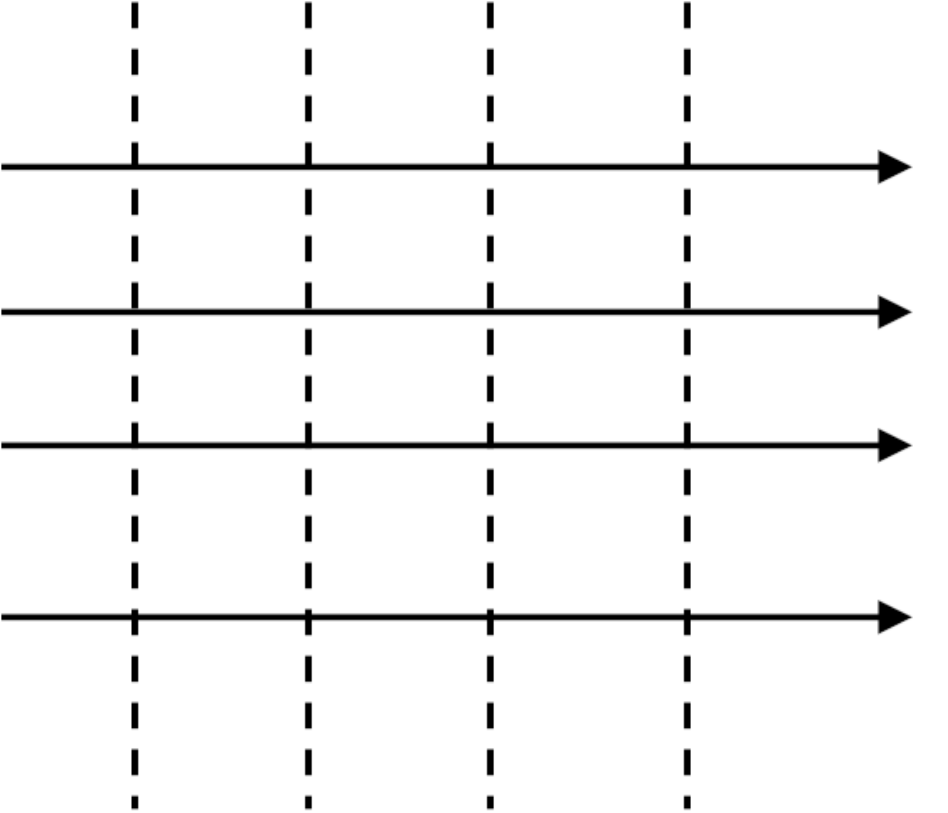} }  
 \subfloat[$D<0$. $\overline{\mathscr{R}}^2_{121} <0$ ] {\label{R1} \includegraphics[height = 3cm,
 width=4cm]{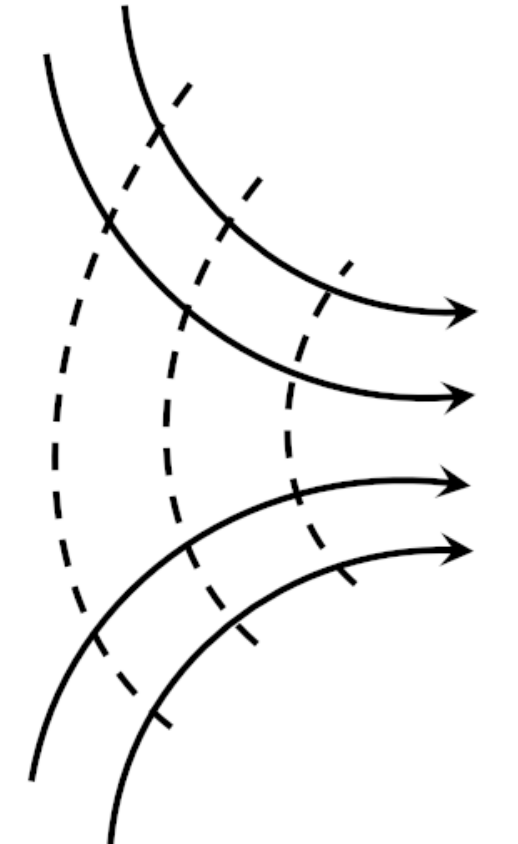} }  
 \subfloat[$D>0$. $\overline{\mathscr{R}}^2_{121} >0$ ] {\label{R1} \includegraphics[height = 3cm,
 width=4cm]{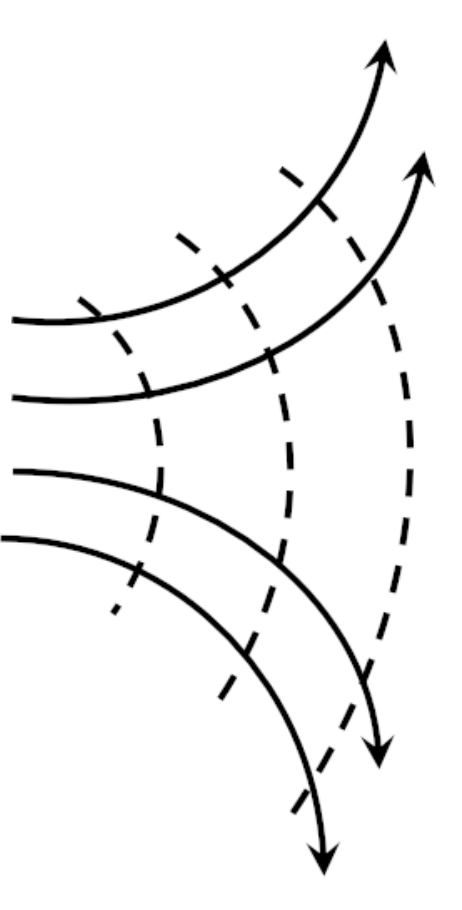} }  
\caption{Intuitive relationship between wavefront curvature and Riemann curvature. Solid line = trajectory, Dashed line = wavefront.}
\label {DandR}
\end{figure}

\section{Stopping criteria in Riemann curvature}
As a large Riemann curvature implies the presence of an object with a large mass in the universe such as a black hole, a large Riemann curvature in the cardiac excitation propagation corresponds to a similar unique feature in the cardiac tissue. In the cardiac excitation propagation, one of the most important analyses of the propagation is related to whether it continues to propagate or stops. The stopping conditions are the most critical knowledge for understanding many cardiac pathologies. However, reliable analytical or computational methodologies remain unknown considering the strong anisotropy and complex shape of the heart. Initial attempts to identify the general stopping conditions were related to the shape of the wavefront, known as the \textit{kinematic analysis}; the linearized version of the propagation velocity ($V$) can be obtained as \cite{Davydov1991}
\begin{equation}
V \approx V_0 - D K ,  \label{Davy}
\end{equation}
where $V_0$ is the conducting velocity determined by the cardiac tissue, and $D$ is the diffusivity constant, and $K$ is the curvature of the wavefront. Eq. \eqref{Davy} implies that the propagation velocity of the cardiac propagation depends on the curvature of the wavefront. If the conducting properties of cardiac tissue are homogeneous, then the only important factor for slowing down the propagation is the shape of the wavefront. If the curvature is too large or above a certain threshold that is determined by the conducting properties of the cardiac tissue, then the propagation velocity dramatically reduces, even to zero.

\begin{figure}[ht]
\centering
\subfloat[ Moving frames ] {\label{Gap10MF} \includegraphics[
 width=4cm]{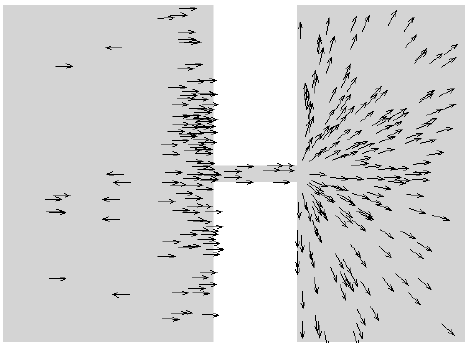} }  
\subfloat[$w_{212}$] {\label{Gap10w122} \includegraphics[
width=4cm]{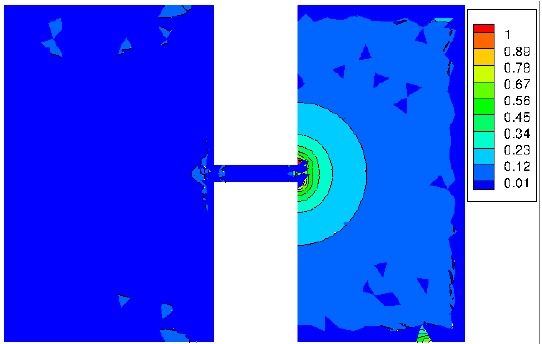} }
\subfloat[$\overline{\mathscr{R}}^2_{121}$ ] {\label{Gap10R1221} \includegraphics[
width=4cm]{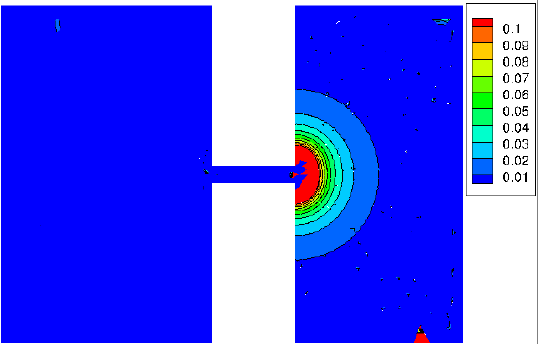} }
\caption{Propagation of the cardiac electric propagation through a small gap.}
\label {Gap10}
\end{figure}

On the cellular scale, the large curvature of the wavefront is related to a phenomenon called \textit{impedence mismatch} or \textit{sink-source mismatch} \cite{Rohr, Xie}, which is a higher ratio between the number of \textit{excited} cells and the number of \textit{excitable} cells than the threshold ratio at which propagation can occur. In this paper, it is assumed that this sink-source mismatch can also be measured by the changes in the connection form or the \textit{Riemann curvature tensor of orthonormal basis} as follows.
\begin{hypo}
If the connection component $w_{212}$ changes abruptly in the direction of $\mathbf{e}^1$, or equivalently the positive $\overline{\mathscr{R}}^2_{121}$ is sufficiently large, then the cardiac action potential propagation slows down dramatically, even to zero.
\end{hypo}
The rigorous validation of this hypothesis is beyond the scope of this paper, but the complementary relationship between the wavefront and trajectory can be intuitively given in Fig. \ref{DandR}. Moreover, Fig. \ref{Gap10} demonstrates the corresponding large quantity of $\omega_{212}$ and $\overline{\mathscr{R}}^2_{121}$ for the classical gap propagation test case to demonstrate a relation between strong wavefront curvature and large positive $\omega_{212}$. From the geometric perspective, the Riemann curvature tensor measures the convergence or divergence of the geodesics of propagational rays \cite{Misner, Dray}. To measure whether trajectories are pushed closer or pulled apart, the first moving frame ($\mathbf{e}^1$) is aligned along the propagational direction and the second moving frame ($\mathbf{e}^2$) is aligned orthogonal to $\mathbf{e}^1$ and the surface normal vector. This implies that $\mathbf{e}^2$ is aligned along the wavefront of the wave, as shown in Fig. \ref{RiemannCurvature}. Therefore, the second covariant differentiation of $\mathbf{e}^2$ along $\mathbf{e}^1$ as it moves along $\mathbf{e}^2$, which is mathematically the same as $\omega_{212}$, is equivalent to the Riemann curvature of the orthonormal basis representing the convergence or divergence of trajectories such as \cite{Misner}
\begin{equation*}
\nabla_{\mathbf{e}^1} \nabla_{\mathbf{e}^1} \mathbf{e}^2 + \mathbf{e}^{j} \overline{\mathscr{R}}^j_{imn} e^1_i e^2_m e^1_n  = 0,  \label{Rieman1}
\end{equation*}
where we used the notation that $\mathbf{e}^i = e_1^i \mathbf{x} + e_2^i \mathbf{y} + e_3^i \mathbf{z}$. This is coincident with the analysis of relative acceleration that was used for the reentry around the pulmonary vein \cite{MyBIOP}.

\section{Aligning along the diffusion propagation}

\begin{figure}[ht]
\centering
\subfloat[AP model ] {\label{ET03} \includegraphics[
width=5cm]{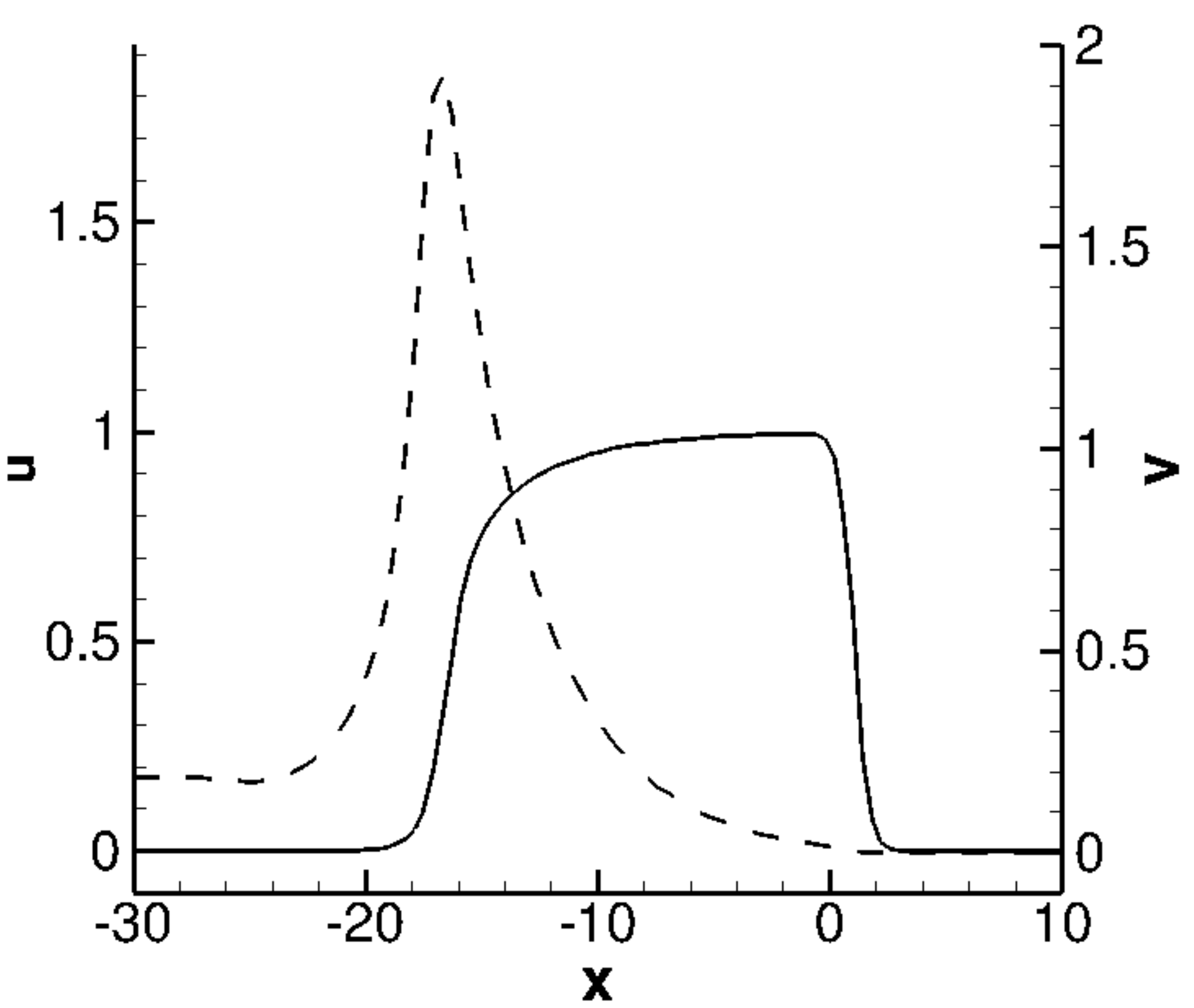} }  
\subfloat[CRN model ] {\label{ET03} \includegraphics[
width=5cm]{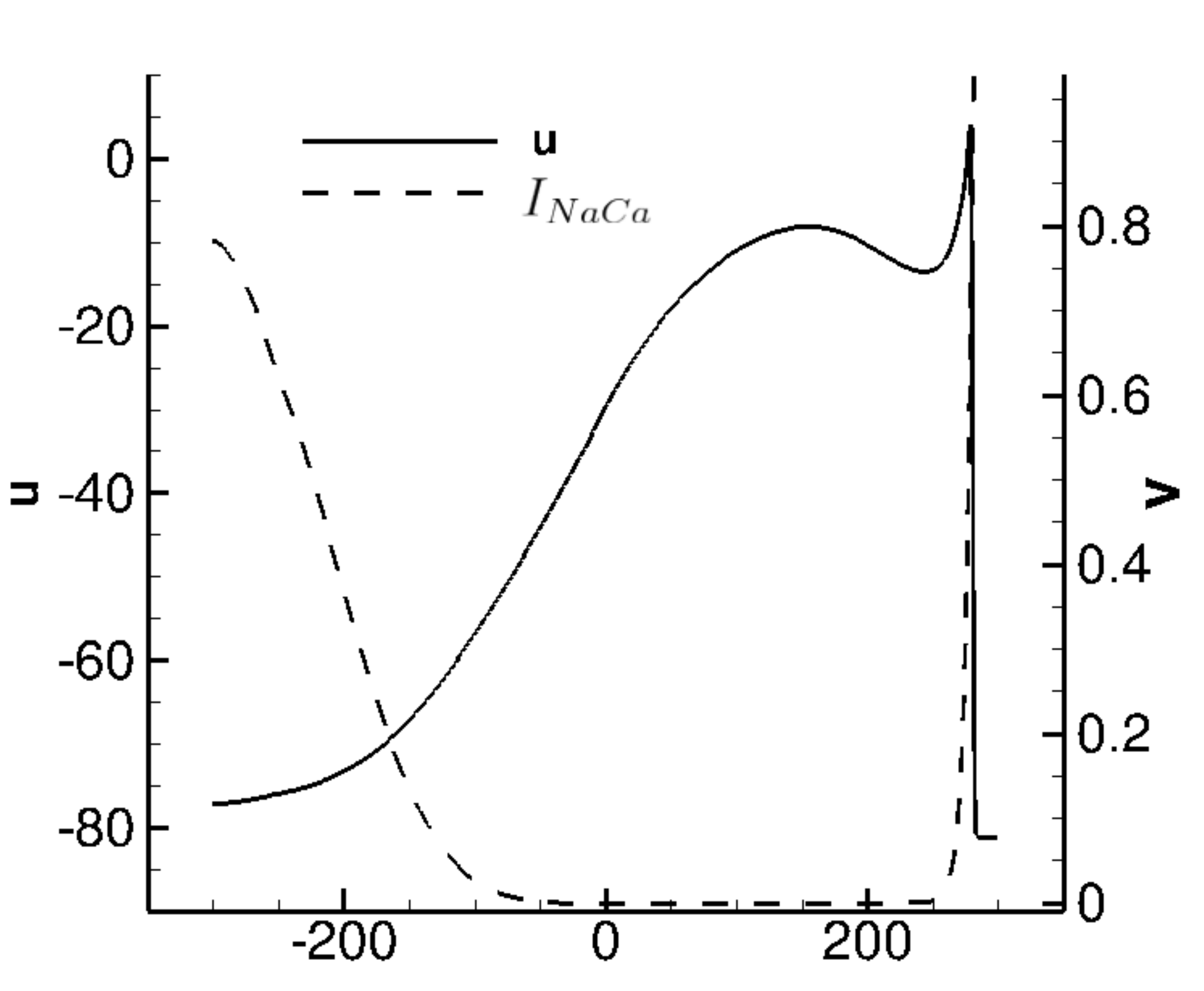} }  
\caption{AP model: $u$ is the solid line, $v$ is the dashed line. CRN model: $u$ is the solid line, $I_{NaCa}$ is the dashed line. The wave travels from left to right. }
\label{ActionPotential}
\end{figure}

The diffusion-reaction model is one popular mathematical model of the cardiac action potential propagation for the simulation and analysis of cardiac action potential propagation on an anisotropic atrium . Since the original FitzHugh-Nagumo model in the 1960's \cite{FitzHugh} as the macroscopic version of the famous Hodgkin-Huxley model \cite{Hodgkin}, several accurate models have been proposed. The two variable Aliev-Panfilov model \cite{Aliev} and twenty-one variable Courtemanche-Ramirez-Nattel (CRN) model \cite{Coutemanche} are two examples. They exhibit similar wave patterns, called the action potentials, as shown in Fig. \ref{ActionPotential}.

In the propagation, the diffusion-reaction process is often represented by the wavefront of the action potential \cite{Zykov}. The trajectory in this process is often considered the \textit{time series of particle traces in the diffusion process} \cite{Tejedor}. However, its counterpart may not correspond to the wavefront on the same scale. In this paper, the trajectory in the diffusion-reaction process is defined as follows.

\begin{mydef}
The tangent unit vector of the trajectory for a diffusion-reaction process is the same as the normalized gradient of the action potential.
\end{mydef}

In an isotropic medium, the gradient of the action potential is the same as the orthogonal to the wavefront. Thus, the above definition implies that the trajectory is orthogonal to the wavefront. However, in an anisotropic medium, the gradient of the action potential is not necessarily orthogonal to the wavefront but depends on the charge distribution of the action potential, and its gradient determines the advancement of propagation. This definition represents the exact mechanism of a diffusion-reaction process but may have additional requirements in computational modeling due to discretization errors. This is called the \textit{validity} of the direction, as will be explained in the next subsection.


 \subsection{Validity}
 \label{MMFrule1}
 For the validity of the direction of the gradient, the independence of the coordinate systems is introduced for the solution of partial differential equations. Consider the following formulation of the diffusion-reaction equation:
\begin{align}
\frac{\partial u}{\partial t} &= \nabla \widehat{\mathbf{d}} \nabla u + \mathcal{R}(u,v, \cdots) , \label{diffusion}
\end{align}
where $\widehat{\mathbf{d}}$ is the diffusivity tensor and $\mathcal{R}$ is the reaction function. Then, the constructed moving frames are the local reference frames, for which the solution obtained using one set of moving frames should not be swapped with the solution using another set of moving frames. This property is defined as the \textit{validity} of moving frames, described as follows.
 
\begin{mydef}
Let $\mathbf{e}^{orig}$ be the initially constructed moving frames that yield the stable and efficient numerical solution of Eq. \eqref{diffusion}. Then, the newly adapted moving frames $\mathbf{e}^{new}$ are called \textbf{valid} if they satisfy
\begin{equation}
\left . \left \|  \frac{\partial u (\mathbf{e}^{orig}) }{\partial t}  -  \frac{\partial  u (\mathbf{e}^{new}) }{\partial t}   \right \|   \right /  \left \|  \frac{\partial u (\mathbf{e}^{orig}) }{\partial t}  \right \|  < \delta,
\end{equation}
for a sufficiently small tolerance $\delta$ or, approximately 
\begin{equation}
\frac{ \left \|   \Delta  u (\mathbf{e}^{orig})   - \Delta u (\mathbf{e}^{new})    \right \|   }{   \left \|   \Delta  u(\mathbf{e}^{orig})  \right \| }  < \delta \Delta t, \label{ValidCond}
\end{equation}
\end{mydef}
where the tolerance $\delta$ is a function of $h$ and $p$, i.e., $\delta = \delta (h, p)$ and is independent of $\Delta t$.

The above condition implies that if the difference of $d u$ with the new moving frames $\mathbf{e}^{new}$ is sufficiently lower than the discretization error, then adapting the new moving frames does not affect the numerical accuracy and stability during time advances. Therefore, the corresponding moving frames can be adapted as a \textit{valid} connection because the direction of $\mathbf{e}^i$ should not affect the computed value of $ \nabla \widehat{\mathbf{d}} \nabla u$ as long as they are orthogonal and lie on the tangent plane. This should always be geometrically true and inscribed in the equation. On the other hand, any invalid moving frames quickly yield unstable and inaccurate computations.

In actual computations, the condition Eq. \eqref{ValidCond} discards any non-differentiable directions of the propagation within an element due to numerical noise or discretization error. For example, suppose that the direction of the gradient at a point is abruptly changed from a smooth distribution of the directions at the other points in an element. Then, the non-differentiable direction at the point is not likely constructed from the diffusion-reaction process because of its dissipative and smooth property. However, the condition Eq. \eqref{ValidCond} is also not satisfied when the curvature of the flow is large, and the spatial resolution is not sufficient in this region. This occurs because the current numerical scheme to compute the covariant derivative requires more grid points in more highly curved regions. Therefore, the parameter $\delta$ must be adjusted to include all the curvatures of the propagation while excluding the gradients of non-realistic propagation.

\subsection{Position of the gradient}
In constructing the trajectory, the unique shape of the wave should be considered. Let \textit{wavefront} be the area where the membrane potential ($u$) increases and \textit{waveback} be the area where $u$ decreases. The AP and CRN models share the generic features of action potential, such as stiffer wavefront, as shown in Fig. \eqref{ActionPotential}, but the increasing or decreasing rates of the wavefront and waveback are different. Thus, integrating the gradient over a time interval should be \textit{sign} sensitive, depending on whether the gradient is on the wavefront or waveback. The location of time integration is chosen through the following methods. 
\begin{enumerate}
 \item Wavefront where ${d u} / {dt}$ is positive. $\nabla u$ is considered as the direction of the tangent vector. (\textsf{WF})
  \item Waveback where ${d u} / {dt}$ is negative. $- \nabla u$ is considered as the direction of the tangent vector. (\textsf{WB})
\end{enumerate}

\subsection{Extracting the gradient}
Moreover, Definition 2 indicates that the trajectory of the diffusion-reaction model is constructed using the unit tangent vector, which is constructed through one of the following methods.
\begin{enumerate}
 \item To choose the unit tangent vector when the magnitude of the gradient is the largest (\textsf{Str}).
\item To integrate all the unit tangent vectors over a time interval when the magnitude of the gradient is non-trivial (\textsf{Int}).
\item To perform weighted integration over all the tangent vectors during a time interval proportional to the magnitude of the gradient (\textsf{Wint}).
\end{enumerate}
The three methods are compared, and one superior method is chosen. Combining the procedures mentioned above for accurate and stable moving frames for Eq. \eqref{diffusion}, we obtain an implementable algorithm for aligning moving frames along the propagational direction, as shown in Algorithm \ref{MMF1stGV}. 

\fbox{
 \begin{algorithm}[H]
 \KwData{Gradient vector (GV) and moving frames (MF)}
 \KwResult{Construct a connection}
  - Initialization\;
 \While{Magnitude of GV is not negligible}{
   - Align MF along GV\;
  \eIf{Newly aligned MF is derived from the following procedure \\
  (1) Validity Check (Eq.\eqref{ValidCond}), \\
  (2) Choose specific position (\textsf{WF}, \textsf{WB}), \\
  (3) Extracting the gradient (\textsf{Str}, \textsf{Int}, \textsf{Wint}) \\}{
   - Accept the newly aligned MF\;
   }{
   - Retrieve the original MF\;
  }
 }
 \caption{Align moving frames along the gradient vector}
 \label{MMF1stGV}
\end{algorithm}
}

\section{Numerical scheme for diffusion-reaction equation}
To apply the alignment scheme along the propagation direction, as explained in the previous section, a stable and accurate numerical scheme to solve diffusion-reaction equations on curved surfaces is required. The numerical scheme to solve Eqs. \eqref{AP1} and \eqref{AP2}, called the \textit{MMF-diffusion} scheme, on an anisotropic and curved surface is achieved in the context of discontinuous Galerkin methods and is implemented at an open-sourced spectral/hp library, called Nektar++ \cite{Nektar++}. The MMF-diffusion scheme is obtained by expanding the gradient $\widehat{\mathbf{d}} \nabla u$ on moving frames as follows. 
\begin{equation}
\widehat{\mathbf{d}} \nabla u = d^{11} q^1 \mathbf{e}^1 + d^{22} q^2 \mathbf{e}^2 = q^1 \mathbf{d}^1 + q^2 \mathbf{d}^2   .\label{dani}
\end{equation}
Note that the new moving frame $\mathbf{d}$ is orthogonal, but it is not necessarily of unit length depending on the diffusion coefficient, particularly in the presence of an anisotropic medium such as a cardiac fiber. Then, we have the following MMF-diffusion scheme \cite{MMF2}
\begin{align}
& \int  f \varphi dx =  \sum_{i=1}^2 \left [ \int   q_i \left ( \nabla \varphi \cdot \mathbf{d}^i \right ) dx - \int_{\partial \mathcal{M}}  \tilde{q}_i  (  \mathbf{n} \cdot \mathbf{d}^i ) \varphi ds \right ]  , \label{WeakDiffusion1} \\
& \int  d^{mm} {q}_m \varphi dx =  - \int  ( \nabla \varphi \cdot \mathbf{d}^m ) u dx - \int  ( \nabla \cdot  \mathbf{d}^m ) u  \varphi dx   + \int_{\partial \mathcal{M}^e} \tilde{u}  (  \mathbf{n} \cdot  \mathbf{d}^m ) \varphi ds .   \label{WeakDiffusion2}
\end{align}
where the tilde sign indicates that the corresponding quantity is the numerical flux at the interface of elements. Owing to the strong restriction of time marching interval due to spatial discretization, an implicit time marching scheme is adapted by using the Helmholtz solver \cite{CGHDG}.

\begin{figure}[ht]
\centering
\subfloat[T=$32.0$] {\label{DiffT30} \includegraphics[
width=5cm]{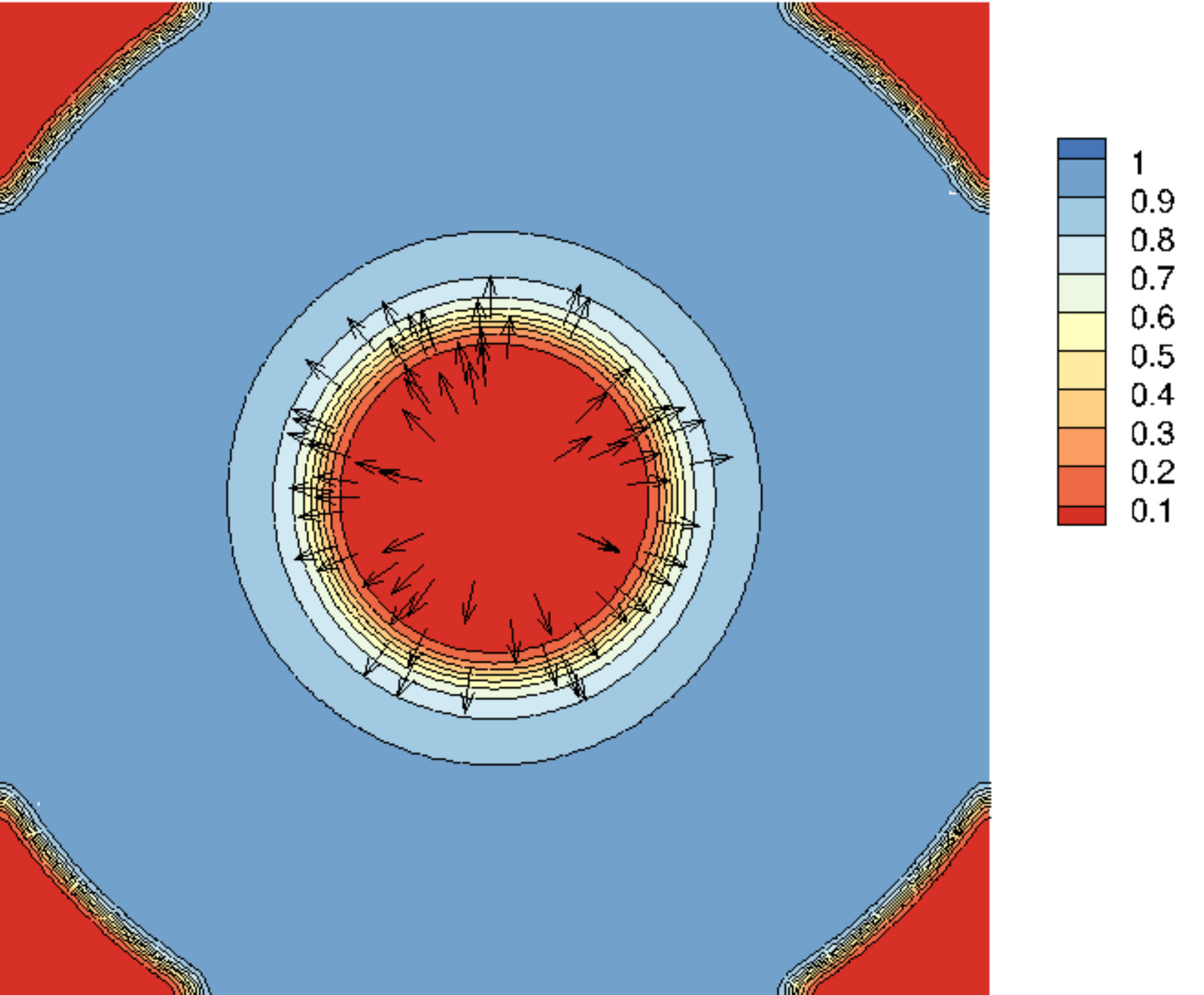} }
\subfloat[T=$40.0$] {\label{DiffMC121} \includegraphics[
width=5cm]{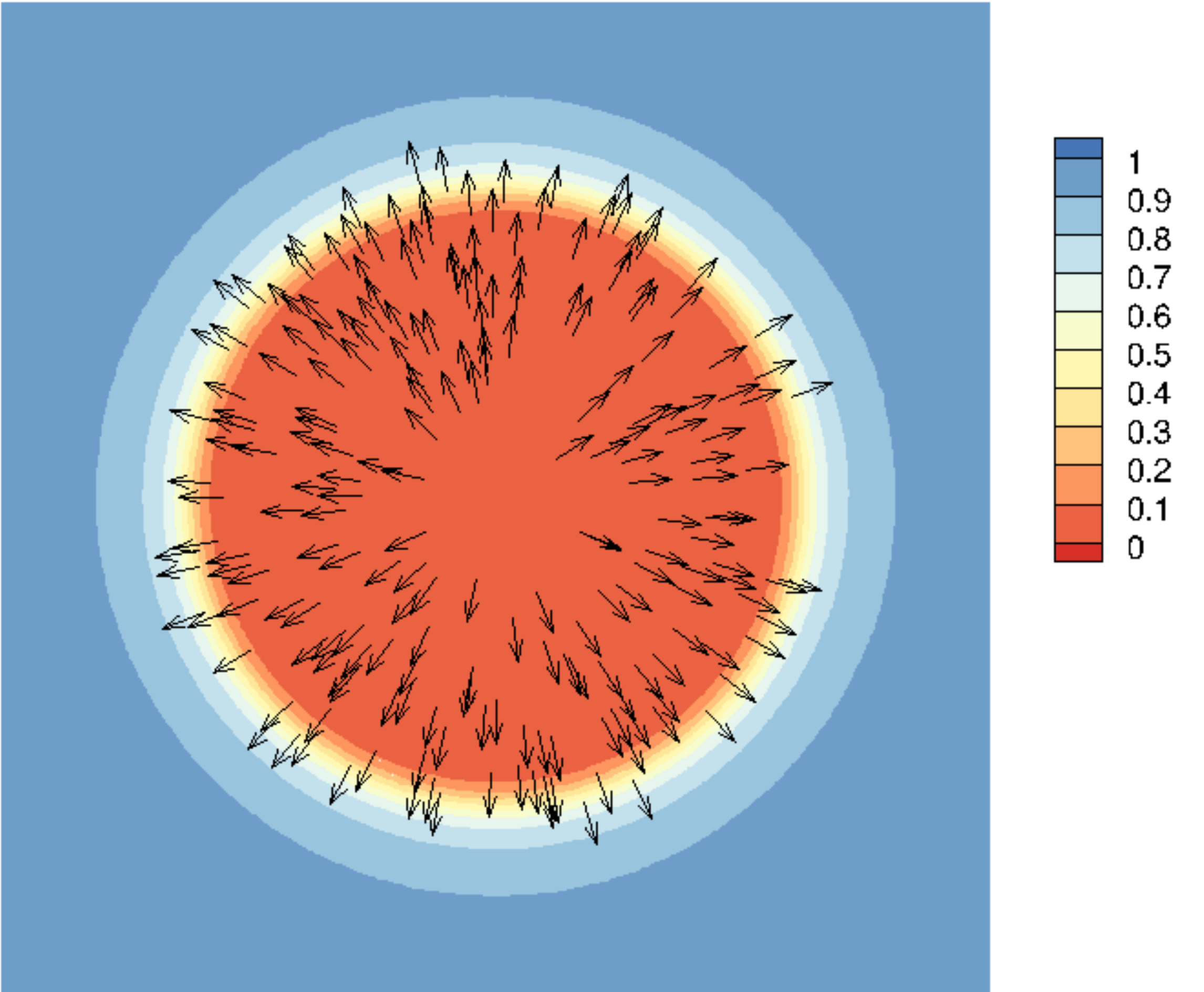} }
\caption{Distribution of the membrane potential (contour) and aligned moving frames (arrow) from the point initialization at $(0,0)$ in the plane. Mesh of $h=4.0$ and $p=4$. WaveBack. Gradient = Str. The domain is $(x,y) = [-40.0, 40.0] \times [-40.0, 40.0]$. }
\label {TestPlane}
\end{figure}

\begin{figure}[ht]
\begin{center}
{
\subfloat[Moving frames ($\mathbf{e}^1$)] {\label{MFPlane}  \includegraphics[height = 4cm, width=4.5cm]{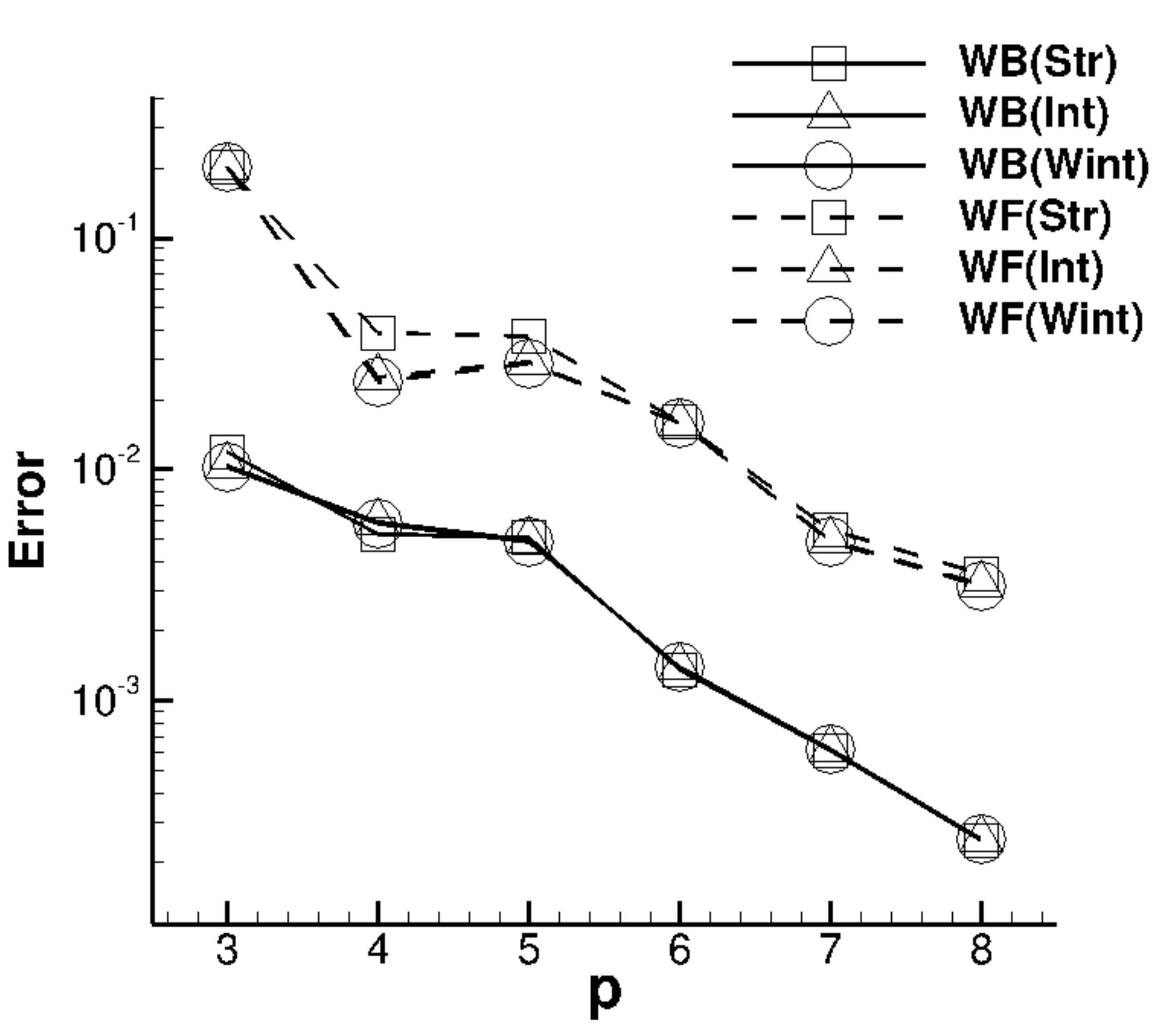}   }\subfloat[Connection ($\omega_{212}$)] {\label{CNPlane}  \includegraphics[height = 4cm, width=4.5cm]{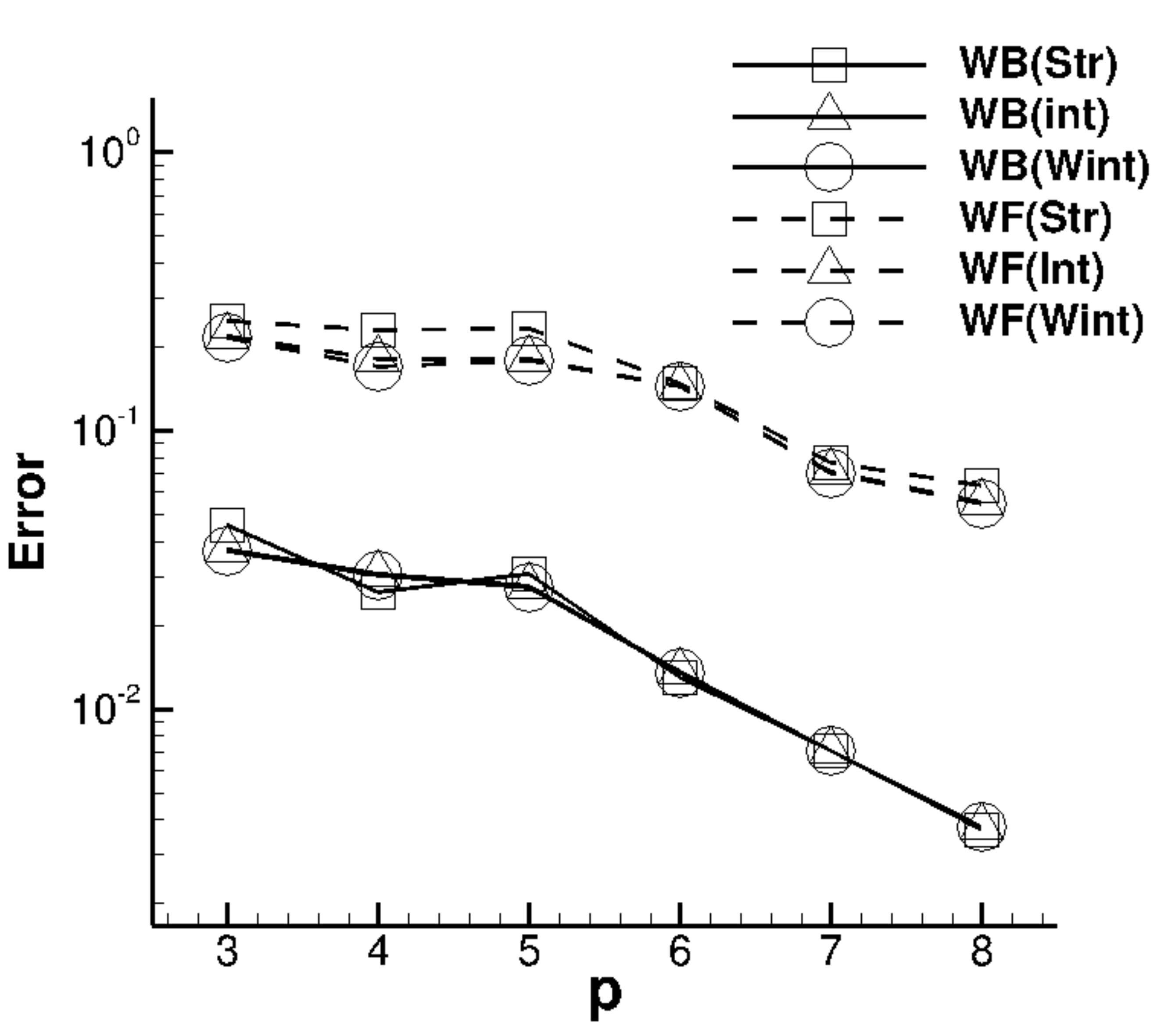}  }\subfloat[Curvature ($\overline{\mathscr{R}}^2_{121}$)] {\label{CVPlane} \includegraphics[height = 4cm, width=4.5cm]{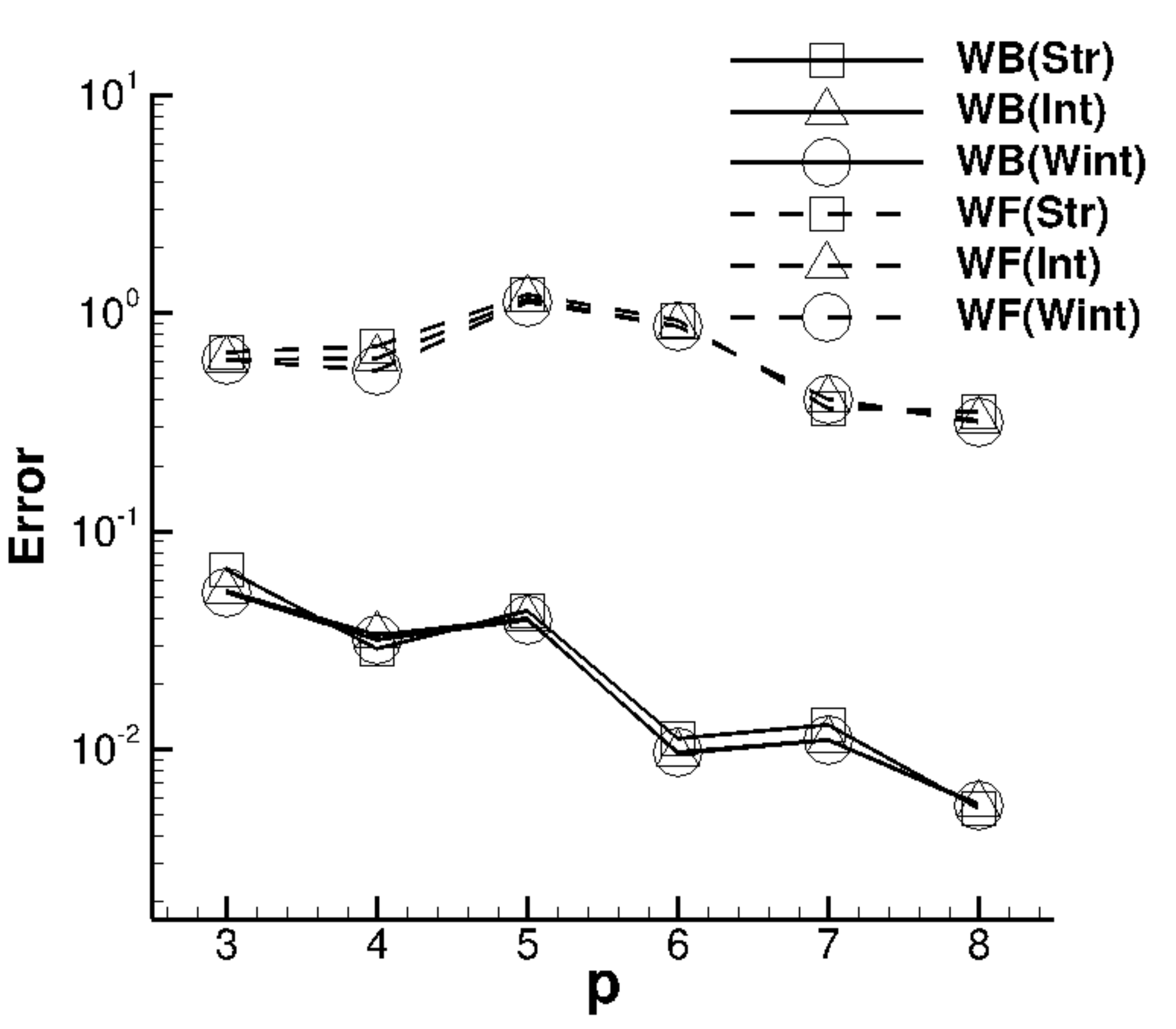}  }  }
\end{center}
\caption{Error $p$-convergence of moving frames, connection, and Riemann curvature for the wave from the point-initiation at $(0,0)$ in the plane. $h$ = $0.4$. Str = Square, Int = Triangle, Wint = Circle. \textsf{WB} = solid line, \textsf{WF} = Dashed line. Measured at $T=40.0$ and $T=20$ for \textsf{WB} and \textsf{WF}, respectively. Area within the radius of $3.0$ from the initiation point is ignored. $\delta = 0.1$.}
\label{Planetestplot}
\end{figure}

\section{Test in the plane}

For the first test in the plane, consider the following Aliev-Panfilov model for the cardiac electric propagation.
\begin{align}
\frac{\partial u}{\partial t} &= \nabla \widehat{\mathbf{d}} \nabla u -  k u ( u - a ) (u - 1) + uv , \label{AP1} \\
\frac{\partial v}{\partial t} &= \left ( \varepsilon_0 + \mu_1 \frac{v}{ u + \mu_2} \right ) ( - v - k u (u - a -1 ) ), \label{AP2}
\end{align}
where $k=8$, $a =0.15$, $\varepsilon_0 = 0.002$. The variable $u$ represents the membrane potential and $v$ represents the collective ion channel openness, including the refractory region. All the variables are dimensionless but can be converted into real measure such as $E = 100 u - 80~(mV)$, $time = 12.9 t~ (ms) $.

Let the wave be point-initialized and propagate in the plane to form a connection, which is the same as the polar coordinate system centered at the initiation point. The moving frames are aligned by Algorithm \ref{MMF1stGV}, and the corresponding connection and Riemann curvature tensor of the moving frames are compared to the exact value. Fig. \ref{TestPlane} displays the propagation of the point-initialized wave with an initial radius of $1.0$ and strength of $1.0$ and the creation of moving frames aligned along the wave trajectory. The gradient is obtained at the waveback when the magnitude of the gradient is the largest.

\begin{table}[htp]
\caption{Plane Test: Error $h $-convergence for the wave from a point-initiation in the plane. $h$ = maximum length of edge. $N_e$ = number of elements. Gradient Computation = \textsf{Wint}. dt=0.01. Measured at $T=40.0$ and $T=20$ for \textsf{WB} and \textsf{WF}, respectively. The area within the radius of $3.0$ from the initiation point is ignored. $\delta = 0.1$.}
\begin{center}
\begin{tabular}{|c|c|c|c|c|c|}
\hline 
  & h ($N_e$) &  4.08 (1946)  & 2.97 (4260)   & 2.21 (7664) & 1.58 (14074) \\
\hline \hline
 \multirow{6}{*} {\textsf{WB}}  &  $\mathbf{e}^1$   & 0.00484769 & 0.000724808 & 0.000153335 & 1.21913e-05  \\
 & order & -  &    5.98 &  5.26  & 7.55\\
 \cline{2-6}
  & $\omega_{212}$ &   0.0271494 & 0.00587723 & 0.00167592  &  0.00019338 \\
  & order & - &   4.82 & 4.25  & 6.44 \\
  \cline{2-6}
 &  $\overline{\mathscr{R}}^2_{121}$ &   0.0394072 & 0.00406901 &  0.0016778  & 0.000226476 \\
 & order & -  &  7.15 & 3.00  & 5.97  \\
 \hline 
   \multirow{6}{*} {\textsf{WF}}  & $\mathbf{e}^1$ &   0.0284854  & 0.00895313  & 0.00228055 & 0.000230573\\
   & order & -  &  3.6449 & 4.6270 & 6.83\\
   \cline{2-6}
 & $\omega_{212}$  &  0.177675  & 0.0853012  & 0.0257 & 0.00416253\\
  & order & -  &   2.3108 & 4.0589 & 5.42\\
  \cline{2-6}
  & $\overline{\mathscr{R}}^2_{121}$ &   1.12586 & 0.24467  & 0.0431835 & 0.00866062\\
  & order & -  &  4.8070 & 5.8682 & 4.79\\
 \hline
\end{tabular} 
\end{center}
\label{TablePlaneTest}
\end{table}%

\begin{table}[htp]
\caption{Sphere test: Error $h$-convergence for a point-initialized wave from the north pole on the sphere. $p$ = $4$. $h$ = maximum length of the edge. $N_e$ = number of elements. Gradient Computation = \textsf{Wint}. dt=0.01. Measured at $T=36.0$ and $T=14$ for \textsf{WB} and \textsf{WF}, respectively. Area within the radius of $3.0$ from the initiation point is ignored. $\delta = 0.1$.}
\begin{center}
\begin{tabular}{|c|c|c|c|c|c|c|}
\hline 
 & h ($N_e$) &   2.48 (1152)  &  1.99 (1152)    & 1.70  (2304)  &1.40  (3072)     \\
 \hline \hline
\multirow{6}{*} {\textsf{WB}}  & $\mathbf{e}^1$ &  0.00148048   &  0.000783078   &  0.000272058  &  7.6965e-05      \\
  & order & -  &    2.9160   &  6.7042  &  6.4565     \\
  \cline{2-6}
  &  $\omega_{212}$ &  0.00589095  &  0.00319238  &  0.00109853   &   0.000619431   \\
  & order & - &  2.8050  &  6.7649  &   2.9296 \\
    \cline{2-6}
 &  $\overline{\mathscr{R}}^2_{121}$  &  0.0222421  &  0.00799031  & 0.00242781 &  0.000995709   \\
  & order & -  & 4.6873   &  7.5540  &  4.5575  \\
 \hline 
\multirow{6}{*} {\textsf{WF}}  & $\mathbf{e}^1$  &   0.0115381  & 0.00717967   & 0.00169167  & 0.000454362   \\
  & order & -  & 2.1721 &   9.1666 &   6.7217    \\
    \cline{2-6}
  & $\omega_{212}$  &  0.0588998  &  0.0234089  &  0.00466264  &  0.00199622    \\
  & order & - &  4.2247 &  10.2319  &  4.3377   \\
    \cline{2-6}
 & $\overline{\mathscr{R}}^2_{121}$ &  0.492037  & 0.329433 & 0.0301431 &  0.0114059   \\
  &  order & -  &   1.8368  & 15.1648  &  4.9691       \\
 \hline
\end{tabular} 
\end{center}
\label{TableSphereTest}
\end{table}

For this wave propagation, the propagational direction of the wave is the same as the radial direction and is expressed as follows.
\begin{equation*}
\mathbf{e}^1 = \frac{\mathbf{r}}{\| \mathbf{r} \| } = ( \cos \theta, \sin \theta ), ~~~\mathbf{e}^2 = ( -\sin \theta, \cos \theta ),
\end{equation*}
where $\theta = \tan^{-1} (y/x)$. The connection should have the following values
\begin{equation*}
\omega_{211} = 0.0, ~~~~ \omega_{212} = \frac{1}{\| \mathbf{r} \| },
\end{equation*}
and the Riemann curvature tensor of orthonormal bases has the following non-trivial component.
\begin{equation*}
 \overline{\mathscr{R}}^2_{121} = \frac{1}{\| \mathbf{r} \|^2 } .
\end{equation*}
The Riemann curvature tensor is the same as the above tensor because the component of the metric tensor for the $\hat{\mathbf{r}}$ axis is the same along the $\hat{\boldsymbol{\theta}}$ axis. The gradient is measured at the wavefront when $T=20.0$ and at the waveback when $T=40.0$, respectively.

Fig. \ref{Planetestplot} displays the exponential $p$-convergence of moving frames (left), connection (center), and Riemann curvature (right). Three different methods are used for computing the propagation direction from the gradient: \textsf{Str} (to choose the gradient when the gradient is the largest), \textsf{Int} (to average the integration of the gradients which is larger than a threshold value), \textsf{Wint} (to perform a weighted integration proportion to the magnitude of the gradient). The result shows that \textsf{Wint} shows slightly better accuracy, however, overall, the difference is very small for this test case. Moreover, the accuracy at the waveback (solid line) is considered better than the wavefront (dashed line), possibly because the waveback is much smoother than the wavefront. In the wavefront, the membrane potential increases dramatically; thus, the differentiation of the insufficiently resolved membrane potential may not yield an accurate derivative of the membrane potential.

Table \ref{TablePlaneTest} shows the exponential $h$-convergence of error for moving frames ($\mathbf{e}^1$), connection ($\omega_{212}$), and Riemann curvature ($\overline{\mathscr{R}}^2_{121}$) by \textsf{Wint}. The convergence at $p$=$5$ on meshes with various $h$ confirms the convergence order of $7.55,~ 6.44,~5.97$ for the waveback and $6.83, 5.42$, and $4.79$ for the wavefront. The exact convergence order should be $p+1$, $p$, $p-1$. The convergence order of connection is one order less than the order of moving frames because the computation of connection requires additional differentiation in space. Similar order reduction from connection occurs for the Riemann curvature.

\begin{figure}[ht]
\centering
\subfloat[T=$28.0$] {\label{SphereT5} \includegraphics[
width=5cm]{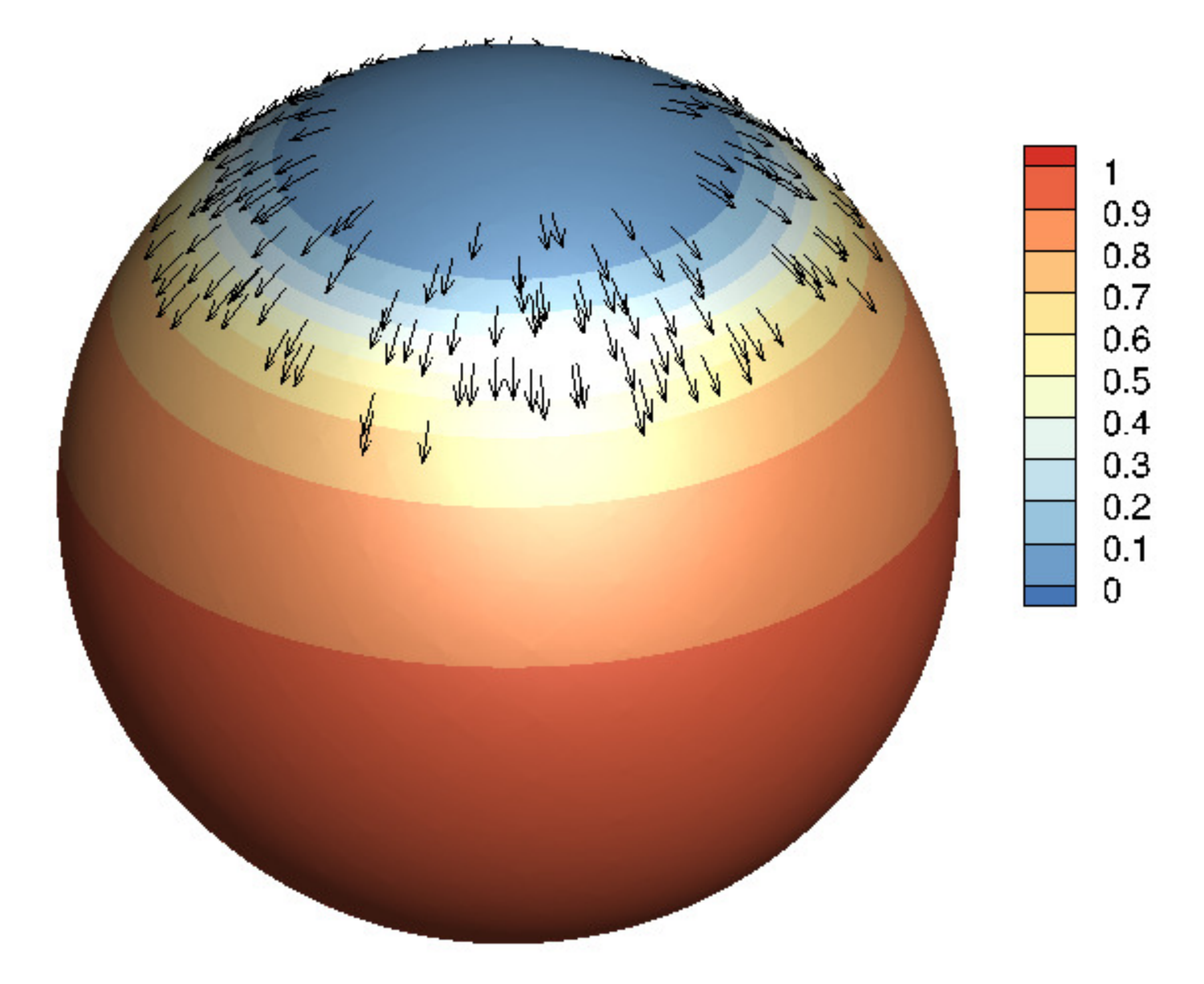} }
\subfloat[T=$36.0$] {\label{SphereT6} \includegraphics[
width=5cm]{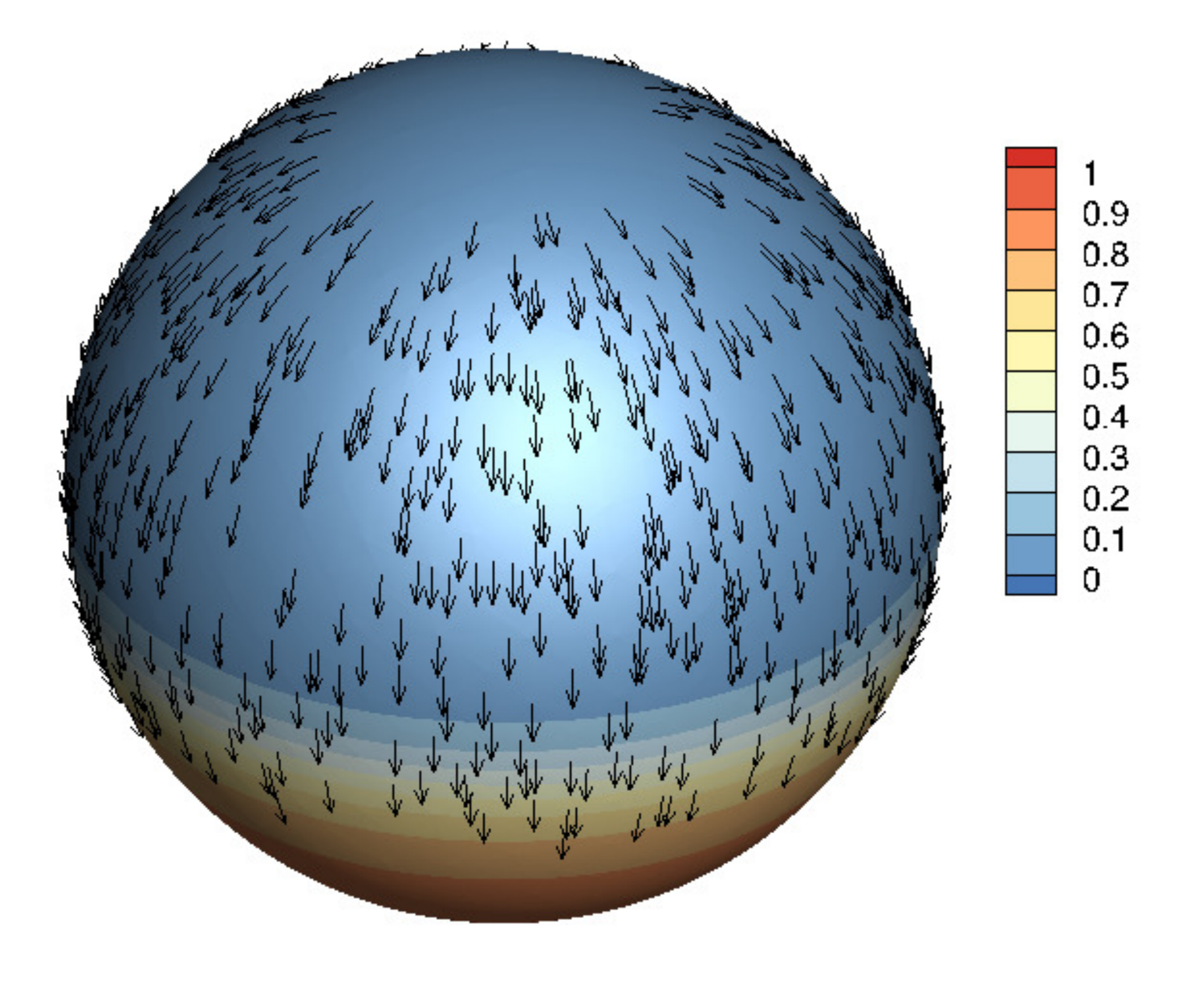} }
\caption{Distribution of the membrane potential ($u$) and moving frames for a point-initialized wave at the north pole. The mesh of $h=4.0$ and $p=4$ is used. WaveBack. Gradient = Wint. The domain is the spherical shell with a radius of $10.0$.}
\label {TestSphere}
\end{figure}

\begin{figure}[ht]
\begin{center}
{
\subfloat[Moving frames ($\mathbf{e}^1$)] {\label{MFSphere}  \includegraphics[height = 4cm, width=4.5cm]{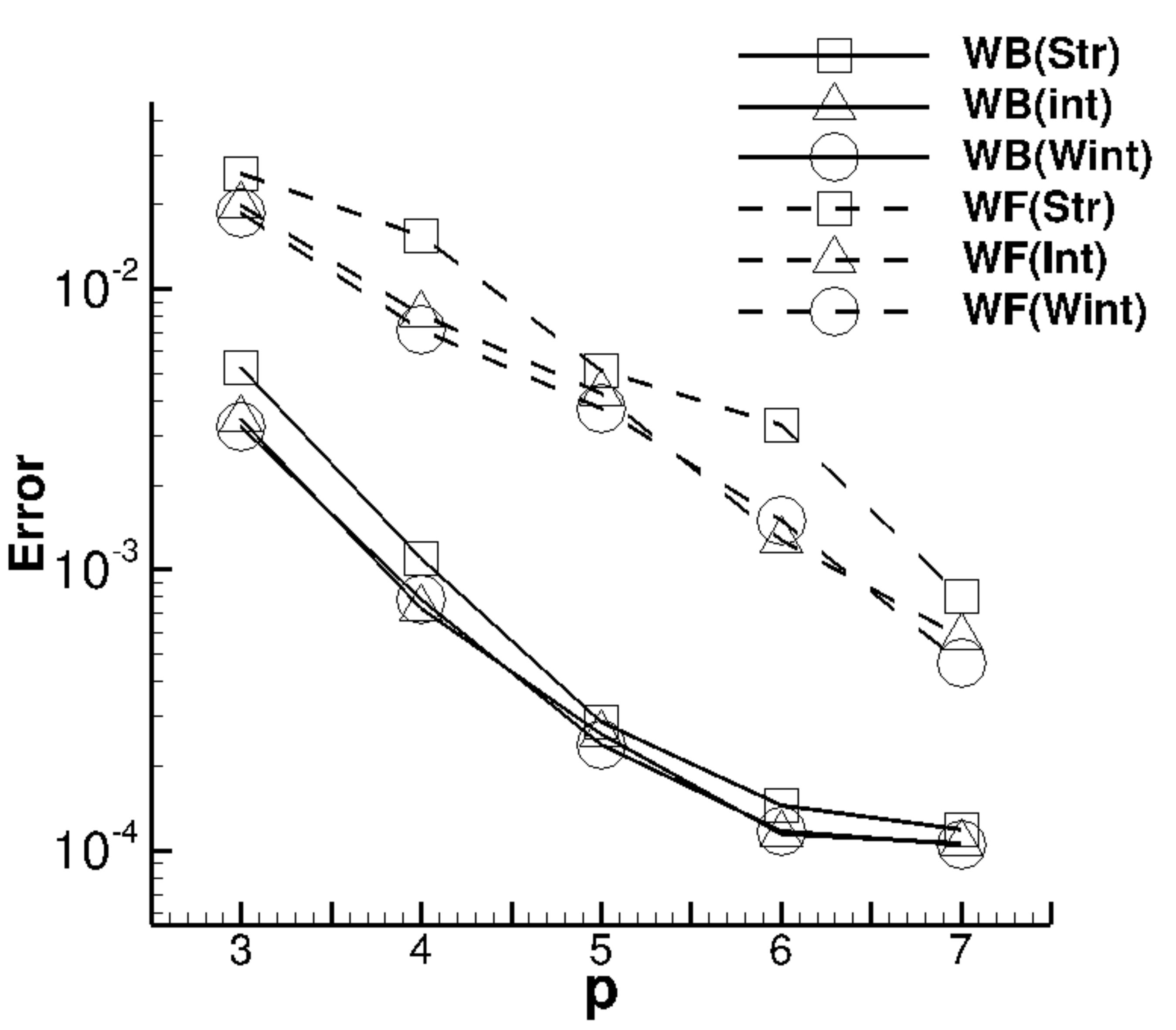}   }\subfloat[Connection ($\omega_{212}$)] {\label{CNSphere}  \includegraphics[height = 4cm, width=4.5cm]{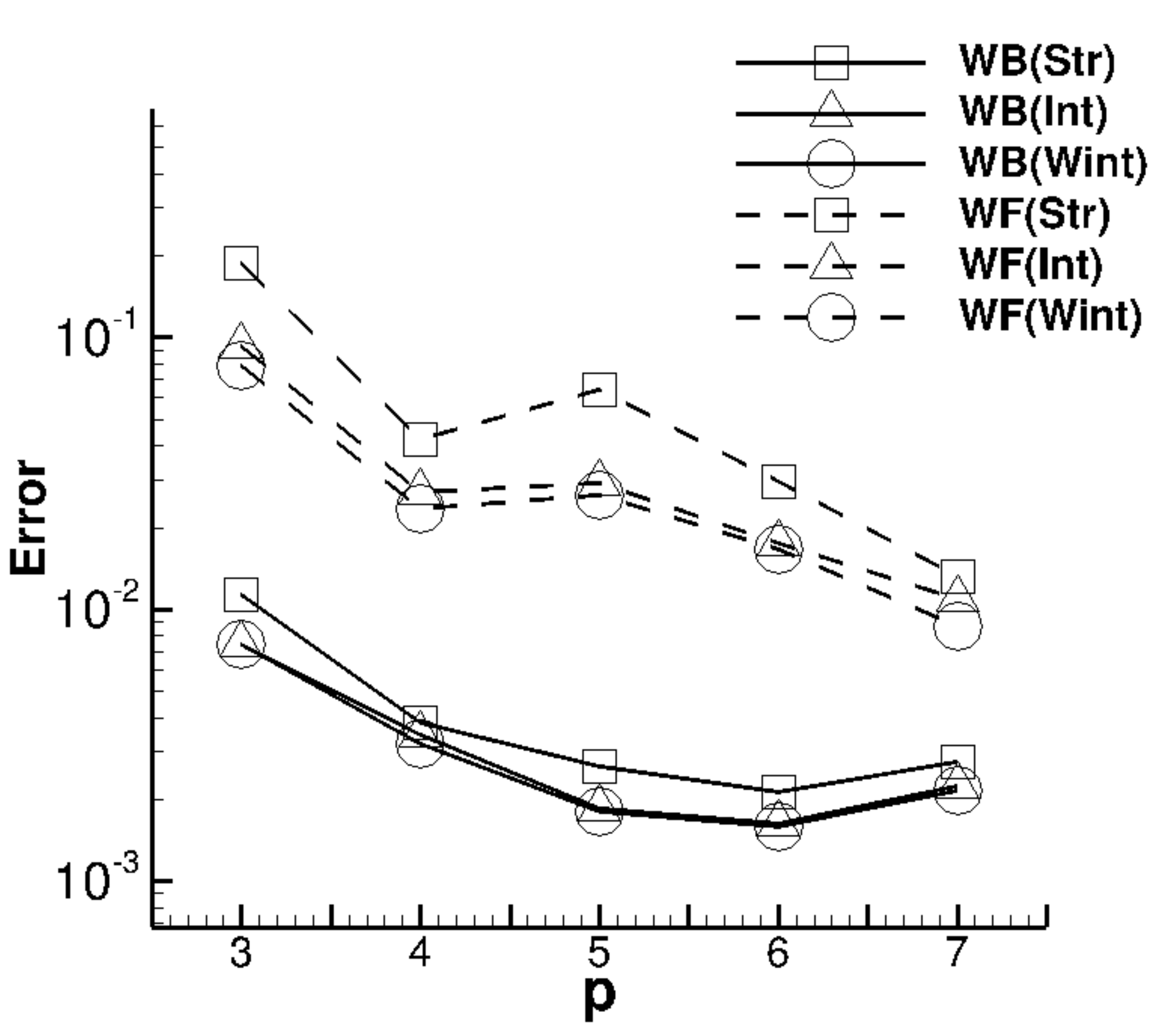}  }\subfloat[Curvature ($\overline{\mathscr{R}}^2_{121}$)] {\label{CVSphere} \includegraphics[height = 4cm, width=4.5cm]{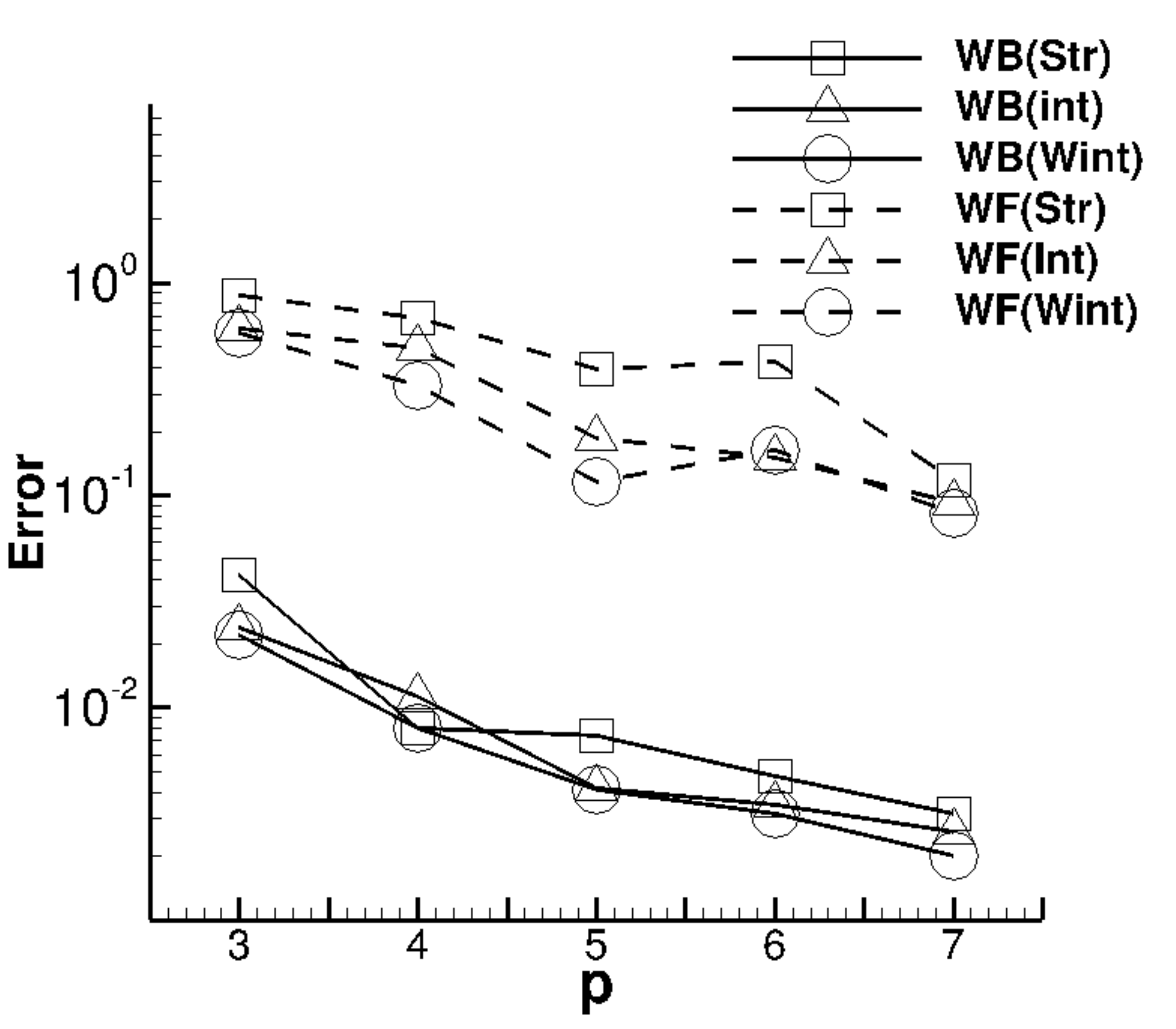}  }  }
\end{center}
\caption{Error $p$-convergence for the wave from a point-initiation on the sphere. $h$ = $0.4$. \textsf{Str} = Square, \textsf{Int} = Triangle, \textsf{Wint} = Circle. \textsf{WB} = Solid line, \textsf{WF} = Dashed line. Measured at $T=36.0$ and $T=14$ for \textsf{WB} and \textsf{WF}, respectively, when the moving frames of $60 \sim 70 \%$ of the area are aligned along the propagational direction. The area within the radius of $3.0$ from the initiation point is ignored. $\delta = 0.1$. }
\label {ErrSphere}
\end{figure}

\section{Test on the sphere}

The second test is to initiate a point-excitation at the north pole on the sphere of radius $R$. Fig. \ref{TestSphere} displays the propagation of the wave from the north pole and the alignment of the moving frames by computing the gradient in the waveback. Then, the sphere is parametrized as follows.
\begin{align*}
(x,y,z) = ( r  \sin \theta \cos \varphi,  r  \sin \theta \sin \varphi,  r \cos \theta )  .
\end{align*}
 The moving frame should be aligned along the latitudinal direction of the sphere, or $\partial / \partial \theta$, as follows.
\begin{align*}
\mathbf{e}^1 &= \partial / \partial \theta = ( \cos \theta \cos \varphi, \cos \theta \sin \varphi, - \sin \theta ) , \\
 \mathbf{e}^2 &=\partial / \partial \varphi = ( - \sin \varphi, \cos \varphi, 0 )  .
\end{align*}
Note that $\mathbf{e}^2$ is a unit vector. The tangent vector of the $\varphi$ axis has the magnitude of $\sin \theta$. The corresponding connection should have the following values.
\begin{align*}
\omega_{211} &= 0.0, ~~~~ \omega_{212} = \frac{ 1 }{R  }  \frac{\cos \theta }{\sin \theta} , \\
\omega_{311} &= - \frac{1}{R } ,~~~~\omega_{322} = - \frac{1}{R } .
\end{align*}
The Riemann curvature of orthonormal basis has the following non-trivial component.
\begin{align*}
\overline{\mathscr{R}}^2_{121} = \frac{1}{R^2 \sin^2 \theta }  . 
\end{align*}
Note that $\overline{\mathscr{R}}^2_{121}$ is not the same as $\mathscr{R}^2_{121}$, that is, equal to $\cos^2 \theta$. Fig. \ref{ErrSphere} shows a similar exponential convergence of moving frames ($\mathbf{e}^1$), connection ($\omega_{212}$), and Riemann curvature ($\overline{\mathscr{R}}^2_{121}$) as the exponential convergence for the plane. For error convergence, \textsf{WB} shows much better accuracy than \textsf{WF}. Out of the three gradient computation methods,  \textsf{Str} shows the worst accuracy and \textsf{Wint} shows the best accuracy. As shown in Fig. \ref{CNSphere}, error stagnation for the connection component $\omega_{212}$ for $p \ge 5$ is possibly due to the geometric approximation error of the spherical mesh \cite{GAerror} as it appears again at a finer mesh for $h$-convergence. However, this error stagnation in $\omega_{212}$ does not affect exponential convergence in $\overline{\mathscr{R}}^2_{121}$ as shown in Fig. \ref{CVSphere}.

Table \ref{TableSphereTest} demonstrates the $h$-convergence for this test problem on the sphere. Similar to the test case in the plane, $\mathbf{e}^1$, $\omega_{212}$, and $\overline{\mathscr{R}}^2_{121}$ shows $p+1$, $p$, and $p-1$ order or higher, respectively. As observed in the $p$-convergence, a low order of $\omega_{212}$ between $h=1.70$ and $h=1.40$ by \textsf{WB} is possibly due to the geometric approximation error of the spherical mesh.  The order of $\omega_{212}$ between $h=1.70$ and $h=1.40$ by \textsf{WF} is $4.3377$, which is close to the exact value of $p$.

\begin{figure}[ht]
\centering
\subfloat[$x^2 + (y+25)^2 = r^2$ ] {\label{Anitype1} \includegraphics[
width=5cm]{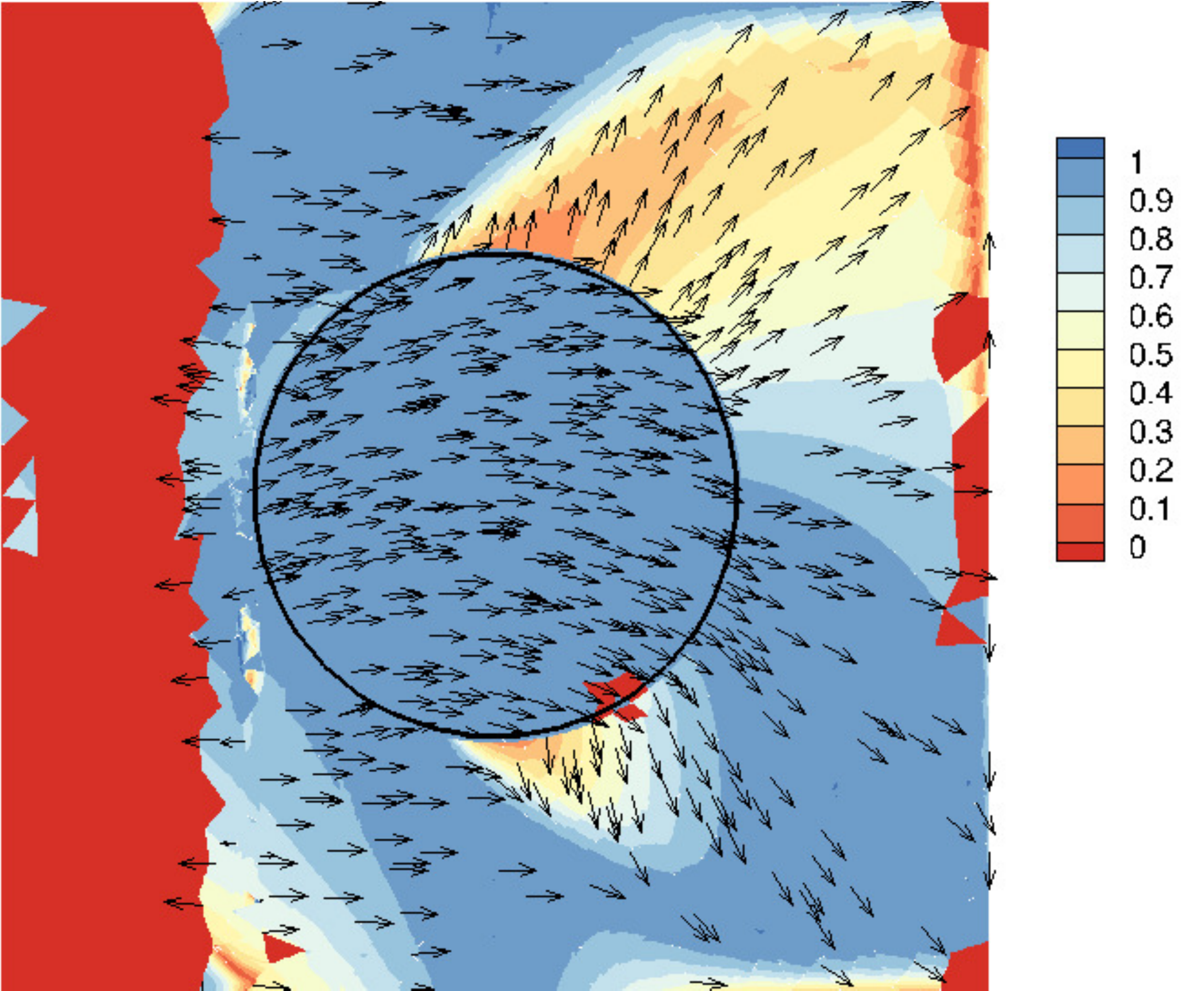} }  
\subfloat[$x^2 + (y-25)^2 = r^2$ ] {\label{Anitype1} \includegraphics[
width=5cm]{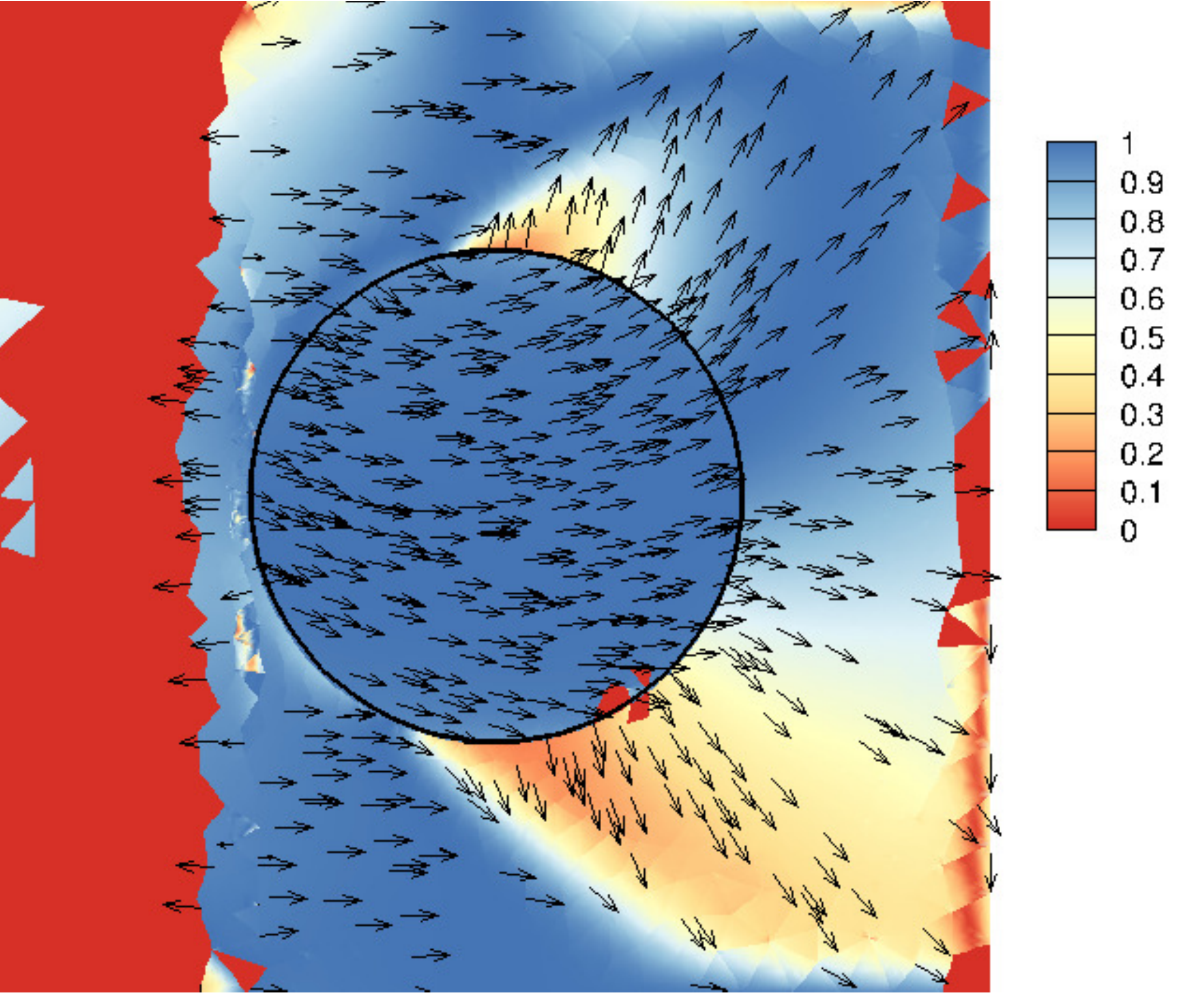} }  
\caption{ Domain = [-20,20] $\times$ [-20,20]. Anisotropy is aligned along the circular arc centered at $(0,-25)$ (a) and $(0,25)$ (a). The plane wave is initiated from the left wall at $x = - 25$.}
\label {Anisotropy1}
\end{figure}

\section{Test in anisotropy medium}

In an isotropic medium, the direction of propagation depends on many factors, especially on the initiation type and point. Even with the same geometry and wave type, the change of initiation point would lead to a different direction of propagation and consequently, a different connection map. However, in the anisotropic medium where one direction is dominantly faster than the other direction in terms of propagation speed, the location of initiation is considered much less significant to the connection map of the propagation. Let us verify this.

In the weak formulation of the MMF-diffusion scheme in the anisotropy medium, the moving frames are aligned along the direction of the anisotropy, and the length of the moving frame is adjusted according to the strength of the anisotropy. Then, the coefficient of the gradient $d^{mm}$ in the mixed formulation in Eq. \eqref{WeakDiffusion2} is modified with the new expression of Eq. \eqref{dani}. A non-uniform length of moving frames yields the same effect as multiplying the diffusivity tensor to the gradient.
 
\begin{figure}[ht]
\centering
\subfloat[Initialized at $(-10,-10)$] {\label{Anitype2} \includegraphics[
width=5cm]{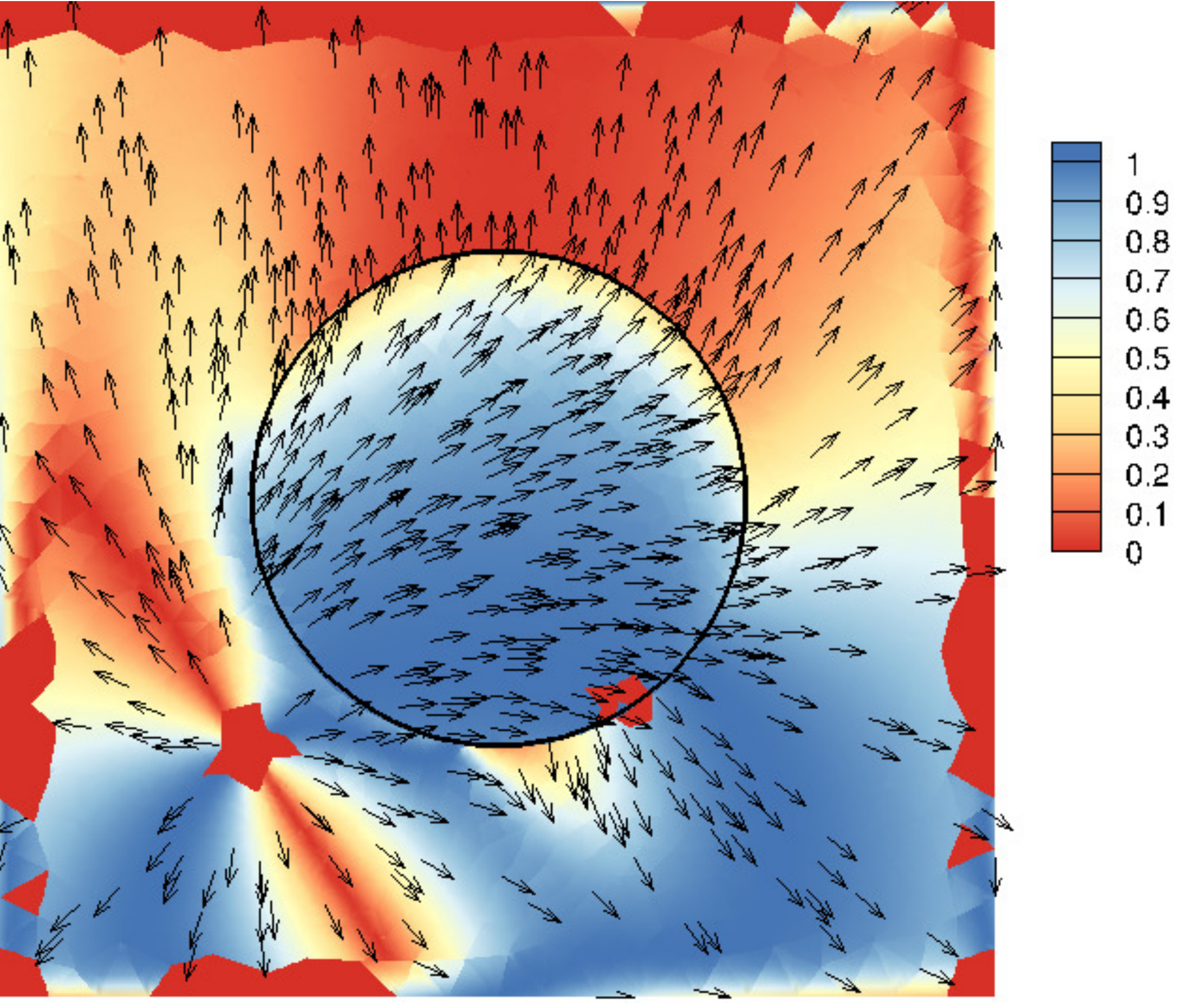} }  
\subfloat[Initialized at $(-10, 10)$] {\label{Anitype3} \includegraphics[
width=5cm]{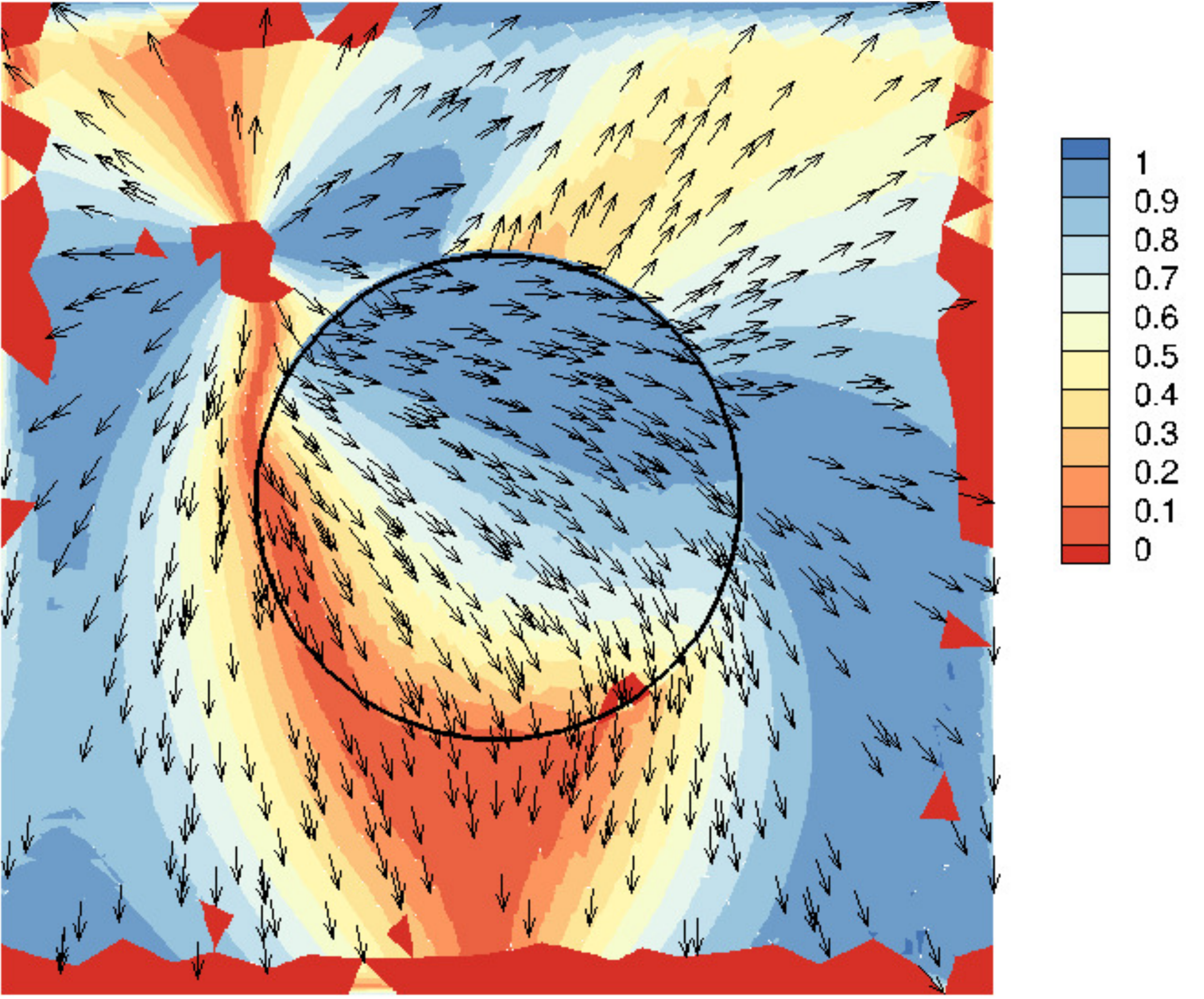} }  
\caption{ Domain = [-20,20] $\times$ [-20,20]. Point-initialized. Anisotropy is aligned along $x^2 + (y+25)^2 = r^2$ within the circle.}
\label {Anisotropy2}
\end{figure}

Fig. \ref{Anisotropy1} presents the map of aligned moving frames along the propagational direction on two different anisotropies within the center circle in the domain of [-20,20] $\times$ [-20,20]: For the left plot, the anisotropy is aligned along the circular arc centered at $(0,-25)$, i.e., the anisotropy is on the line of $x^2 + (y+25)^2 = r^2$. For the right plot, the anisotropy is aligned along the circular arc centered at $(0,25)$, i.e., the anisotropy is on the line of $x^2 + (y-25)^2 = r^2$. The plane wave from the left wall follows the direction of the anisotropy within the circle. Observe the coincidence (blue color) of the propagational direction and anisotropy within the circle in Fig. \ref{Anisotropy1}.

The main difference between cardiac electric propagation and electromagnetic wave is that cardiac electric propagation does not always follow the direction of anisotropy. Fig. \ref{Anisotropy2} displays the propagation after the point initialization near the circle where the anisotropy is aligned along the line $x^2 + (y+25)^2 = r^2$. In the left plot, the wave initiates in the low left corner, i.e., at $(-10,-10)$. In the right plot, the wave initiates in the top left corner, i.e., at $(-10,10)$. To follow the anisotropy within the circle, the propagation should take a very sharp turn or the propagational velocity along the orthogonal to the fiber should be almost zero. However, in some areas within the circle, the wave fails to follow the anisotropy. In Fig. \ref{Anisotropy2},  non-blue color indicates a non-trivial angle between the anisotropy and the direction of propagation. The deviation is more significant for the initiation, which started from the top-left corner because anisotropy is curved toward the bottom of the domain. The directional difference between the propagational direction and anisotropy could be the key feature in analyzing and predicting the abnormal propagation of the cardiac electric flow in the heart.

\begin{figure}[ht]
\centering
\subfloat[Original fiber map] {\label{Orgfiber} \includegraphics[
width=5.5cm]{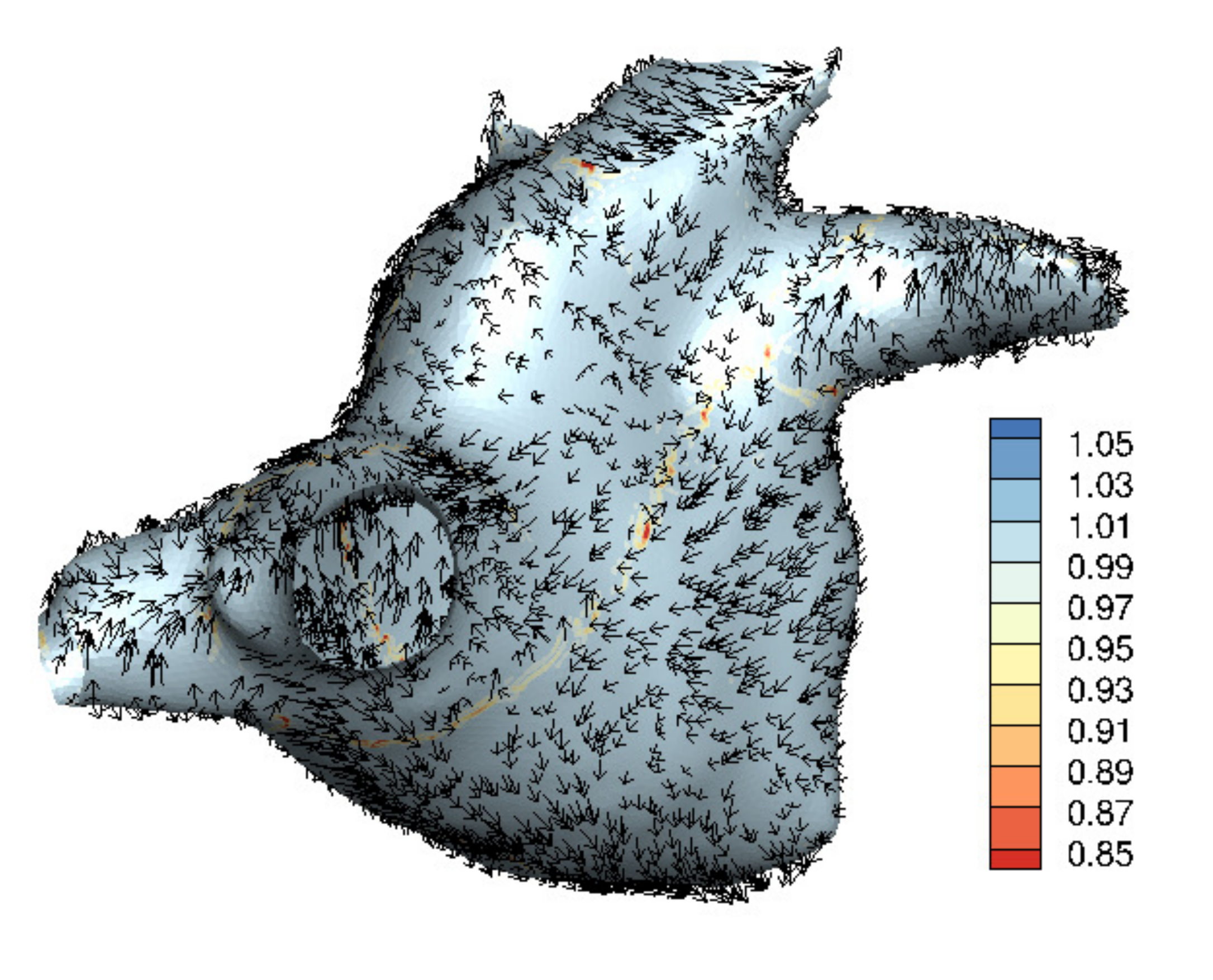} }
\subfloat[Projected fiber map] {\label{Projfiber} \includegraphics[
width=5.5cm]{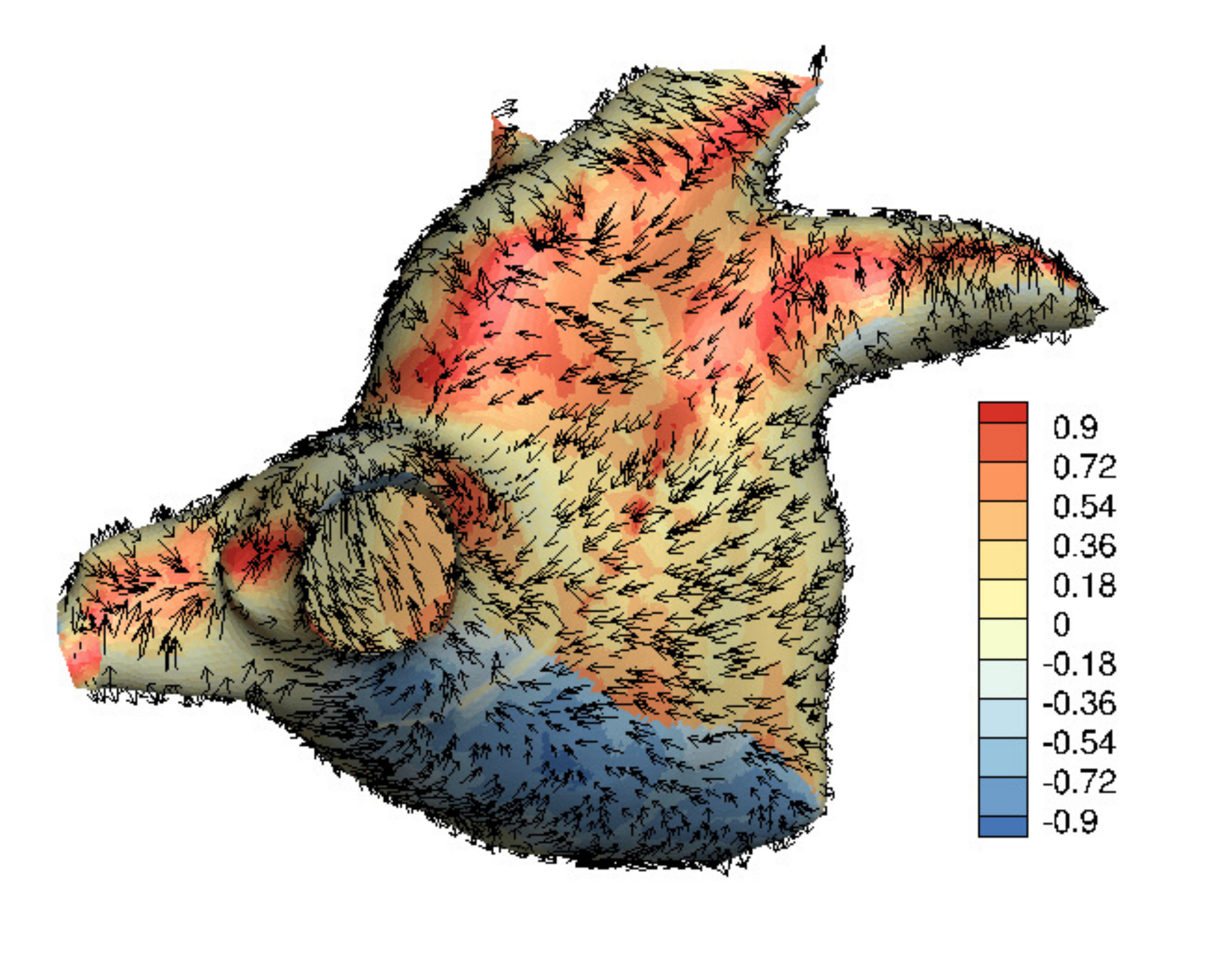} }
\caption{In (a), the contour represents the magnitude of the cardiac fiber. In (b), the contour represents $\sin \theta$ for the angle $\theta$ between the cardiac fiber and the surface normal vector ($\mathbf{k}$).}
\label {Cardiacfiber}
\end{figure}

\section{Applications to Atlas of atrium}

This section describes how the proposed algorithm can be applied to the generation of an Atlas of a 2D atrium, or the connection map of the cardiac electric signal propagation. The most significant benefit of the Atlas is to show the spatial-temporal map of the divergence and convergence of the propagation, possibly to show the location where the propagation stops according to Hypothesis 1. In contrast to the contour plots of the wavefront, this Atlas is unique and valid in the presence of anisotropy; thus, it can be used to analyze the pattern of propagation in a complex cardiac tissue.

In the following sections, we demonstrate the several possible uses of the Atlas for the analysis of the cardiac electric signal propagation in the atrium. For this purpose, the three-dimensional shape and cardiac fiber of the atrium are constructed from magnetic resonance angiography images and histological examination of ex-vivo atria, respectively \cite{Cantwell}, as shown in Fig. \ref{Orgfiber}. These mesh and fiber data are not the actual data of the human heart but serves the purpose as an exemplary use of the proposed scheme to the realistic atrium.

\subsection{Preprocessing}

The length of the cardiac fiber representing the ratio of the strength of conductivity is not constant, but slightly variable between $0.78$ and $1.04$. The cardiac tissue of the atrium is anatomically known to be relatively thin. It can be modeled as a curved surface, but still has the features of a three-dimensional structure. Thus, the direction of the cardiac fiber may not be orthogonal to the surface normal vector. In other words, if the cardiac fiber is represented as a vector $\mathbf{f}$, then in a two-dimensional model of atrium, the vector $\mathbf{f}$ should be orthogonal to the surface normal vector $\mathbf{k}$ such as $\mathbf{f} \cdot \mathbf{k} = 0$. For this reason, some of the original fiber that has an oblique angle $\theta$ with $\mathbf{k}$ is projected on the tangent plane, but the conductivity velocity was also adjusted by the factor of $\sin \theta$ because the velocity parallel to the tangent plane is only considered in the 2D model such as
\begin{equation*}
\mathbf{f}_{proj} = \mathbf{f} - (\mathbf{f} \cdot \mathbf{k}) \mathbf{k} =   \sin \theta  \mathbf{f}.
\end{equation*}
Fig. \ref{Projfiber} displays the projected fiber (arrow) and $\cos \theta$ (contour) between the 3D cardiac fiber and the surface normal vector. To simplify notation, let $\mathbf{f}$ mean $\mathbf{f}_{proj}$ in the remainder of this paper.

\section{Analysis of 2D cardiac fiber}

The following analysis of a 2D cardiac fiber yields the two fundamental features of the cardiac fiber on \textit{controlling the propagation}: The first is to block the propagation, and the second is to remove the sources of self-initiation.

\subsection{To block the propagation}
\begin{figure}[ht]
\centering
\subfloat[$\omega_{212}$] {\label{FiberCN} \includegraphics[
width=5cm]{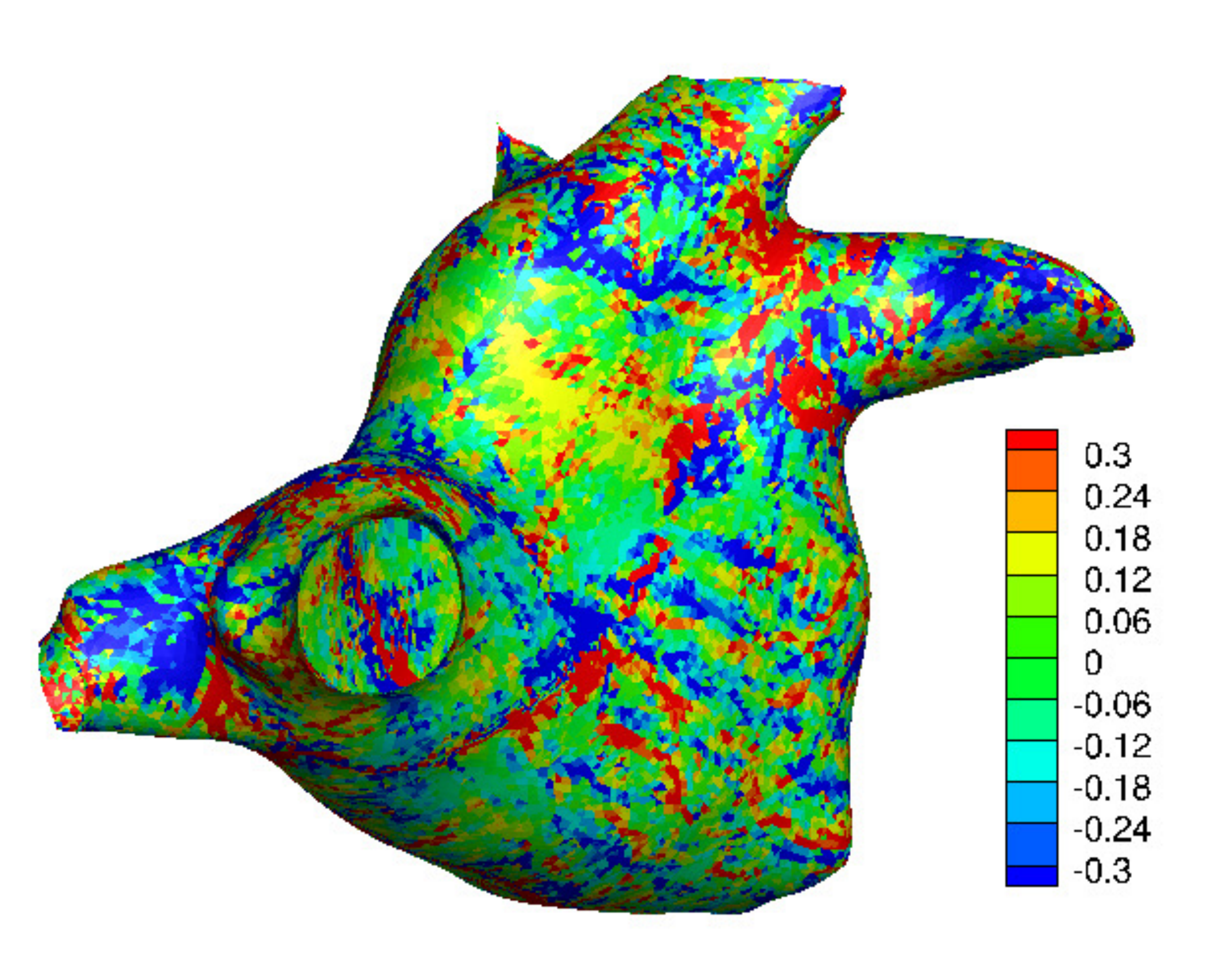} }
\subfloat[$\overline{\mathscr{R}}^2_{121}$] {\label{FiberRC} \includegraphics[
width=5cm]{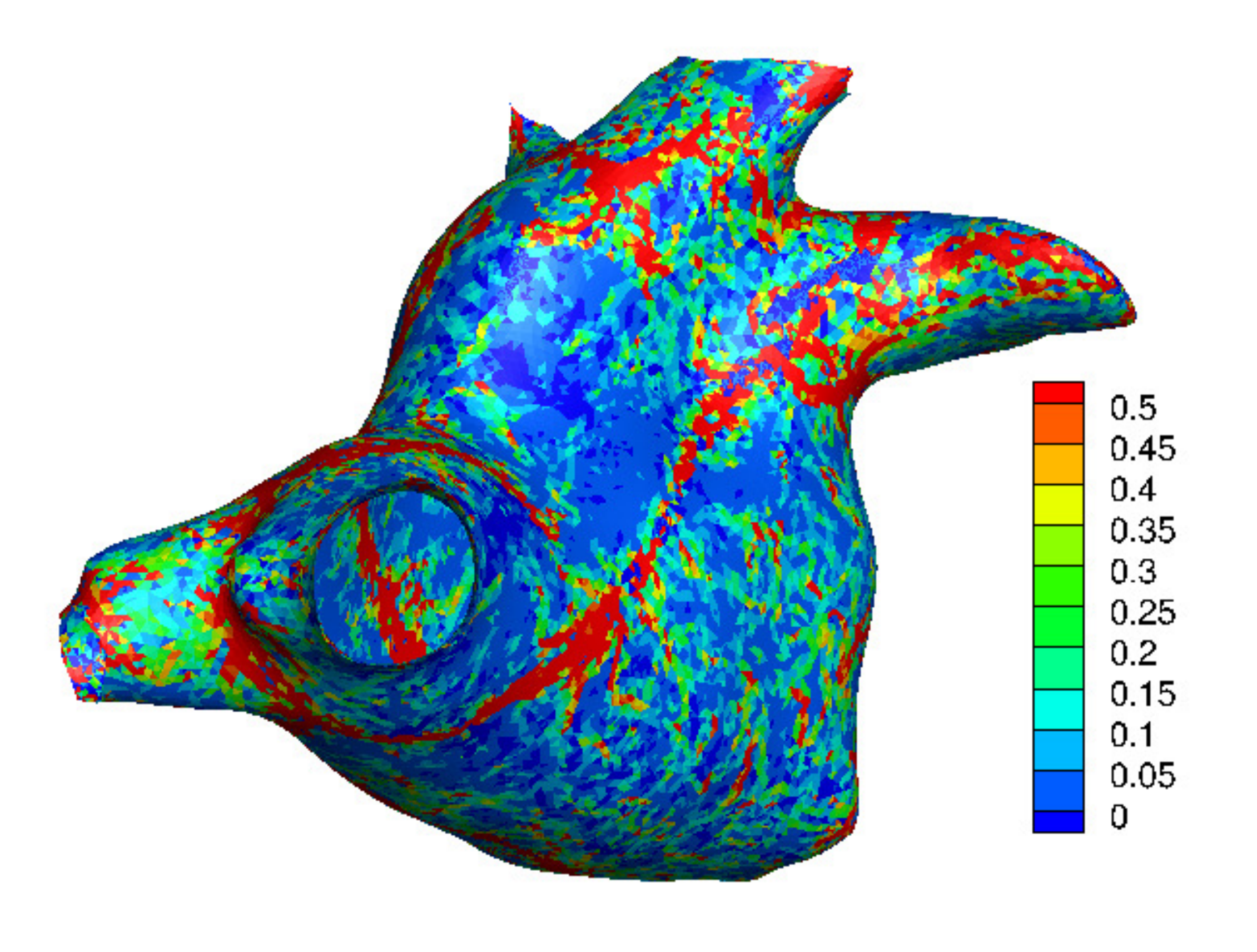} }
\caption{Projected cardiac fiber: Distribution of the connection component (a) $\omega_{212}$ and (b) $\overline{\mathscr{R}}^2_{121}$.}
\label {CardiacfiberCNRC}
\end{figure}

Applying Eq. \eqref{omegaijk} to the obtained cardiac fiber reveals the connection map and the Riemann curvature of orthonormal bases, as shown in Fig. \eqref{CardiacfiberCNRC}. Fig. \ref{FiberCN} shows the distribution of $\omega_{212}$, which is overall relatively low, i.e., below $0.3$, in most of the atrium except the entrance of each vein. Observe that the region with a large magnitude of $\omega_{212}$ is approximately overlapped with the region with a large positive $\overline{\mathscr{R}}^2_{121}$ as shown in Fig. \ref{FiberRC}. A positive $\overline{\mathscr{R}}^2_{121}$ means a diverging fiber, and a negative $\overline{\mathscr{R}}^2_{121}$ means a converging fiber. Thus, only the red area is related to the possible stopping area for the cardiac electric propagation. If the cardiac electric signal propagates along the cardiac fiber, then the propagation may stop or dramatically slowdown in the red area. However, if the cardiac signal flows in a different direction than the fiber, then this connection map does not tell anything for certainty, and it is possible that the electric flow passes through the red area in another direction.

\begin{figure}[ht]
\centering
\subfloat[ ] {\label{fiberDiv2} \includegraphics[
width=5cm]{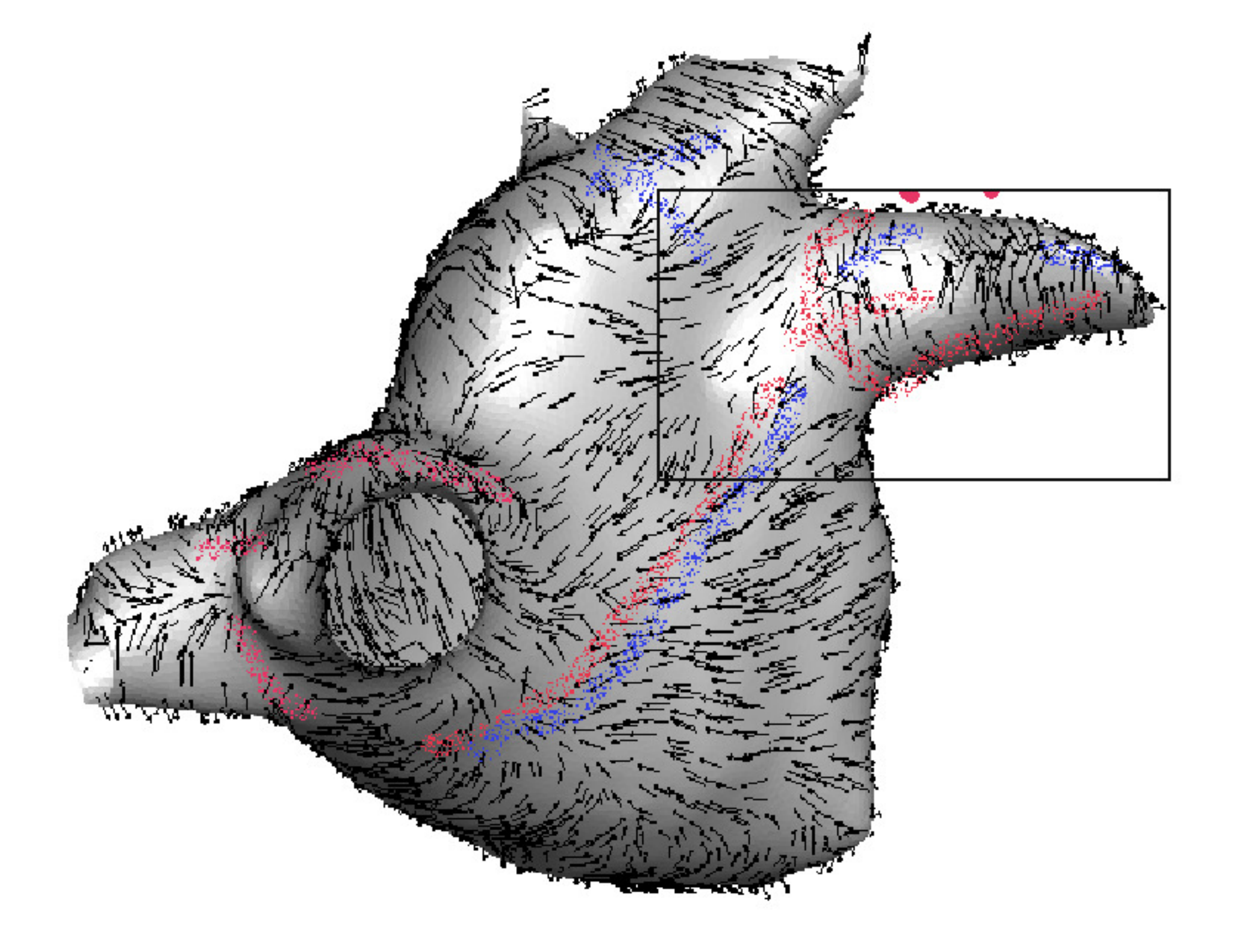} }
\subfloat[ ] {\label{FiberRC2} \includegraphics[
width=5cm]{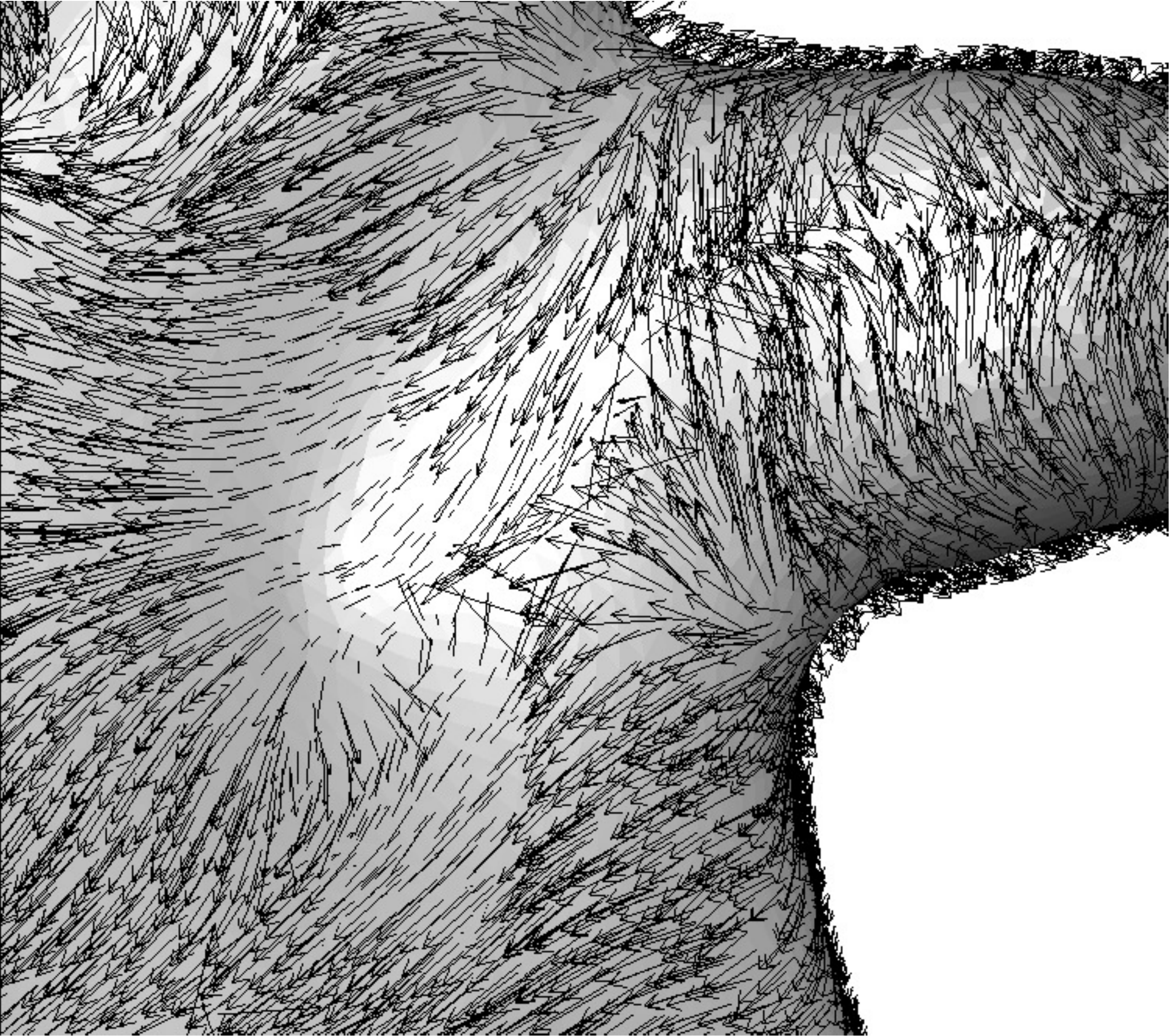} }
\caption{Magnified view for cardiac fiber with high curvature: (a) region of interest, (b) magnified around the vein.}
\label {Cardiacfiber2}
\end{figure}

Fig. \ref{Cardiacfiber2} presents the magnified view of the cardiac fiber around the vein. The red and blue highlights roughly correspond to the region with high $\omega_{212}$ and $\overline{\mathscr{R}}^2_{121}$. Fig. \ref{FiberRC2} shows that the red area for $\omega_{212}$ is related to the abrupt change of the cardiac fiber around the vein. The abrupt change of microscopic fiber is anatomically observed at the canine heart \cite{Hamabe}, sheep posterior left atrium \cite{Klos}, and the human heart \cite{Tan}. The connection map confirms the causal relationship between the anatomical structure and propagational pattern for a unidirectional block around the pulmonary vein to initiate the atrial reentry \cite{MyBIOP, Nattel2008,Quan1990}.

\subsection{To remove the sources}

Other interesting geometric properties of the cardiac fiber can be better understood by computing the convariant divergence and curl of the cardiac fiber, which can be computed as follows \cite{MMFCovariant}.
\begin{align*}
\nabla \cdot \mathbf{f} &= \left ( \nabla f_1 \cdot  \mathbf{e}^1 + \Gamma^1_{21} f_2 + \nabla f_2 \cdot  \mathbf{e}^2 + \Gamma^2_{21} f_1  \right ) , \\
( \nabla \times \mathbf{f} ) \cdot \mathbf{k}   &=   \left ( \nabla f_2 \cdot  \mathbf{e}^1 + \Gamma^2_{12} f_2 - ( \nabla f_1 \cdot  \mathbf{e}^2 + \Gamma^1_{21} f_1 )  \right ), 
\end{align*}
where $\Gamma^i_{jk}$ is the Christoffel symbol of the second kind and is computed using the connection of the connection form $\omega^i_j$ because of the equality $\Gamma^i_{jk} = \omega^i_j \langle \mathbf{e}^k \rangle$ for orthonormal basis $\mathbf{e}^i$.

\begin{figure}[ht]
\centering
\subfloat[$\nabla \cdot \mathbf{e}^1$] {\label{fiberDiv} \includegraphics[
width=5cm]{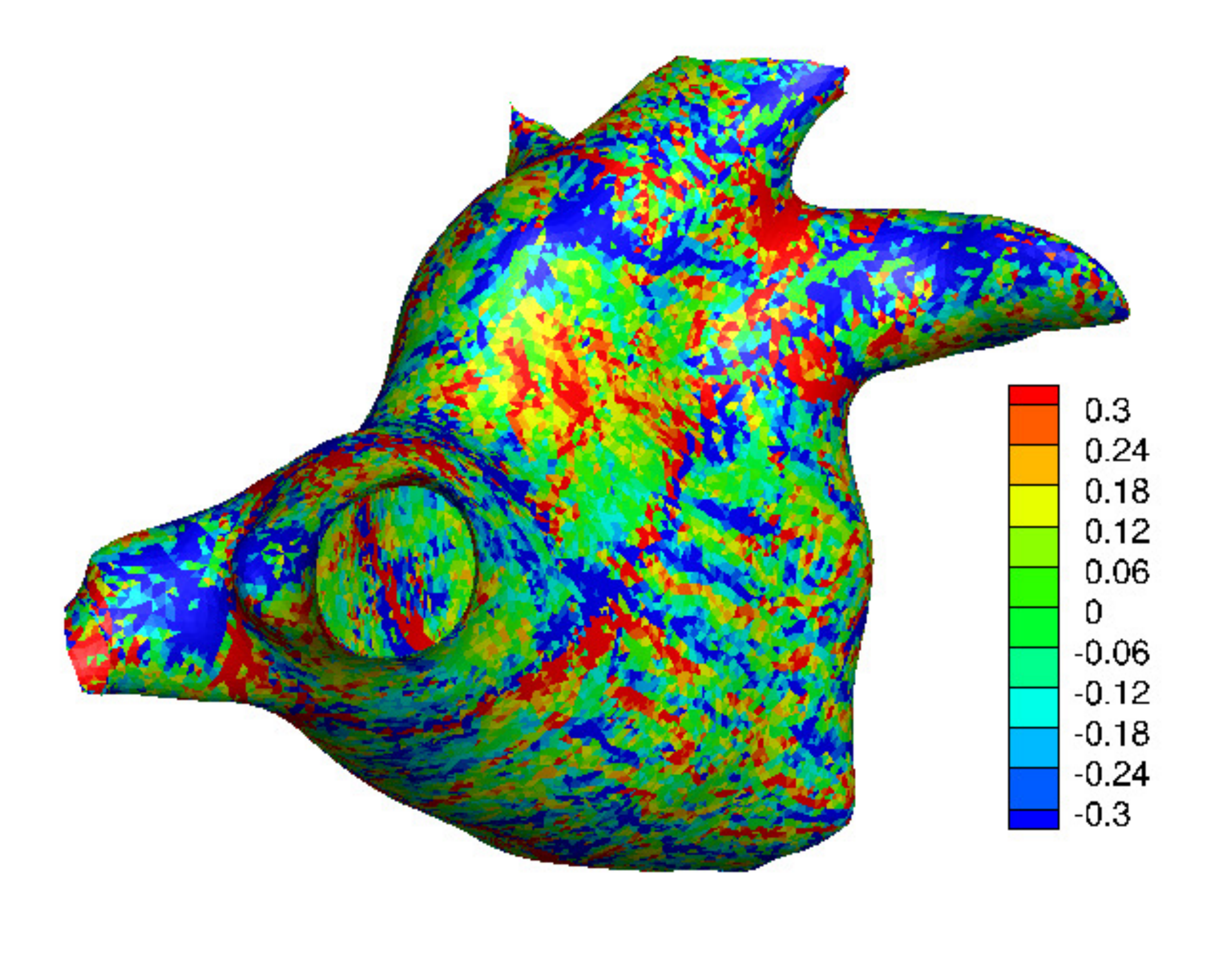} }
\subfloat[$\mathbf{k} \cdot \nabla \times \mathbf{e}^1$] {\label{FiberRC} \includegraphics[
width=5cm]{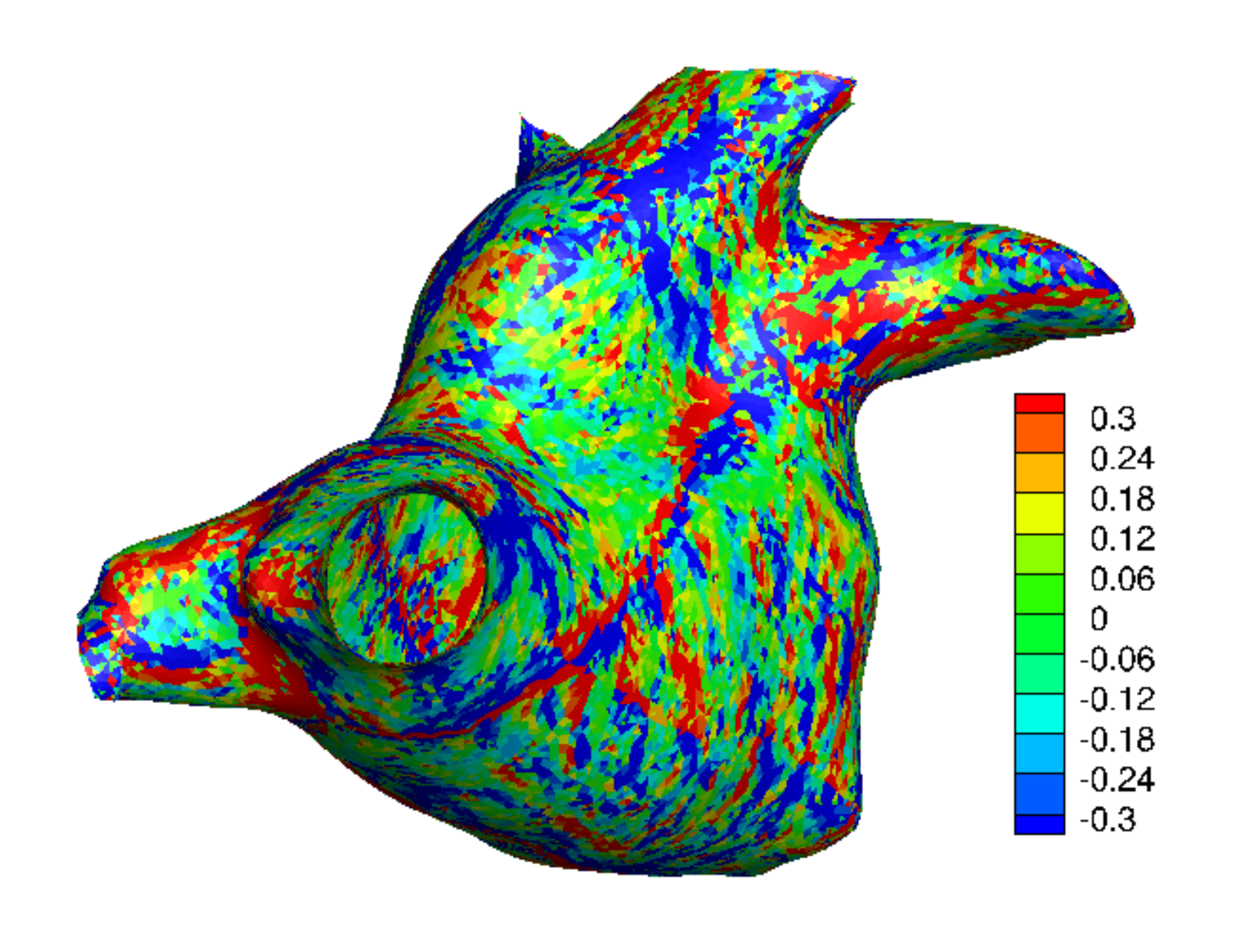} }
\caption{(a) Divergence of the fiber. (b) Projection of the curl of the fiber onto the surface normal direction.}
\label {CardiacfiberDivCurl}
\end{figure}

Fig. \ref{fiberDiv} shows a similar distribution of the divergence of the cardiac fiber except around the veins with a large magnitude of $\omega_{212}$. With regard to electromagnetics, the cardiac electric flow is related to the flow of charged ions to induce the electric field ($\mathbf{E}$) along the fiber. The divergence of the fiber, or the divergence of the electric field if the electric propagation completely follows the fiber, means the amount of charge density ($\rho$) as shown by Gauss's law as
\begin{equation*}
\nabla \cdot \mathbf{E}  = \rho \approx 0  .
\end{equation*}
Thus, the negligible divergence of the cardiac fiber eliminates the possibility of self-excitation due to a high charge density by the fiber structure. Because many cardiac electric disorders are initiated by self-excitation such as a spiral wave, it is obvious that the cardiac fiber of a healthy heart should avoid a structure with a high covariant divergence.

Another interesting geometric feature arises concerning suppressing the magnetic field under the structure of the cardiac fiber. According to Maxwell's equations, the magnetic field ($\mathbf{B}$) is generated by the non-negligible curl of the electric field by Faraday's law such as 
\begin{equation*}
\frac{\partial \mathbf{B}}{\partial t} = - \nabla \times \mathbf{E}. 
\end{equation*}
The above equation implies that if $\mathbf{E}$ is curl-free, the magnetic field remains zero if no magnetic field is initially present. Fig. \ref{FiberRC} displays the negligible amount of $\nabla \times \mathbf{f}$ that is approximately the same as $\nabla \times \mathbf{E}$ if the electric propagation completely follows the fiber. The suppression of the magnetic field in the propagation of the cardiac electric flow is crucial because the magnetic field could change the propagational direction of charged ions to deviate from the cardiac fiber \cite{Jackson}. Observe that the area with a strong magnetic field is again approximately overlaps with the area with a high $\overline{\mathscr{R}}^2_{121}$. This means that an abnormal cardiac electric flow, such as a spiral wave, could generate a relatively significant amount of magnetic field.

\section{Connection map by trajectory tracing}
To derive the connection map of the cardiac electric signal propagation in the obtained atrium, the twenty-one variable Courtemanche-Ramirez-Nattel (CRN) model is solved on the mesh of the atrium with a initialization from a point close to the sinoatrial node (SAN). The CRN model is known to be more accurate in the atrium than the Aliev-Panfilov model. However, similar results would be obtained by the two-variable Aliev-Panfilov model. The CRN model is described as follows \cite{Coutemanche}: Consider the following diffusion-reaction equation. 
\begin{equation*}
\frac{d u}{d t} = \nabla \cdot \widehat{\mathbf{d}} \nabla u  + I_{ion} ,
\end{equation*}
where $I_{ion}$ is the total ionic current which is given as
\begin{align*}
I_{Ion} &= I_{Na} + I_{K1} + I_{to} + I_{Kur} + I_{Kr} + I_{Ks} + I_{Ca, L} \\
&~~~~~~ + I_{p,Ca} + I_{NaK} + I_{NaCa} + I_{b,Na} + I_{b,Ca} ,
\end{align*}
where $I_{Na}$ is the fast inward $Na^+$ current, $I_{K1}$ is the inward rectifier $K^+$ current, $I_{to}$ is the transient outward $K^+$ current, $I_{Kur}$ is the ultrarapid delayed rectifier $K^+$ current,  $I_{Kr}$ is the rapid delayed rectifier $K^+$ current, $I_{Ks}$ is the slow delayed rectifier $K^+$ current, $I_{Ca, L}$ is the L-type inward $Ca^{2+}$ current, $I_{p,Ca}$ is the sarcoplasmic $Ca^{2+}$ pump current, $I_{NaK}$ is the $Na^+$-$K^+$ pump current, $I_{NaCa}$ is the $Na^+$/$Ca^{2+}$ exchanger current, $I_{b,Na}$ is the background $Na^+$ current, and $I_{b,Ca}$ is the background $Ca^{2+}$ current. For exact parameter values, refer to \cite{Coutemanche}.

For numerical simulation, we used the following specifications: the order of solution polynomial is $p=4$ with $39,388$ regular elements. The diffusion operator is marched in time implicitly by the IMEX third-order time marching scheme. The weak formulation is adapted by the MMF-diffusion scheme in the context of discontinuous Galerkin methods \cite{MMF2}. The gradient of propagation is computed at \textsf{WB} by \textsf{Wint}. The time step is set as $0.01$, and the final time is $T=100.0$. $\delta$ is given as $10.0$.

\subsection{What is the role of the fiber?}

\begin{figure}[ht]
\centering
\subfloat[ Aligned moving frames] {\label{APIsoMF} \includegraphics[
width=4.5cm]{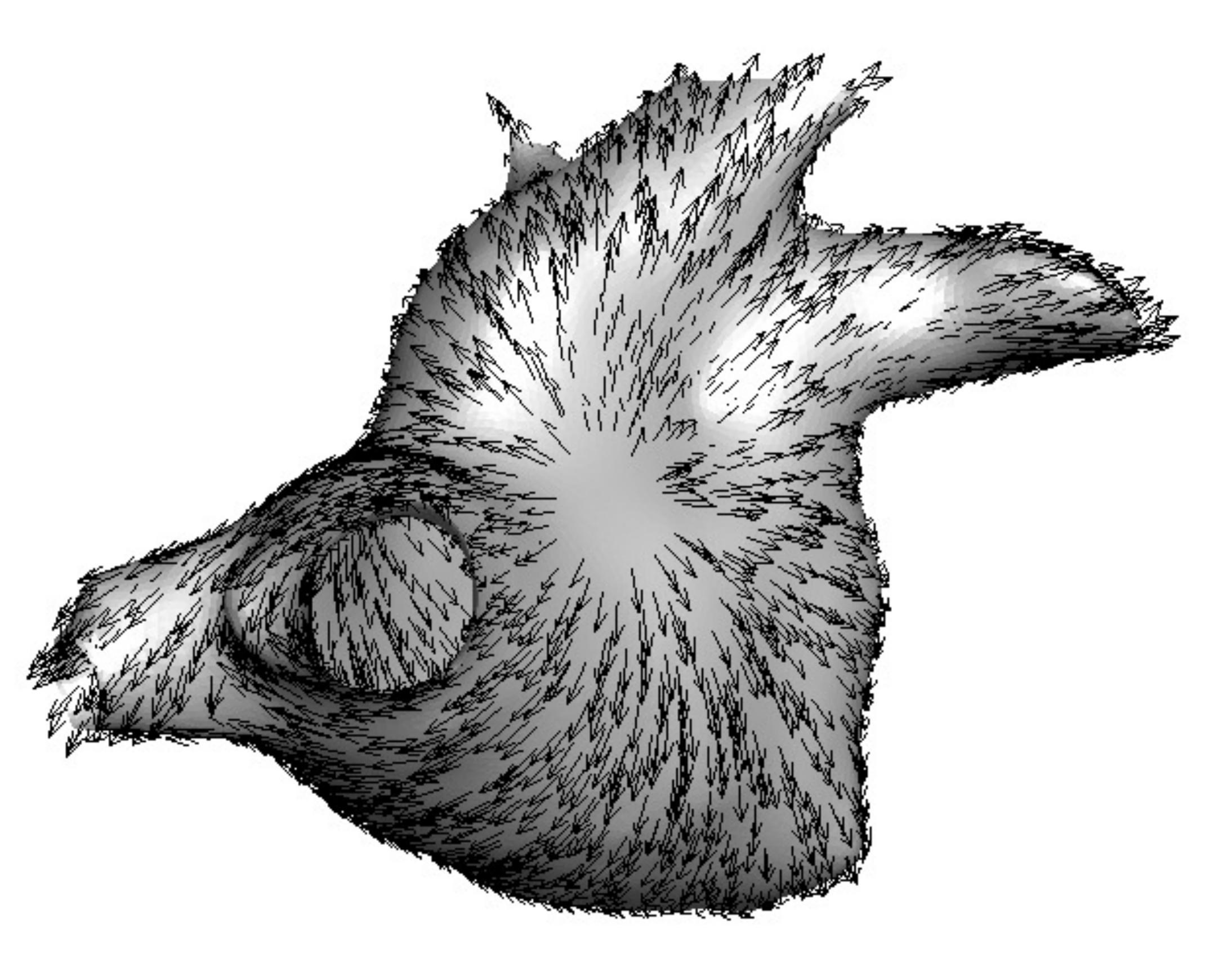} }
\subfloat[$\omega_{212}$ ] {\label{APIsoCN} \includegraphics[
width=4.5cm]{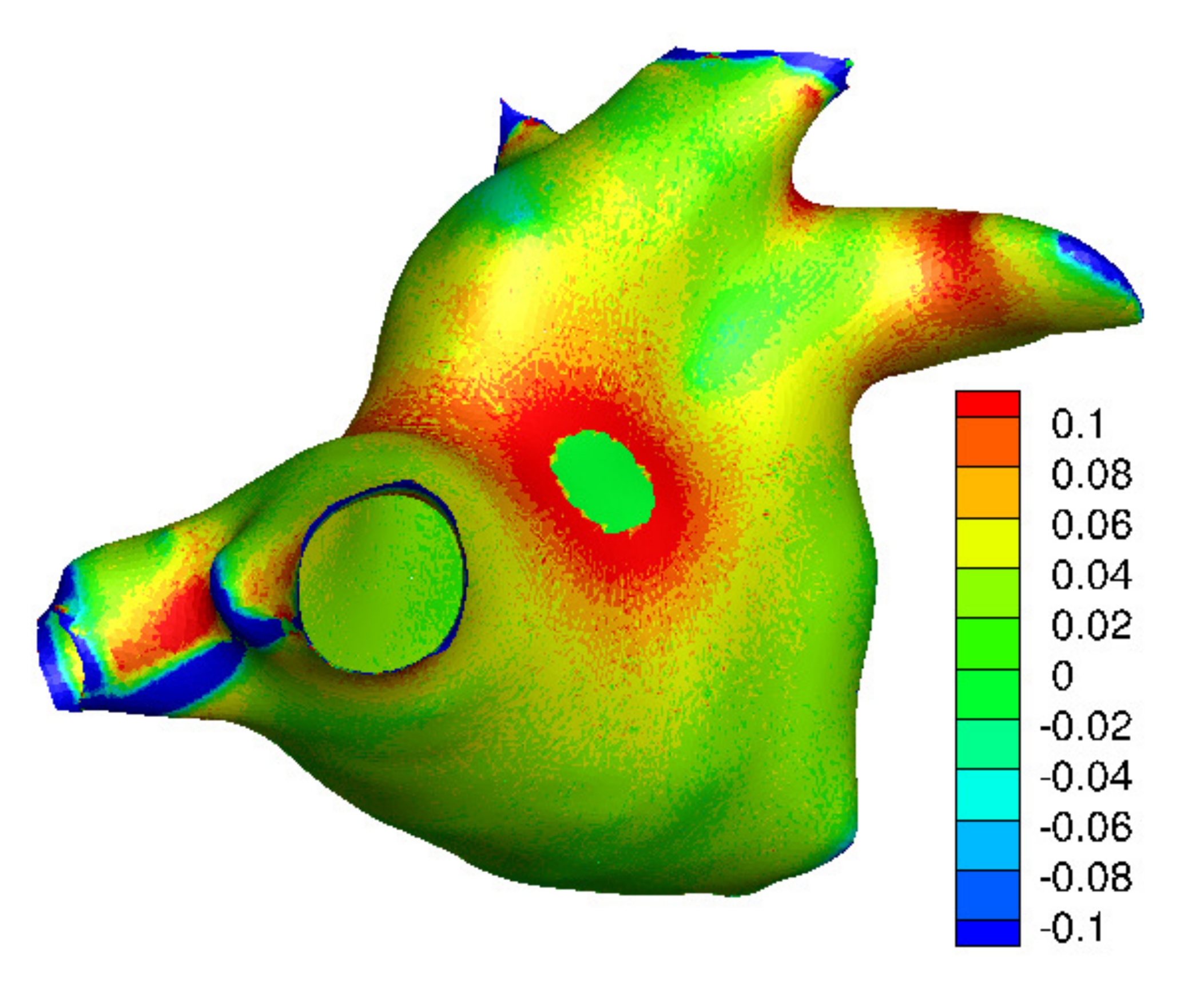} }
\subfloat[$\overline{\mathscr{R}}^2_{121}$ ] {\label{APIsoRM} \includegraphics[
width=4.5cm]{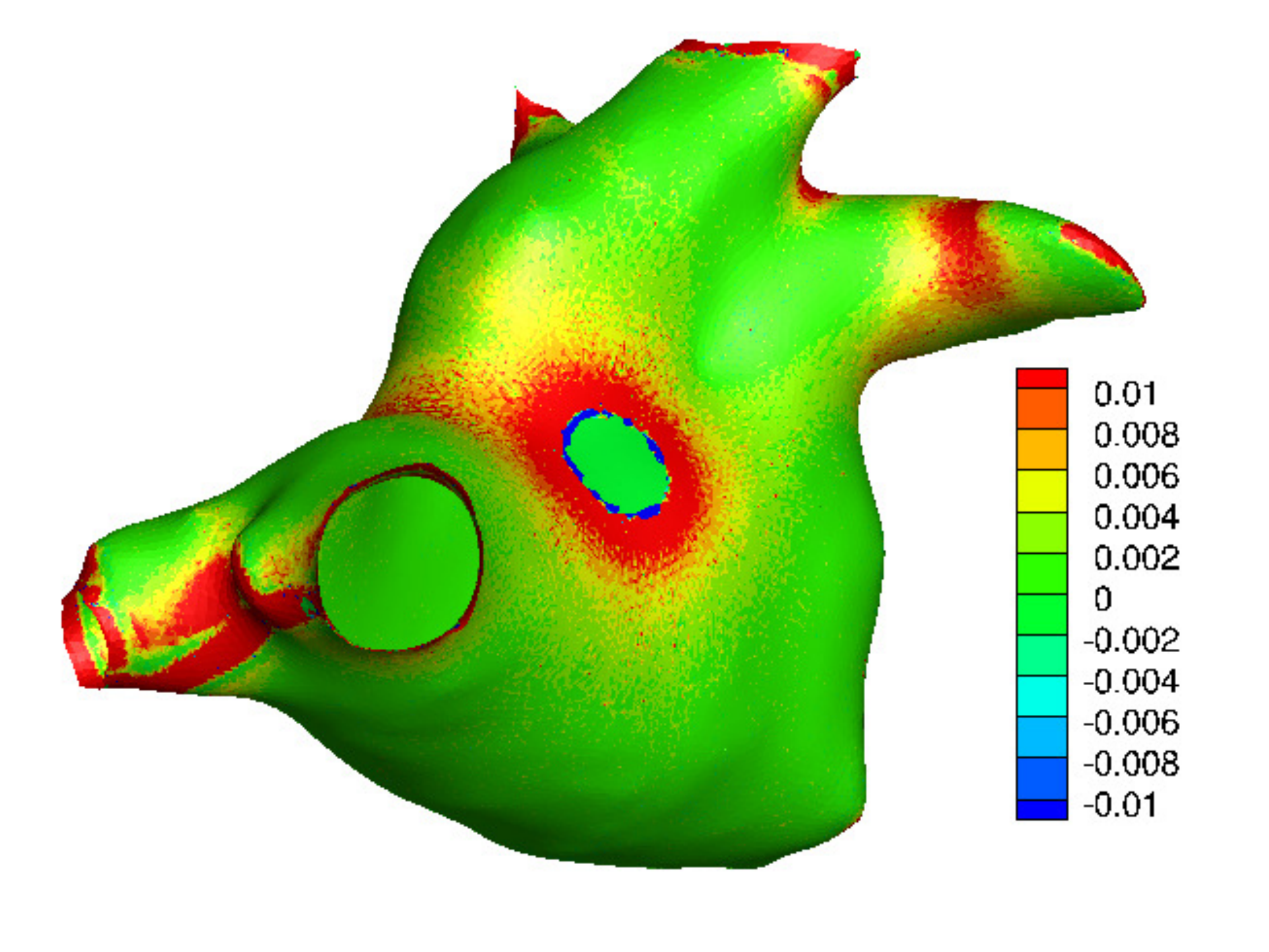} }
\caption{Without the fiber: (a) Aligned moving frames along the propagational direction. Distribution of (b) the connection form $\omega_{212}$ and (c) $\overline{\mathscr{R}}^2_{121}$.}
\label {IsoMF}
\end{figure}

\begin{figure}[ht]
\centering
\subfloat[ Front view ]  {\label{APIsoMFPVex}  \includegraphics[
width=4.5cm]{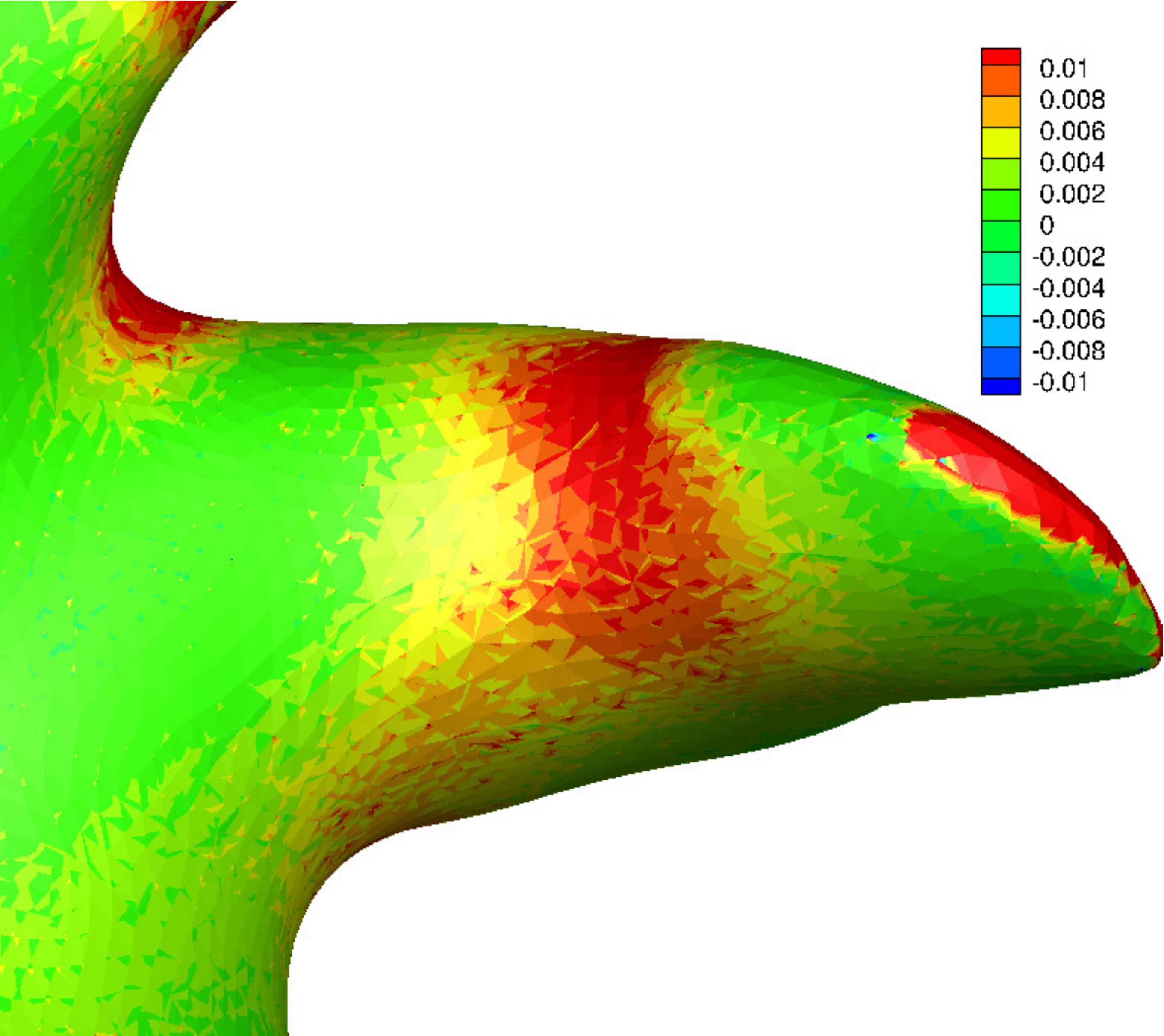} }
\subfloat[ Back view ] {\label{APIsoMFPVey}  \includegraphics[
width=4.5cm]{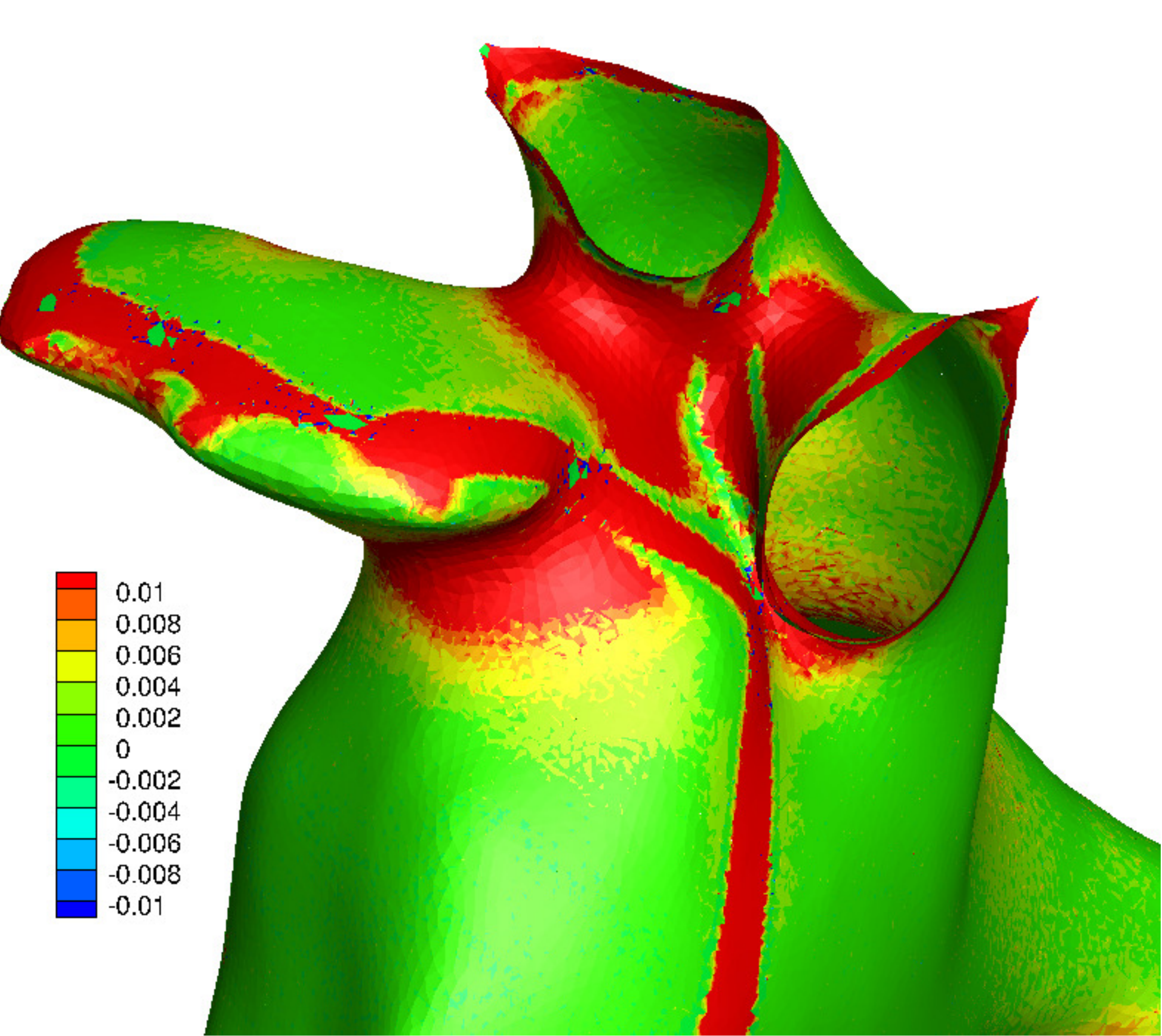} } 
\caption{Magnification of the region of the veins where $\overline{\mathscr{R}}^2_{121}$ is relatively high. A corresponding high curvature occurs due to the indented shape of the vein. }
\label {IsoMFPV}
\end{figure}

The first simulation is to run the CRP model on the atrium, but without the fiber to reveal the role of the cardiac fiber in the cardiac electric signal propagation. Fig. \ref{APIsoMF} displays the aligned moving frames along the propagation direction in the isotropic cardiac tissue. The heart without the cardiac fiber would produce a significantly decreased mechanical efficiency in pumping the oxygenated blood without twisting and squeezing. Electrophysiologically speaking, Fig. \ref{APIsoCN} and \ref{APIsoRM} reveal the presence of strong curvature (1) around the initiation point and (2) in the middle of the veins. Strong $w_{212}$ and consequently strong $\overline{\mathscr{R}}^2_{121}$ around the initiation point requires a larger magnitude or larger radius of initial excitation to overcome the initial strong curvature restraint. Moreover, the strong $w_{212}$ in the middle of the veins could be closely linked to atrial reentry because a small change in the direction to the veins can stop or allow the propagation into the vein. The magnified view in Fig. \ref{IsoMFPV} confirms a relatively strong distribution of $\overline{\mathscr{R}}^2_{121}$ around the veins (red color) in isotropic atrium.

\begin{figure}[ht]
\centering
\subfloat[ Aligned moving frames] {\label{APAniMF} \includegraphics[
width=4.5cm]{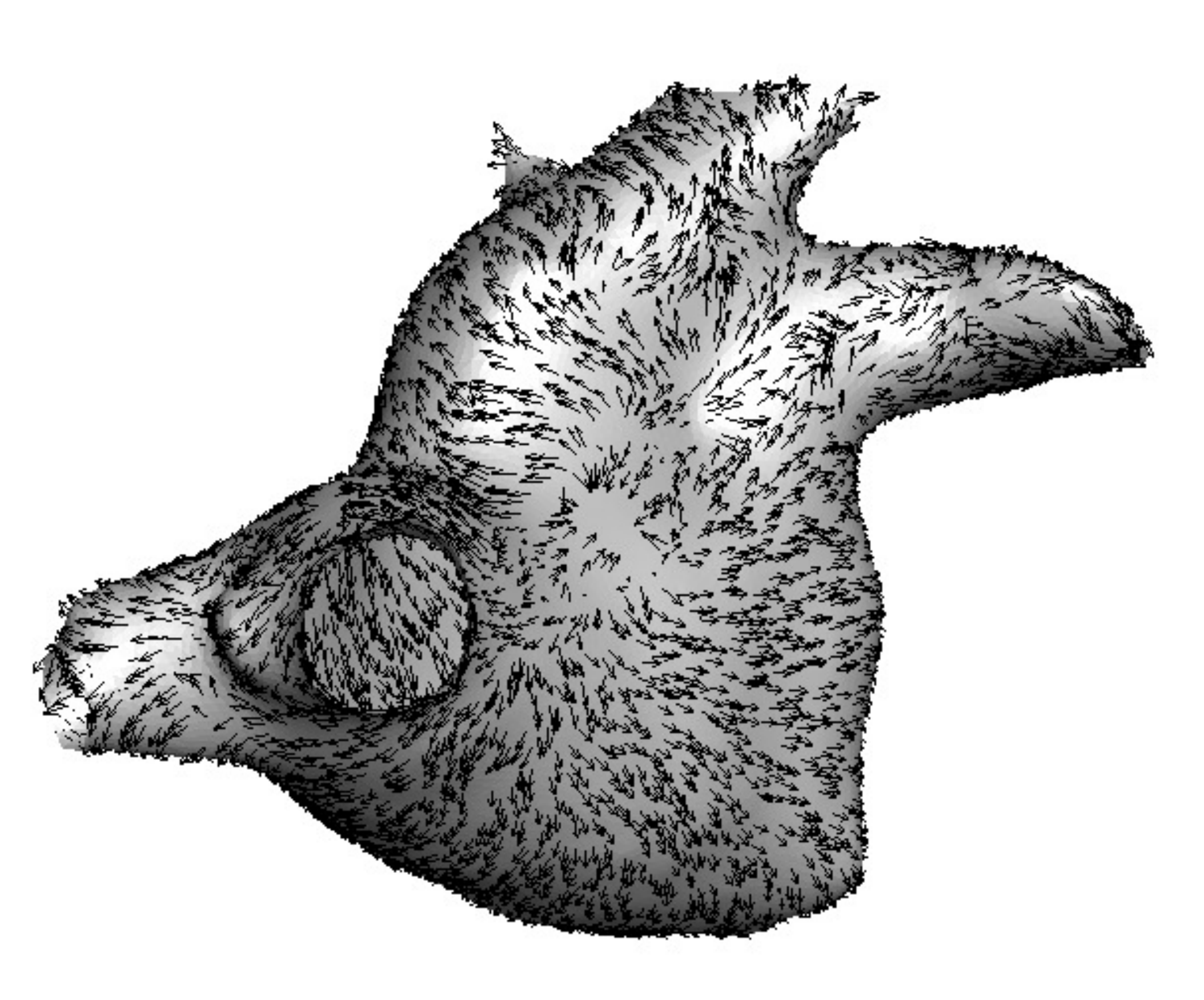} }
\subfloat[$\omega_{212}$ ] {\label{APAniCN} \includegraphics[
width=4.5cm]{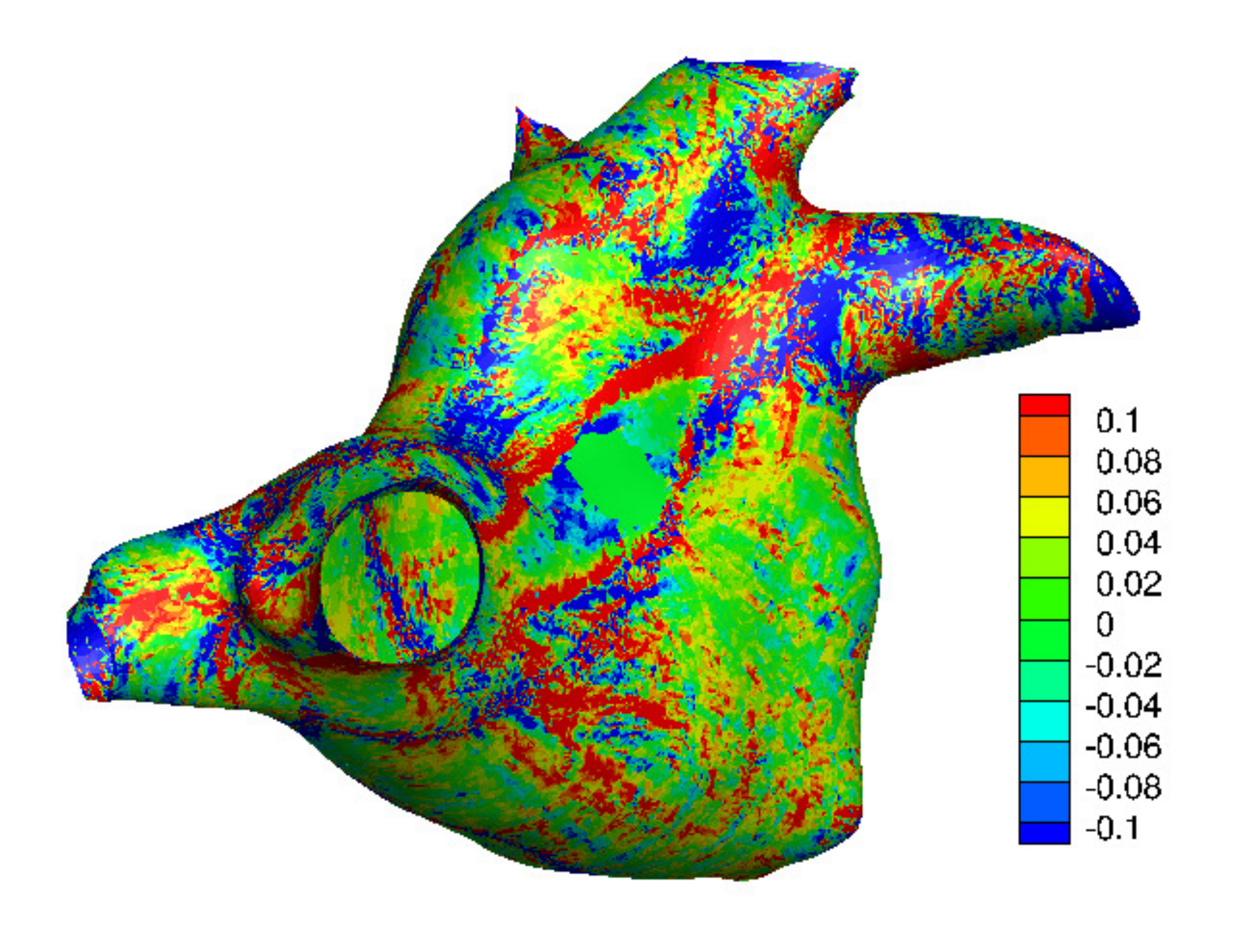} }
\subfloat[$\overline{\mathscr{R}}^2_{121}$ ] {\label{APAniCN} \includegraphics[
width=4.5cm]{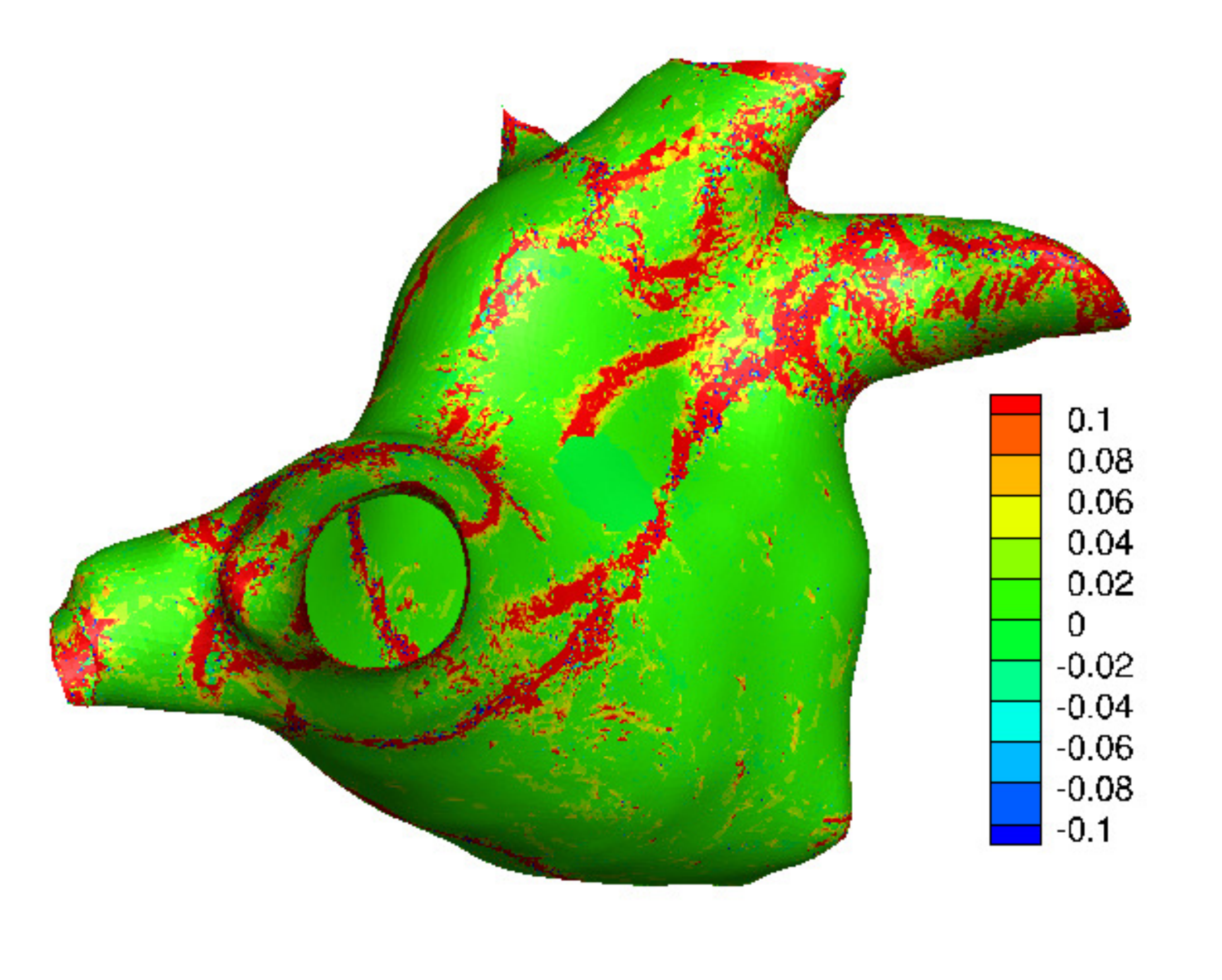} }
\caption{With the fiber: (a) Aligned moving frames along the propagational direction. Distribution of (b) the connection form $\omega_{212}$ and (c) $\overline{\mathscr{R}}^2_{121}$.}
\label {AniMF}
\end{figure}

\begin{figure}[ht]
\centering
\subfloat[  ]  {\label{AniMFPV1}  \includegraphics[
width= 4.5 cm]{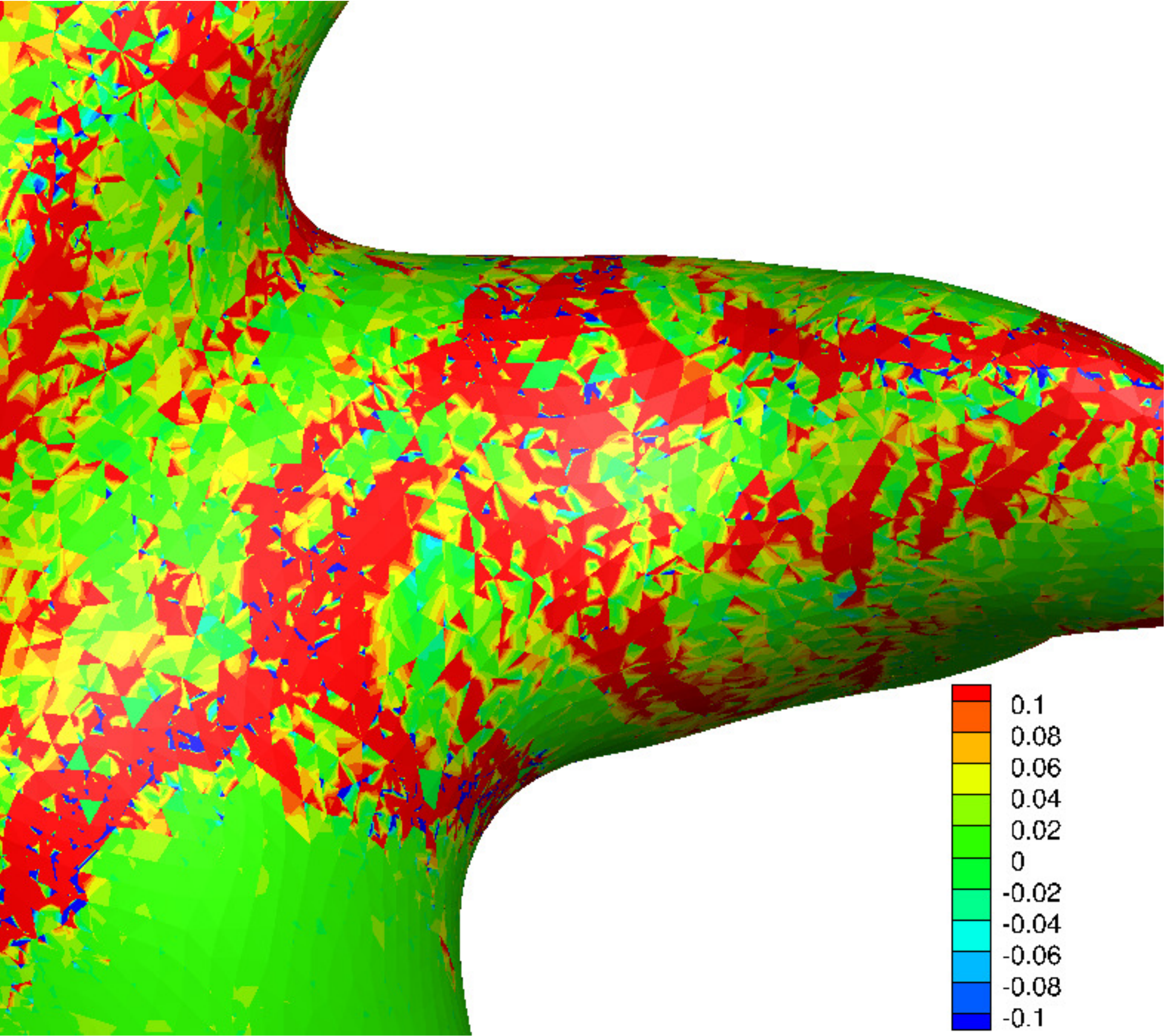} }
\subfloat[  ]  {\label{AniMFPV2}  \includegraphics[
width= 4.5 cm]{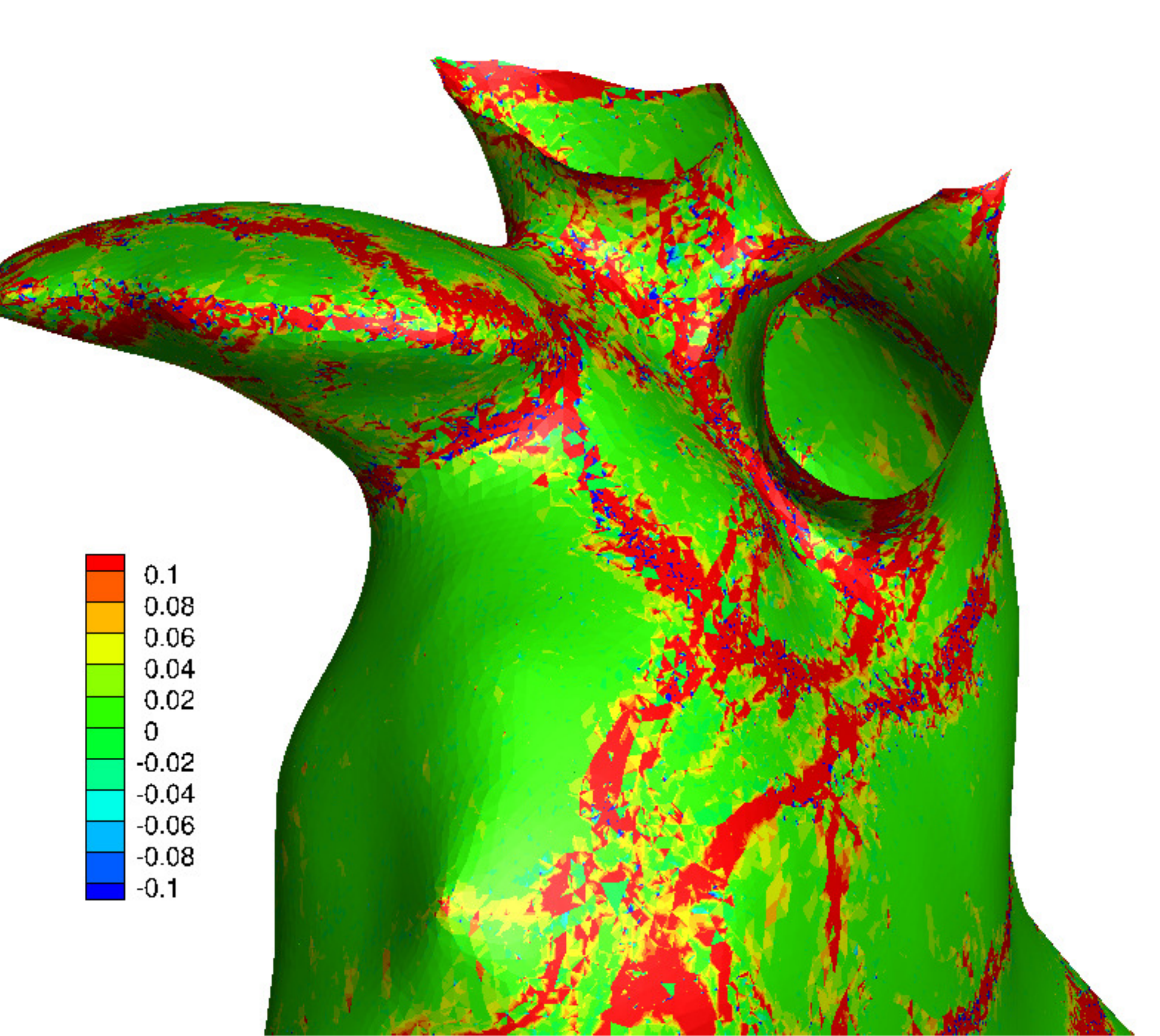} }
\caption{Magnification of the region of the veins where $\overline{\mathscr{R}}^2_{121}$ is relatively high.}
\label {AniMFPV}
\end{figure}

Adding the anisotropic cardiac fiber into the mesh of the atrium dramatically changes the trajectory of the cardiac electric propagation. Fig. \ref{APAniMF} displays the propagational direction that is approximately the same as the fiber map of Fig. \ref{Projfiber}. A quantitative comparison between the propagational direction and the fiber will be shown later. First, the strong curvature around the initiation point disappears due to anisotropy. The initial excitation by SAN could be less intensive than the isotropic medium. Second, there is a strong $\overline{\mathscr{R}}^2_{121}$ in a larger area in the veins, including the roots. This broader and stronger $\overline{\mathscr{R}}^2_{121}$ at the entrance and the veins prevent the cardiac electric signal from propagating into the veins to prevent atrial reentry and subsequent fibrillation. This new high-curvature distribution can significantly reduce the risk of generating a unidirectional block for atrial reentry by avoiding angle-dependent pathways of the veins. Fig. \ref{AniMFPV} shows the stronger Riemann curvature $\overline{\mathscr{R}}^2_{121}$ throughout the veins and the roots of the veins. The presence of the cardiac fiber dramatically changes the curvature, especially around the veins. However, the critical roles in the prevention of atrial reentry should be studied more closely in the future.

\subsection{Cardiac fiber validation}

This subsection briefly demonstrates how cardiac fiber data can be validated by the proposed scheme. The cardiac fiber is optimally designed to deliver the electric signal along the fiber direction efficiently. If the cardiac fiber is not along the propagational direction, then many factors are to be considered. A few of those have been mentioned below. 

\begin{enumerate}
\item The data of the obtained cardiac fiber contain noise.
\item The electric propagation models or their parameters for the atrium are not accurate.
\item The atrium of the patient is in an abnormal condition. 
\end{enumerate}

\begin{figure}[ht]
\centering
\subfloat[$d^{11} = 2.0$ ] {\label{fibermatchd2}  \includegraphics[
width=4.5 cm]{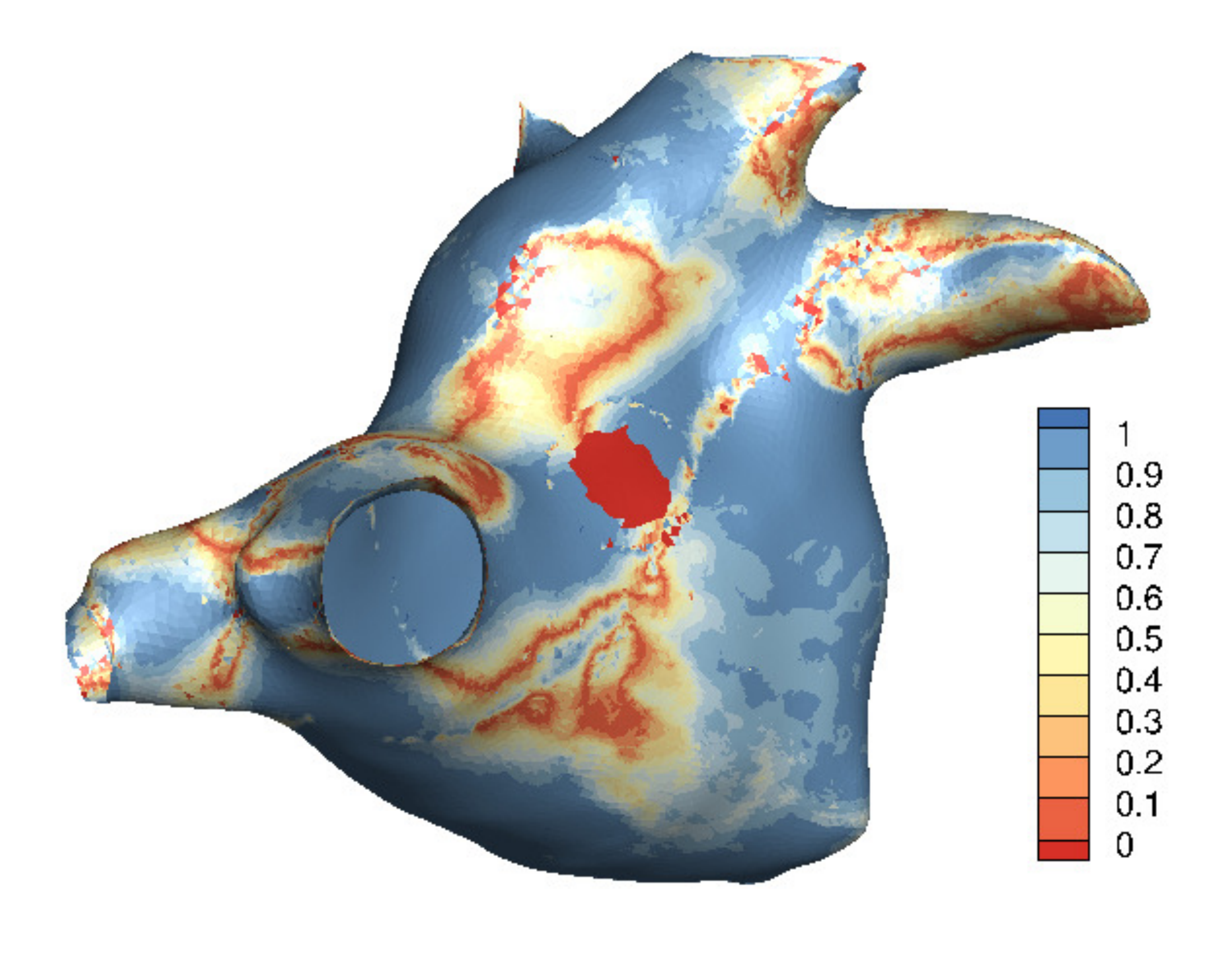} } 
\subfloat[$d^{11} = 3.0$  ]  {\label{fibermatchd3}  \includegraphics[
width= 4.5 cm]{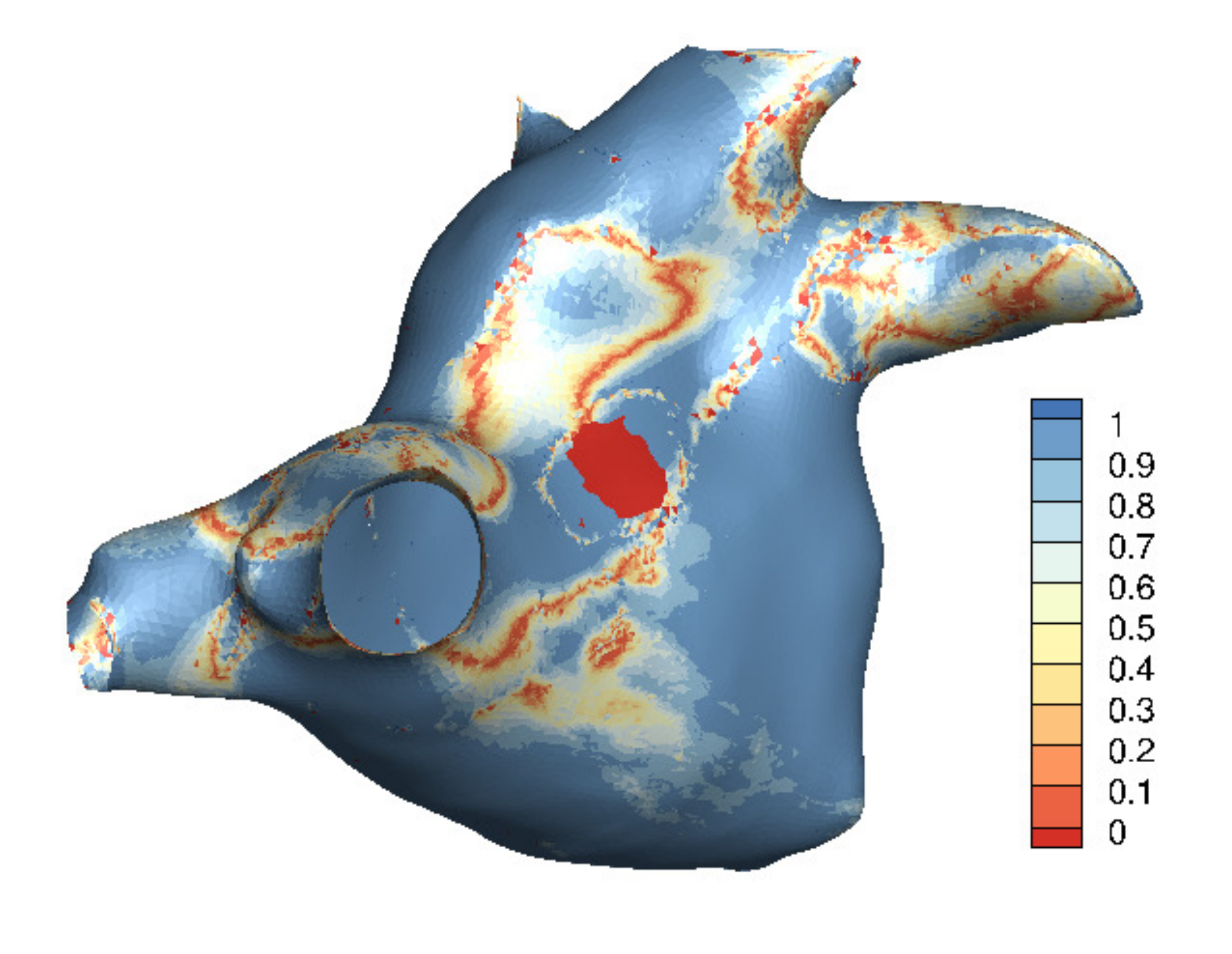} }
\subfloat[$d^{11} = 4.0$  ]  {\label{fibermatchd4}  \includegraphics[
width= 4.5 cm]{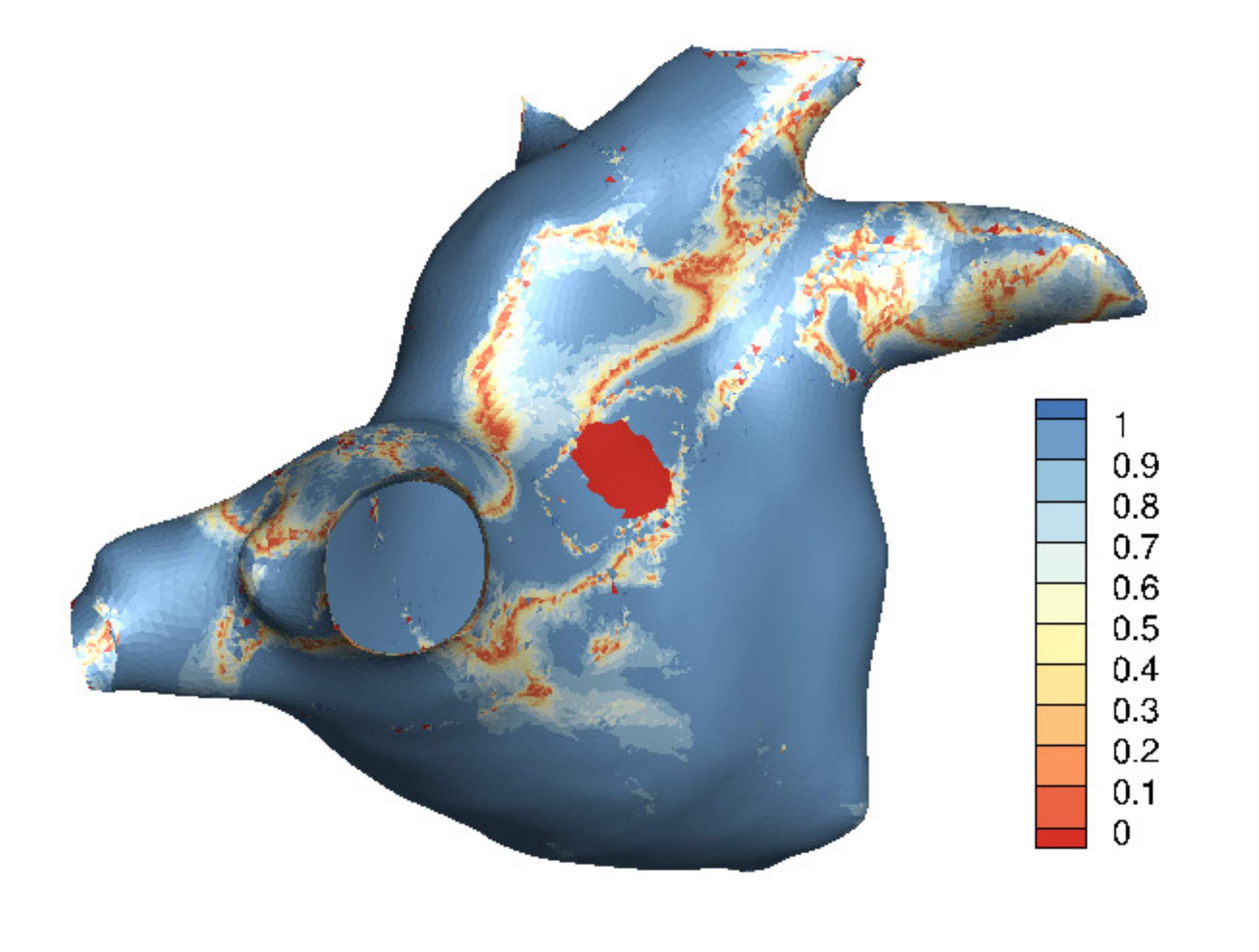} }
\caption{ Distribution of $| \cos \theta |$ where $\theta$ is the angle between the cardiac fiber and the propagational direction. The contour value of $1$ (blue) means parallel (agreement), $0$ (red) means orthogonal (disagreement).}
\label {fibermatch}
\end{figure}

Fig. \ref{fibermatch} presents the match between the cardiac fiber and the propagational direction in the atrium to display the distribution of $| \cos \theta |$ where $\theta$ is the angle between the cardiac fiber and the propagational direction. This figure presents the sufficient coincidence of the two directions except for some regions, especially around the veins. As the anisotropic strength increases as $2.0$ (Fig. \ref{fibermatchd2}) , $3.0$ (Fig. \ref{fibermatchd3}), and $4.0$ (Fig. \ref{fibermatchd4}), the cardiac electric propagation follows the cardiac fiber more closely. This computational scheme can be used to analyze the personalized atrium fiber and shape for clinical studies and prevention of atrial fibrillation.

\subsection{Effect of fibrosis}
\begin{figure}[ht]
\centering
 \includegraphics[width=6 cm]{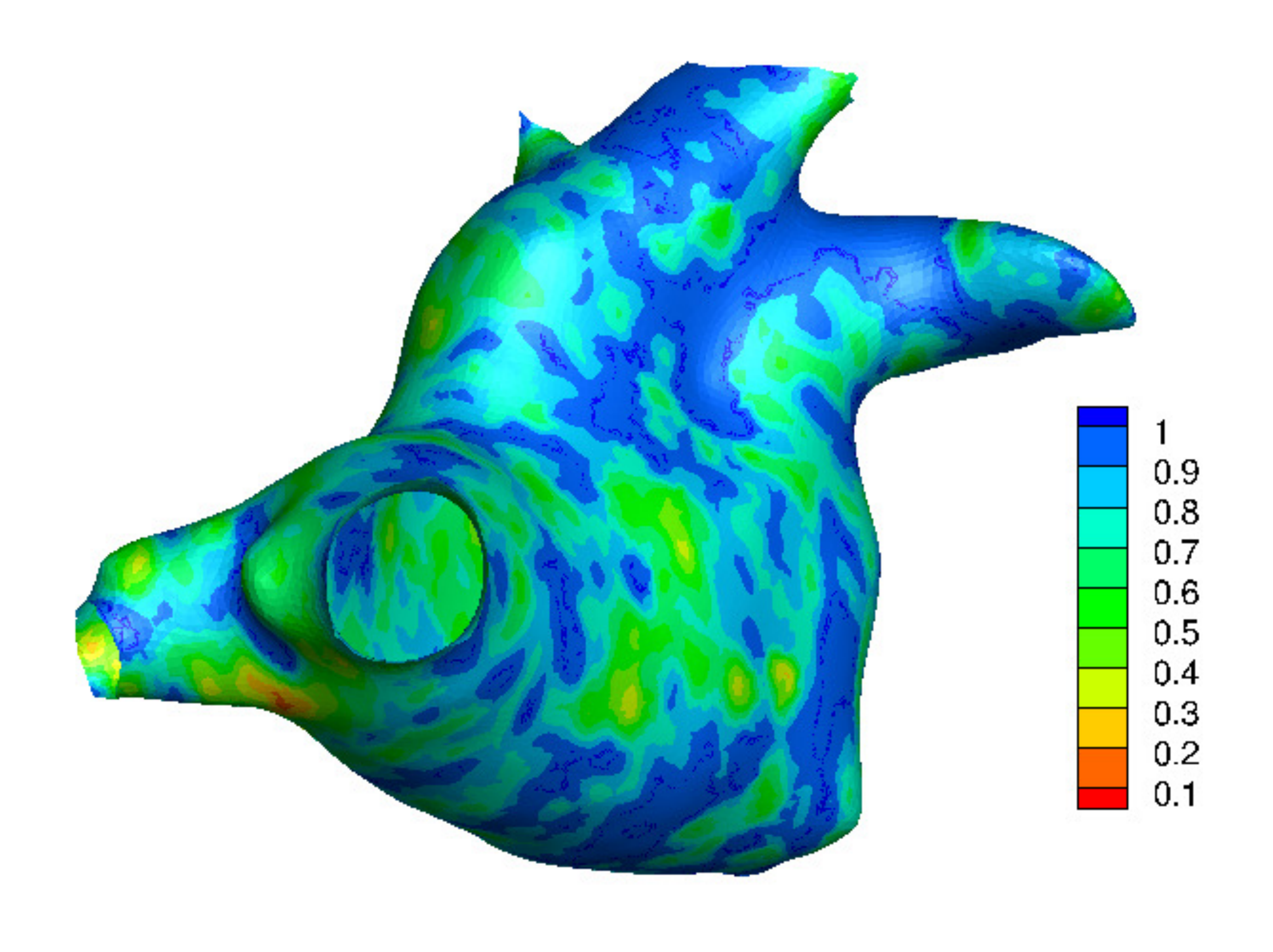}
\caption{Distribution of the conductivity of cardiac tissue. The scale of $1$ corresponds to the normal conductivity and $0$ corresponds to no conductivity. }
\label {Scartissue}
\end{figure}

Another important application of this computational scheme is the analysis of cardiac electrical signal propagation in the atrium with fibrosis, or even scar tissue. From the point of views of electrophysiology, a fibrosis corresponds to low or even zero conductivity in the cardiac tissue. Fig. \ref{Scartissue} presents an artificially generated map of the various conductivities from $0$ to $1$. A scale $1$ implies normal conductivity and $0$ means no conductivity. We use this conductivity map to simulate cardiac electrical propagations both with and without the cardiac fibers.

\begin{figure}[ht]
\centering
\subfloat[Moving frames ] {\label{IsoCmapMF}  \includegraphics[
width=4.5 cm]{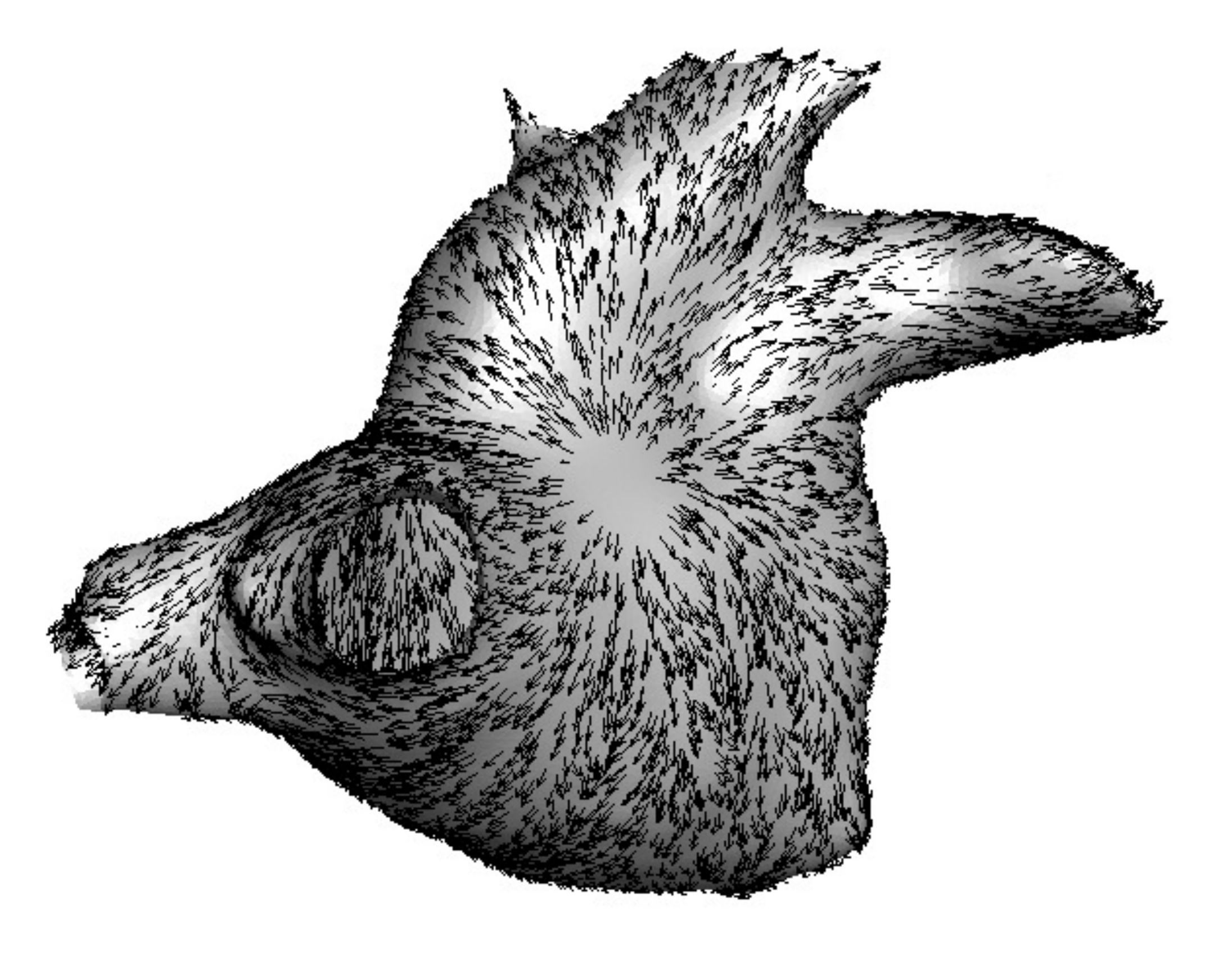} } 
\subfloat[$\omega_{212}$ ] {\label{IsoCmapCN}  \includegraphics[
width=4.5 cm]{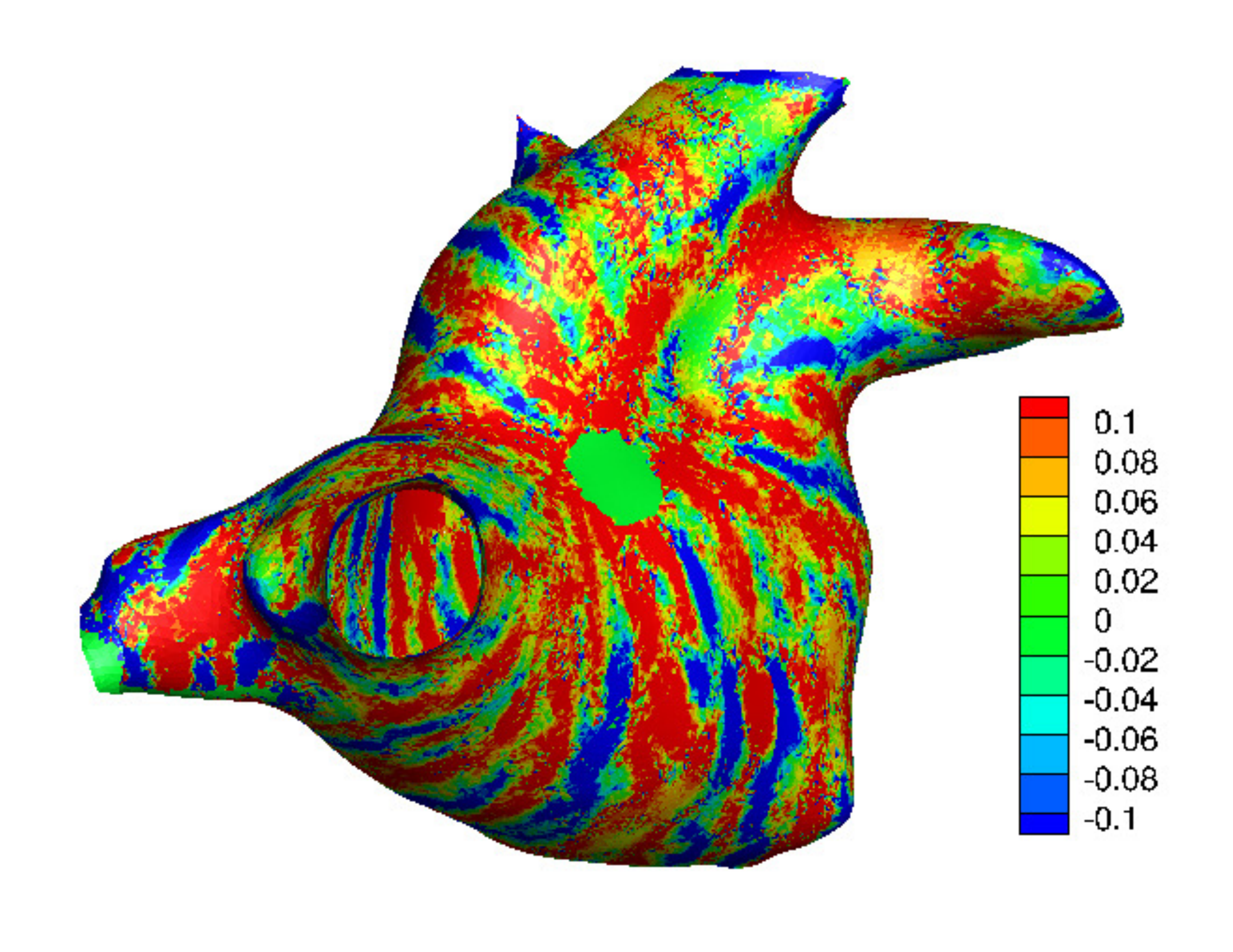} } 
\subfloat[$\overline{\mathscr{R}}^2_{121}$ ] {\label{IsoCmapCR}  \includegraphics[
width=4.5 cm]{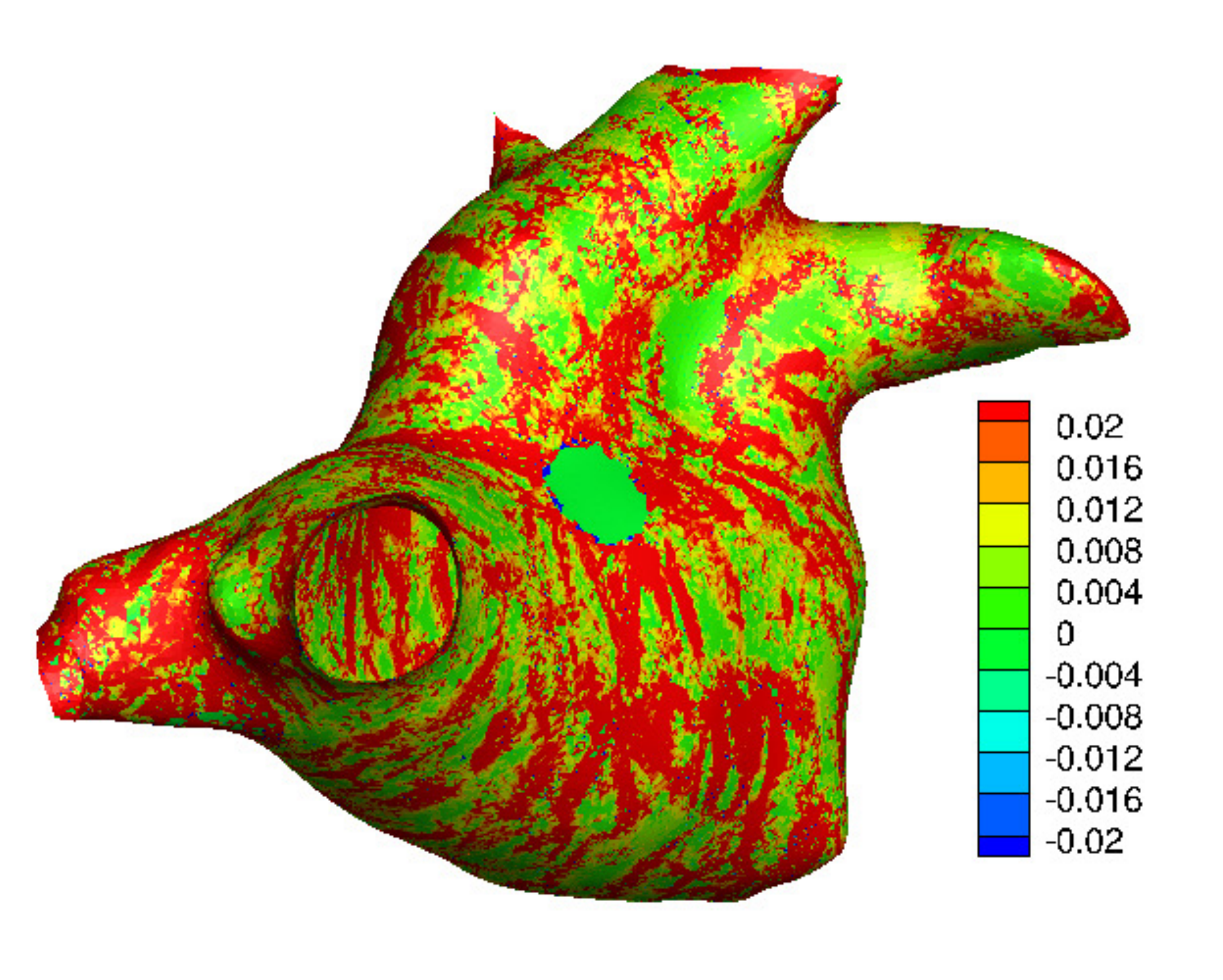} } 
\caption{Atrial Isotropy: Propagation in the presence of scar tissue: (a) aligned moving frames (b) connection component $\omega_{212}$, and (c) $\overline{\mathscr{R}}^2_{121}$. }
\label {IsoCmap}
\end{figure}

\begin{figure}[ht]
\centering
\subfloat[Moving frames ] {\label{AniCmapMF}  \includegraphics[
width=4.5cm]{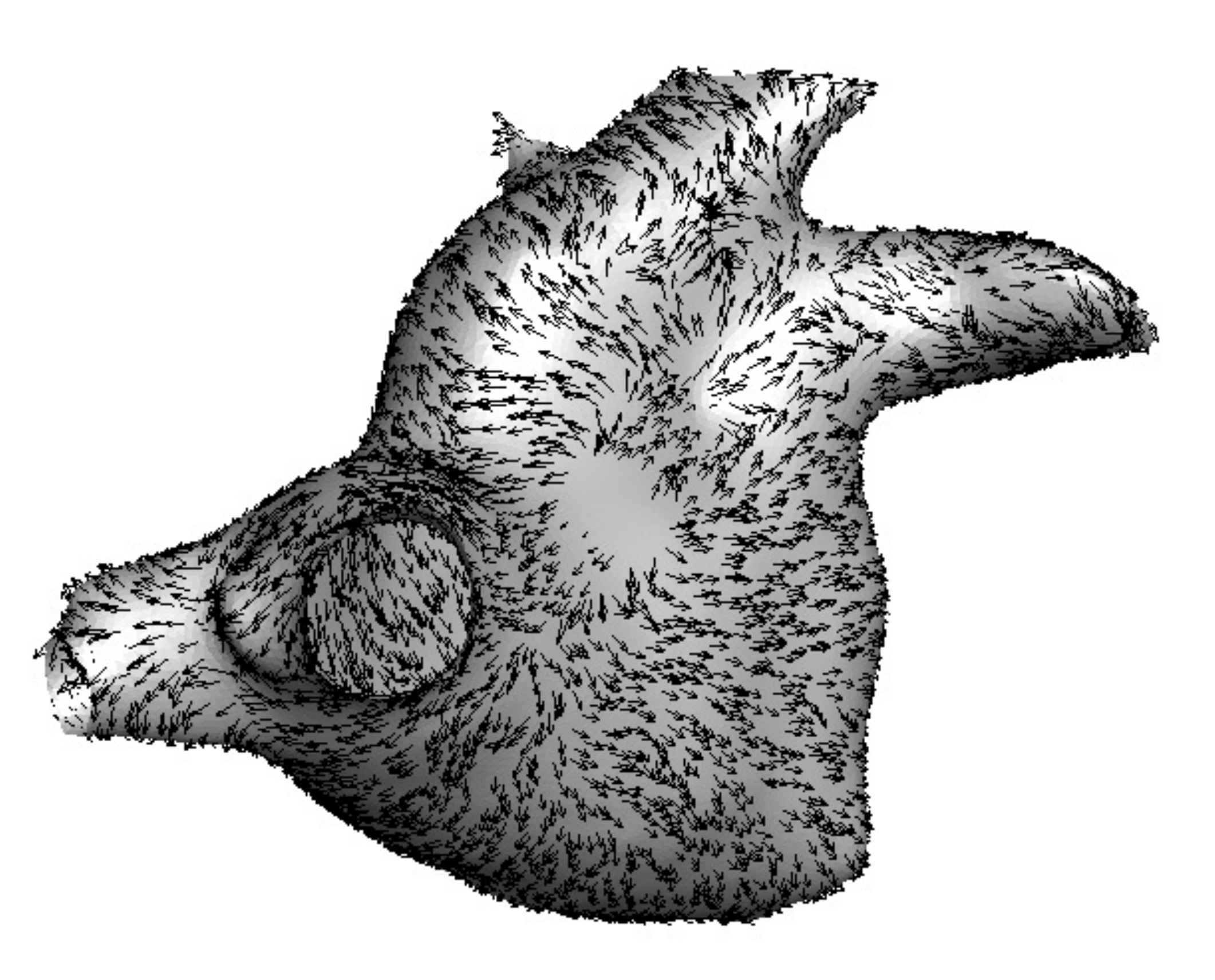} } 
\subfloat[$\omega_{212}$ ] {\label{AniCmapCN}  \includegraphics[
width=4.5 cm]{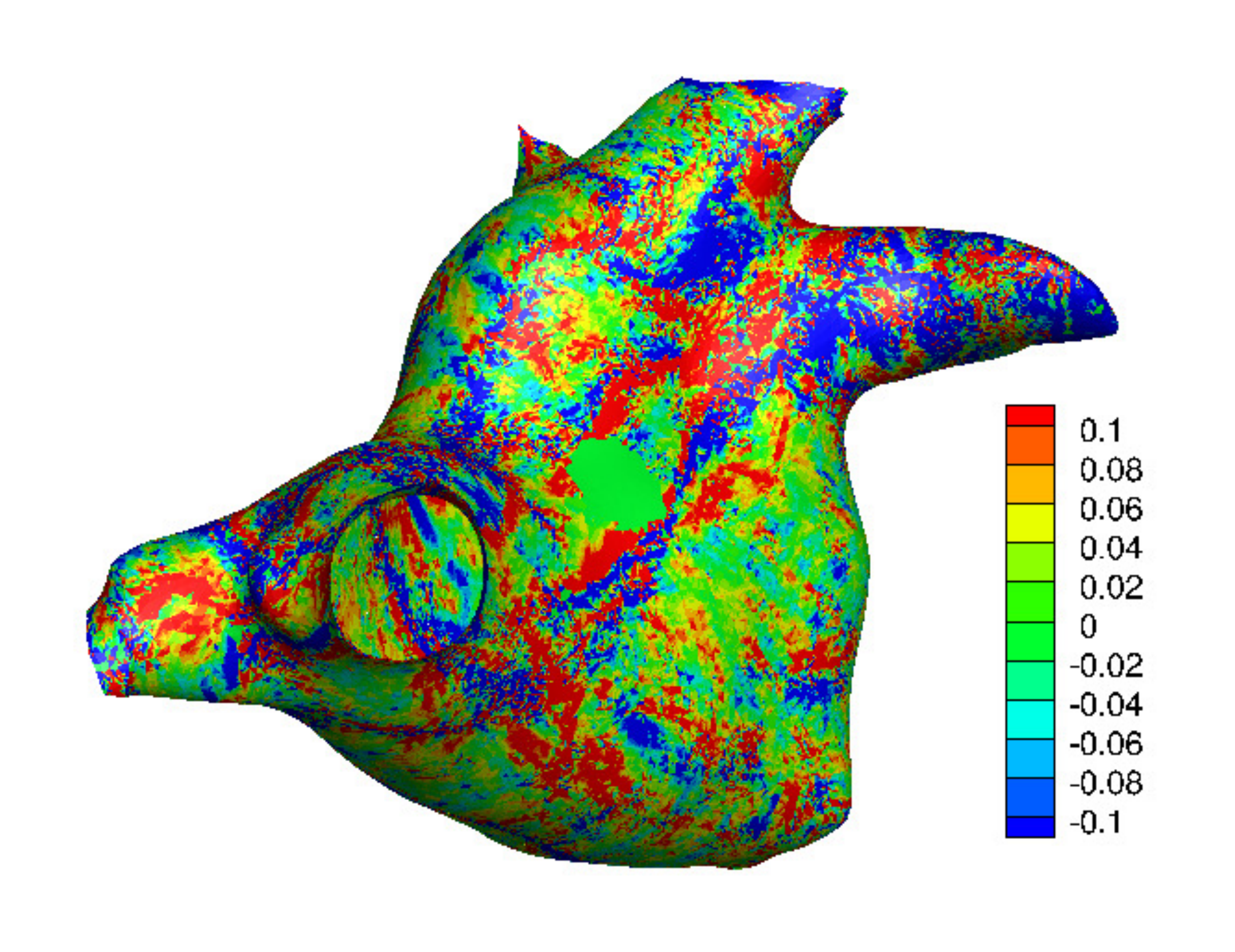} } 
\subfloat[$\overline{\mathscr{R}}^2_{121}$ ] {\label{AniCmapCR}  \includegraphics[
width=4.5 cm]{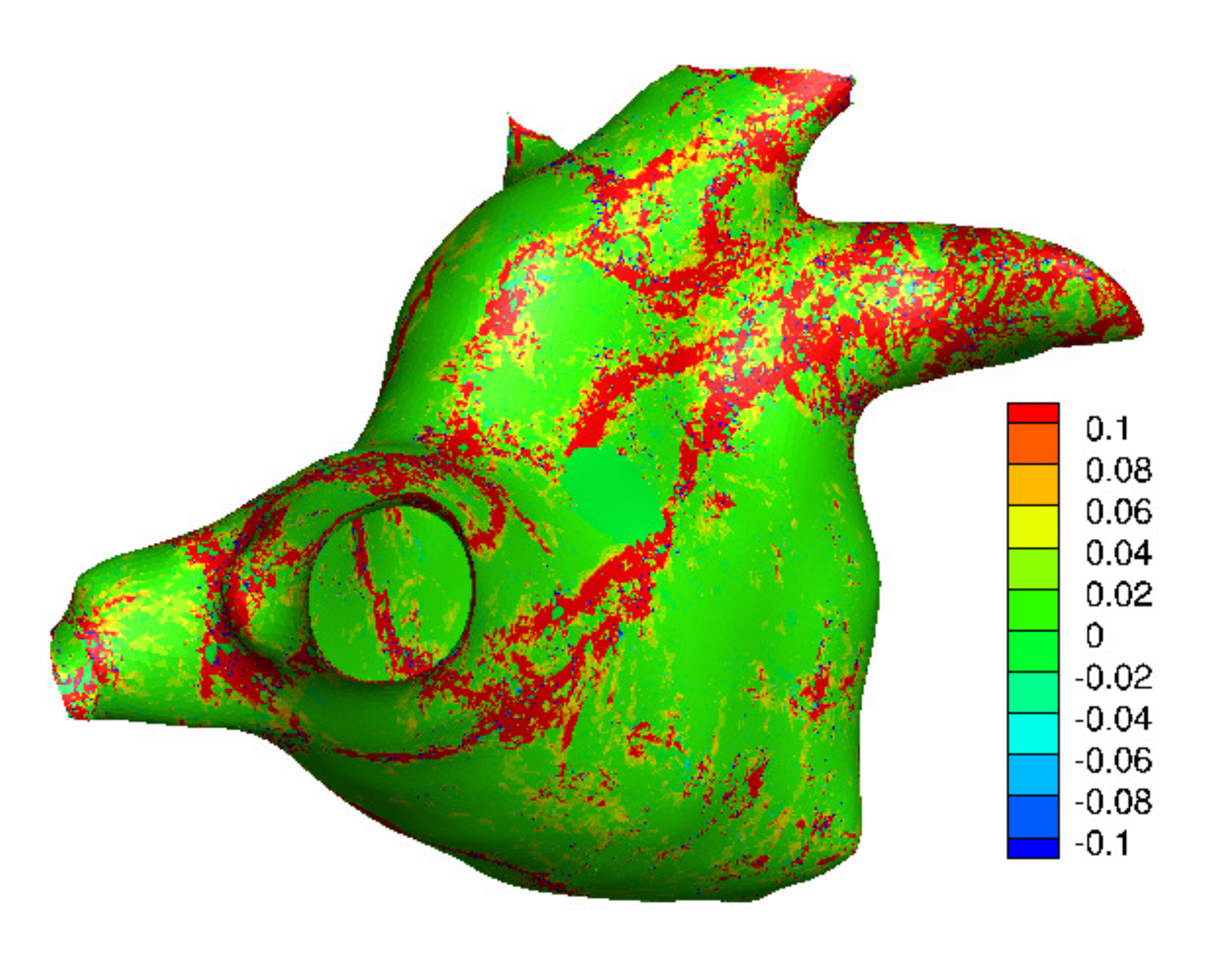} } 
\caption{Atrial anisotropy: Propagation in the presence of scar tissue: (a) aligned moving frames (b) connection component $\omega_{212}$, and (c) $\overline{\mathscr{R}}^2_{121}$. }
\label {AniCmap}
\end{figure}

Fig. \ref{IsoCmap} displays the aligned moving frame along the propagational direction and its corresponding $\omega_{212}$. In comparison with Fig. \ref{IsoMF}, where the conductivity is all $1$, we observe a notable change in the propagational direction as well as the corresponding $\omega_{212}$ and $\overline{\mathscr{R}}^2_{121}$. This implies that the propagation in an isotropic atrium is very sensitive to fibrosis, thereby yielding a considerable change in the propagation and, consequently, its excitation sequence for mechanical pumping. On the other hand, the same fibrosis with the cardiac fiber does not experience a much different propagational direction or the corresponding $\omega_{212}$ and $\overline{\mathscr{R}}^2_{121}$. This is clearly observed in Fig. \ref{AniCmap}, which is only slightly different from Fig. \ref{AniMF}.

\section{Discussion}
The connection map provides an efficient quantification of the stopping conditions of the cardiac electrical signal propagation. This map includes both spatial and temporal components when the propagation is steady. Thus, all the features of the time-dependent propagation can be viewed in a single plot. The aligned moving frames can be obtained from computational simulations, equivalent mapping, or clinical mapping data. However, the analysis of this map could be more multidimensional and contain geometric information regarding propagation. Therefore, we expect that this analysis would complement the current analysis techniques of cardiac electrical signal propagation in the map format, such as the kinematic analysis (wavefront analysis) or phase map analysis.

The first drawback of this method is that the construction of the aligned moving frames is achieved by solving the two-dimensional diffusion-reaction equations for the atrium. Thus, it is computationally expensive and not available in clinic timeline. If the numerical solution of the diffusion-reaction equations could be obtained from a set of solutions such as the inertial manifold \cite{Mallet-Paret}, then the construction time could be significantly reduced. Second, certain parameters need to be fitted for every simulation, such as $\delta$ and the minimum magnitude of the gradient, which may affect the final connection map. In the future works, these parameters could be entirely eliminated or be guided for each case individually. Third, the computation of the Riemann curvature tensor is the third derivative of the aligned moving frames, which requires the sufficient smoothness of the geometry and cardiac fiber data.

The future work would therefore involve the (1) use of the anatomically-accurate atrial and ventricular data, (2) study of the representations of spiral waves in terms of the connection form and Riemann curvature tensor, and (3) finding the differences in the connection map for cardiac restitution and cardiac memory.

\section*{Acknowledgements}
The authors thank Professor Eun-Jae Park (Dept. of Computational Science and Engineering, Yonsei University) for inspirational and encouraging discussion. This research was supported by the Basic Science Research Program through the National Research Foundation of Korea (NRF) and funded by the Ministry of Education, Science and Technology (No. 2016R1D1A1A02937255).

\bibliographystyle{elsarticle-num}
\bibliography{MMFConnection_ArXiv}

\end{document}